\documentclass[12pt,leqno]{amsart}
\usepackage{amsmath}
\usepackage{amssymb}
\usepackage{amsthm} 
\usepackage{stmaryrd,mathrsfs}
\usepackage{epsf}
\usepackage{graphicx}
\usepackage{esint} 
\usepackage{dsfont}
\usepackage[pagebackref, colorlinks]{hyperref} 
\usepackage{color}
\usepackage{enumerate}

\numberwithin{equation}{section}

\newtheorem{lem}{Lemma}[section]
\newtheorem{prop}[lem]{Proposition}
\newtheorem{defi}[lem]{Definition}
\newtheorem{thm}[lem]{Theorem}
\newtheorem{cor}[lem]{Corollary}
\theoremstyle{remark}
\newtheorem{rmk}[lem]{Remark}

\newcommand{\nn}{\nonumber}
\newcommand{\ms}{\medskip}

\newcommand{\R}{\mathbb{R}}
\newcommand{\C}{\mathbb{C}}
\newcommand{\N}{\mathbb{N}}

\renewcommand{\Re}{\mathscr{R}\hspace{-1.5pt}e}
\renewcommand{\Im}{\mathcal{I}\hspace{-1.5pt}m}

\newcommand{\dr}{\partial}
\newcommand{\dist}{\,\mathrm{dist}}
\newcommand{\sm}{\setminus}
\newcommand{\supp}{\mathrm{supp}}
\newcommand{\diam}{\mathrm{diam}}

\newcommand{\wt}{\widetilde}
\newcommand{\wh}{\widehat}

\newcommand{\cL}{{\mathcal  L}}
\newcommand{\cD}{{\mathcal  D}}

\newcommand{\cJ}{{\mathcal  J}}

\newcommand{\cA}{{\mathcal  A}}
\newcommand{\cB}{{\mathcal  B}}
\newcommand{\cC}{{\mathcal  C}}

\newcommand{\W}{W}

\newcommand{\vW}{v_\W}
\newcommand{\cH}{\mathcal H}

\newcommand{\1}{{\mathds 1}}

\DeclareMathOperator{\Tr}{Tr}
\DeclareMathOperator*{\essinf}{ess\, inf}

\newcommand{\ds}{\displaystyle}

\theoremstyle{plain}

\newcommand{\bp}{\noindent {\it Proof}.\,\,}
\newcommand{\ep}{\hfill$\Box$ \vskip 0.08in}

\def\div{\mathop{\operatorname{div}}}

\definecolor{violet}{cmyk}{0.91,0,0.88,0.12}

\begin{document}

\title{Dirichlet problem in domains with lower dimensional boundaries}

\author{J. Feneuil}
\address{Temple University, Philadelphia PA, USA}
\email{joseph.feneuil@temple.edu}

\author{S. Mayboroda}
\address{University of Minnesota, Minneapolis MN, USA}
\email{svitlana@math.umn.edu}

\author{Z. Zhao}
\address{Institute for Advanced Study, Princeton NJ, USA}
\email{zzhao@ias.edu}

\thanks{The second author was  supported  by the Alfred P. Sloan Fellowship, the NSF INSPIRE Award DMS 1344235, NSF CAREER Award DMS 1220089, the NSF RAISE-TAQ grant DMS 1839077, and the Simons Foundation grant 563916, SM. The third author was partially supported by NSF grant numbers DMS-1361823, DMS-1500098 and DMS-1664867.}

\date{October 25, 2018}

\subjclass[2010]{35J25, 35J70}

\keywords{Degenerate elliptic operators, Dirichlet problem, complex coefficients, Moser-type estimates, Non-tangential maximal function, Square functionals, Carleson measures.}

\maketitle

\begin{abstract} The present paper pioneers the study of the Dirichlet problem with $L^q$ boundary data for second order operators with complex coefficients in domains with lower dimensional boundaries, e.g.,  in $\Omega := \R^n \setminus \R^d$ with $d<n-1$. Following the first results in \cite{DFMprelim}, \cite{DFMAinfty}, the authors introduce an appropriate degenerate elliptic operator and show that the Dirichlet problem is solvable for all $q>1$ provided that the coefficients satisfy the small Carleson norm condition.  

Even in the context of the classical case $d=n-1$, (the analogues of) our results are new. The conditions on the coefficients are more relaxed than the previously known ones (most notably, we do not impose any restrictions whatsoever on the first $n-1$ rows of the matrix of coefficients) and the results are more general. We establish local rather than global estimates between the square function and the non-tangential maximal function and, perhaps even more importantly, we establish new Moser-type estimates at the boundary and improve the interior ones. 
\end{abstract}

\tableofcontents

\section{Introduction}

\subsection{State of the art}

The main objective of the present article is the study of the well-posedness  of the Dirichlet problem. Formally, given an open domain $\Omega \subset \R^n$, an elliptic operator $L$ in divergence form, and $q\in (1,+\infty)$, we say that the Dirichlet problem (D$_q$) is well-posed if for any $g\in L^p (\dr \Omega, \sigma)$, we can find a unique function $u$ such that 
\begin{equation} \label{Dp1}
 L u = 0 \text{ in } \Omega,
\end{equation}
\begin{equation} \label{Dp2}
 u = g \text{ on } \partial \Omega,
\end{equation}
and 
\begin{equation} \label{Dp3}
\|N(u)\|_{L^q(\partial \Omega,\sigma)} \leq C \|g\|_{L^q(\partial \Omega,\sigma)},
\end{equation}
where $N$ is a non-tangential maximal function and the constant $C>0$ is independent of $g$. Obviously, the statement above is not complete as we need to make precise the meaning of \eqref{Dp1}--\eqref{Dp3}. We naturally expect \eqref{Dp1} to be taken in the weak sense. 
But what is the meaning of \eqref{Dp2} when $u$ is defined on the open set $\Omega$? If $\Omega$ is very irregular, a proper choice of the measure $\sigma$ on $\dr\Omega$ and the definition of $N(u)$ is unclear as well. 
What \eqref{Dp1}--\eqref{Dp3} means in our context will be carefully explained later, but first, let us give a brief (and somewhat narrowly focused) presentation of the relevant history of the Dirichlet problem. 
Due to the huge literature on the topic, we will be unable to cite all the works, and we apologize in advance for the omissions.

We start the history on the topic with a work of Dahlberg (see \cite{Da}). Let $\Omega$ be a bounded Lipchitz domain. Dahlberg proved that the harmonic measure (for the Laplace operator) defined on the boundary $\dr \Omega$ is $A^\infty$-absolutely continuous\footnote{ $A^\infty$ absolute continuity is a quantitative version of the mutual absolute continuity} with respect to the surface measure on $\dr \Omega$. 
This property is known to imply that the Dirichlet problem (D$_q$)  - associated with the domain $\Omega$ and the elliptic operator $\Delta$ - is well-posed for large enough $q$. It was proved just a little later that D$_q$ is well-posed for the Laplacian on Lipschitz domains for any $2-\varepsilon<q<\infty$
and that the range is sharp in the sense that for any $q<2$, we can find a Lipschitz domain such that (D$_q$) is false \cite{Da2, KenigB}. 

Consider the Laplacian in the domain $\Omega$ lying above the graph of a Lipschitz function $\varphi : x\in \R^{n-1} \to \R$. Let us keep it mind that a bi-Lipschitz change of variable $\rho: \Omega \to \rho(\Omega)$ preserves the well-posedness of the Dirichlet problem. 
In \cite{JK} Jerison and Kenig use the change of variable 
$$\rho : \, (x,t) \in \R^n  \to (x,t-\varphi(x))$$
 that flattens the domain $\Omega$ and that maps the Laplacian to another elliptic operator $\cL = - \div A_0 \nabla$ with bounded, measurable, symmetric, and $t$-independent coefficient. By using a Rellich identity, they establish that these conditions on $\cL$ are sufficient to ensure that (D$_2$) is well-posed, and hence, they extend the result of Dahlberg to the case where the Laplacian is replaced by an elliptic operator $L= - \div A \nabla$ where $A$ has real, bounded, symmetric, and $t$-independent coefficients. Analogous results for real non-symmetric operators have been proved much later in \cite{KKPT}, \cite{HKMP}  using square function/non-tangential maximal function estimates and elements of the solution of the Kato problem. The complete situation for operators with complex coefficients is still not clear, although the solution of the Kato problem \cite{Kato} and later developments allowed to treat block matrices and some of their generalizations \cite{AMM}, \cite{M}. 

The next breakthrough we shall talk about in the study of the Dirichlet problem is a result from Kenig and Pipher (see \cite{KePiDrift}, which use methods developed in \cite{KKPT}). Consider again the Laplacian and $\Omega = \{t > \varphi(x)\}$ a domain that lies above the graph of the function $\varphi : x\in \R^{n-1} \to \R$. 
The change of variable
$$\rho: \, (x,t) \in \R^n \to (x, ct-\varphi_t(x)),$$
where $c$ is a large constant and $\varphi_t$ is the convolution of $\varphi$ by a smooth mollifier,  also sends $\Omega$ to $\R^{n}_+$, but maps now the Laplacian $-\Delta$ to an elliptic operator $\cL := -\div A_0 \nabla$ where $A_0$ satisfies the conditions that $|\nabla A_0(x,t)| \leq C/t$ and $|t\nabla A_0(x,t)|^2 dx \, \frac{dt}{t}$ is a Carleson measure. Kenig and Pipher showed that the two latter conditions are enough to ensure that (D$_q$) is well-posed if $q$ is large enough, hence extending the result of Dahlberg to a new class of elliptic operators.

Dindo\v s, Petermichl, and Pipher studied in \cite{DPP} the conditions needed for well-posedness (D$_q$) when $q>1$ is small. They established that, for a given $q>1$, the Dirichlet problem (D$_q$) - associated to the Lipschitz domain $\Omega$ and the elliptic operator $L=-\div A \nabla$ - is well-posed if both the Lipschitz constant of $\Omega$ and the Carleson norm of $|t\nabla A|^2 dx\, \frac{dt}{t}$ are smaller than $\epsilon(q) \ll 1$.

One has to also mention a number of perturbation results, in $L^\infty$ and in Carleson measure norm, which we shall not review here. 

Our focus is on operators with coefficients whose gradient satisfies the Carleson measure condition, as above. All the previous results that we mentioned in this context were established in the case where $L=-\div A \nabla$ has real coefficients. In \cite{DP}  Dindo\v s and Pipher introduced a notion of $q$-ellipticity based on a notion of $L^q$-dissipativity (see \cite{CD,CM}), and cleverly used this notion of $q$-ellipticity to obtain ``$q$-Cacciopoli's inequalities'' and ``reverse H\"older inequalities'' (see Subsection \ref{ssES} for the precise statement) which can be seen as a weakened version of Moser's estimates. 
They used these partial estimates to get the well-posedness of the Dirichlet problem (D$_q$) whenever the $q$-ellipticity, in addition to appropriate Carleson measure estimates and some structural conditions on coefficients,  holds.

\medskip

In all the above works, the boundary of the domain $\Omega \subset \R^n$ has Hausdorff dimension $n-1$. In \cite{DFMprelim}, Guy David and the first two authors of the present article launched an elliptic theory adapted to domains $\Omega$ that are the complement in $\R^n$ of sets $\Gamma$ with dimensions $d < n-1$. 
Since the boundary $\Gamma$ of the domain $\Omega$ is too thin to be ``seen'' by the Laplacian (or the general elliptic operators), the operators $L:=-\div A \nabla$ used in this theory are degenerate and satisfy the ellipticity condition with a different homogeneity. Precise definition are given later, see \eqref{defEllip2}, \eqref{defBdd2}. We also mention that similar degenerate operators have been considered before, notably in \cite{FKS, FJK}, but the well-posedness or any related boundary estimates have never been attacked. 

Subsequently in \cite{DFMAinfty}, David, Feneuil, and Mayboroda established that if $\Gamma$ is the graph of a function $\varphi : \, \R^d \to \R^{n-d}$ with small Lipschitz constant, we can find a particular degenerate elliptic operator $L:=-\div A \nabla$ such that the harmonic measure on $\Gamma$ (and associated to $L$) is $A^\infty$-absolutely continuous with respect to the $d$-dimensional Hausdorff measure on $\Gamma$.  While this operator, being the simplest one that we can treat, can be thought of as an analogue of the Laplacian in the domains with lower dimensional boundaries, it already carries most of the difficulties exhibited by the operators whose coefficients satisfy the aforementioned Carleson condition. Indeed, by necessity, it is not a constant coefficient operator, and it cannot be $t$-independent either: a rather delicate dependence of the coefficients on the distance to the boundary is exactly what makes the problem well-posed in the higher co-dimensional context. 
The  proof in \cite{DFMAinfty} uses some ideas from \cite{KKPT,KKiPT,DPP2}, and relies on a {\it new} change of variable $\rho$ that sends $\Omega$ to $\R^n \setminus \R^d = \R^d \times (\R^{n-d} \setminus \{0\})$ and that is almost - up to Carleson measure - an isometry in the last $n-d$ variables. As shown in \cite{MZ}, the $A^\infty$-absolute continuity of the harmonic measure implies the well-posedness of the Dirichlet problem (D$_q$) when $q$ is large. 
\bigskip

The main aim of this article is to prove that, given $q>1$ and $d<n-1$, the Dirichlet problem (D$_q$) is well-posed in the domain $\R^n\setminus \R^d$ for any degenerate elliptic operator (in the full generality of possibly complex coefficients) with conditions in the spirit of \cite{KePiDrift}. As a consequence, whenever $\Gamma$ is the graph of a function with small Lipschitz constant and $L$ is given as in \cite{DFMAinfty}, the Dirichlet problem (D$_q$) is well-posed.

While the article is written with $d<n-1$ in mind, all our computations can be adapted with very light changes to the case where $d=n-1$ and the domain is the upper half plane $\R^{d+1}_+$. In that context it can be viewed as an alternative to \cite{DP}. We build on their ideas and introduce new tools. As a result, even in the classical setting, our conditions on the operator are weaker and our results are somewhat stronger. Most notably, we do not impose any restrictions on the first $n-1$ rows of the coefficient matrix. We will make these statements more precise below. 

\subsection{Main result} \label{SSMR}

Let $d<n-1$  and $\Omega = \R^n \sm \R^d = \{(x,t) \in \R^n, \, x\in \R^d \text{ and } t\in \R^{n-d}\sm \{0\}\}$; the boundary $\dr \Omega$ is assimilated to $\R^d$. We write $X = (x,t)$ or $Y = (y,s)$ for the running points in $\Omega$. The notation $B_l(x)$ is used for the ball in $\R^d$ with center $x$ and radius $l$.

\ms

First, we introduce the square function and the averaged version of the non-tangential maximal function. Let $a$ be a positive number. For $x\in \R^d$, we define the regular cones in $\R^{d+1}_+:= \R^d \times (0,+\infty)$ as
$$\Gamma_a(x) := \left\{(z,r) \in \R^{d+1}_+,\,  |z-x| < ar \right\},$$
and the higher co-dimensional cones as
$$ \wh \Gamma_a(x) := \{(y,s) \in \Omega, \, |y-x|<a|s|\}.$$
For $(z,r)\in \R^d \times (0,\infty)$, we write $W_a(z,r)$ for the Whitney box
$$W_a(z,r) = \{(y,s) \in \Omega, \, y\in B_{ar/2}(z),\, r/2 \leq |s| \leq 2r \}.$$
If, in addition, $q\in (1,+\infty)$,  we define the $q$-adapted square function $S_{a,q}$ as
\begin{equation}\label{defSq}
S_{a,q}(v)(x) := \left( \iint_{\wh \Gamma_a(x)} |\nabla v(y,s)|^2 |v(y,s)|^{q-2} dy \frac{ds}{|s|^{n-2}} \right)^\frac1q,
\end{equation}
where $v$ is a measurable function that satisfies $|v|^{q/2-1} v \in W^{1,2}_{loc}(\Omega)$. Note that the definition of $S_{a,q}$ makes sense even when $q<2$, because, in particular, $\nabla v \equiv 0$ almost everywhere on $|v|=0$. For any function $v\in L^q_{loc}(\Omega)$ and any $x\in \R^d$
the non-tangential maximal function (in the average sense) $\wt N_{a,q}$ is 
\begin{equation}\label{defNq}
\wt N_{a,q}(v)(x) = \sup_{\Gamma_a(x)}  v_{W,a,q}
\end{equation}
where $v_{W,a,q}$ is defined on $\R^{d+1}_+$ by
\begin{equation} \label{defw}
v_{\W,a,q}(z,r) := \left( \frac1{|W_a(z,r)|} \iint_{W_a(z,r)} |v(y,s)|^q dy \, ds   \right)^\frac1q.
\end{equation}

\medskip

We say that $\cL=-\div A \nabla$ is an elliptic operator \textit{with weight $|t|^{d+1-n}$} (we will always omit this weight for brevity) if there exists $C>0$ such that the complex matrix $A$ satisfies
\begin{equation} \label{defEllip}
 \Re( A(X)\xi \cdot \bar \xi) \geq C^{-1} |t|^{d+1-n} |\xi|^2 \qquad \text{ for } X\in \R^n \setminus \R^d, \, \xi \in \C^n,
\end{equation}
and
\begin{equation} \label{defBdd}
|A(X)\xi \cdot \bar \zeta| \leq C |t|^{d+1-n} |\xi| |\zeta| \qquad \text{ for } X\in \R^n \setminus \R^d, \, \xi,\zeta \in \C^n.
\end{equation}
Or alternatively, if the reduced matrix $\cA := |t|^{n-d-1} A$ satisfies the classical elliptic and boundedness condition. We say that $L$ is $q$-elliptic if the coefficient matrix $A(X)$ satisfies \eqref{defEllip}--\eqref{defBdd} and $\cA(X)$ satisfies the condition \eqref{defpellipt} given in the next section. We refer the reader to the next section for the discussion of $q$-ellipticity, but let us cite one key property: $\cL$ is $q$-elliptic for all $q\in (1,+\infty)$ if and only if the matrix $\cA$ is real-valued. It means that $q$-ellipticity is a notion intrinsically linked to matrices with complex coefficients. 

\medskip

\medskip

The next assumption we put on $\cL$ will use Carleson measures. Let us first introduce the following definition.

\begin{defi} \label{defCMI}
We say that $f$ satisfies the Carleson measure condition if
$$d\mu(x,r) := \sup_{W_a(x,r)} |f|^2 dx \frac{dr}r$$
is a Carleson measure, that is, if there exists $C_a>0$ such that 
$$\|f\|_{CM,a}:= \sup_{x\in \R^d, \, r>0} \fint_{y\in B(x,r)} \int_{0}^r \sup_{W_a(y,s)} |f|^2 dy \frac{ds}s \leq C_a.$$
\end{defi}

Let us make a few remarks. It is easy to check that the fact that $f$ satisfies the Carleson measure condition doesn't depend on the choice of $a$, and that for any $a,b>0$
$$\|f\|_{CM,a} \leq C_{a,b} \|f\|_{CM,b}.$$
The quantity $f$ in the above definition can be a measurable function, the gradient of a function, or a matrix-valued function. So $|f|$ is respectively the absolute value, the vector norm, or the matrix norm. Furthermore, the Carleson measure condition forces $f$ to be locally Lipschitz.

\begin{defi} \label{defH1}
An elliptic operator $\cL = - \div |t|^{d+1-n} \cA \nabla$ with complex coefficient satisfies the hypothesis $(\cH^1)$ if the matrix $\cA$ can be decomposed as
\begin{equation} \label{decompA}
\cA = \begin{pmatrix} \cA_{1} & \cA_2 \\ \cB_3 & bI \end{pmatrix} + \cC,
\end{equation}
where
\begin{enumerate}[(i)]
\item $I$ is an $(n-d)$-identity matrix, $b$ is a real scalar function, and $\cB_3$ is a real matrix in $M_{(n-d)\times d}$,
\item $b$ is uniformly bounded from above and below,
\item the quantities $|t|\nabla_x \cB_3$, $|t|^{n-d}\div_t [|t|^{d+1-n}\cB_3]$, $|t|\nabla b$ and $\cC$ satisfy the Carleson measure condition, that is there exists $C>0$ such that
\begin{equation} \label{Ckappa}
\sup_{z\in \R^d, \, r>0} \fint_{y\in B(z,r)} \int_{0}^\infty \sup_{(x,t) \in W_1(y,s)} \left[  |t|^2 |\nabla_x \cB_3|^2 + |t|^2|\nabla b|^2 + |\cC|^2  \right] \frac{ds}s dy \leq C
\end{equation}
and for any $1\leq j\leq d$,
\begin{equation} \label{Ckappabis}
\sup_{z\in \R^d, \, r>0} \fint_{y\in B(z,r)} \int_{0}^\infty \sup_{(x,t) \in W_1(y,s)} |t|^{2(n-d)} \left|\sum_{\ell >d} \dr_{t_\ell} \left[(\cB_3)_{\ell j} |t|^{n-d-1} \right] \right|^2 \frac{ds}s dy \leq C
\end{equation}
where we write the entries of $\cB_3$ as $(\cB_3)_{d < \ell \leq n \atop{1\leq j \leq d} }$ and $\dr_{t_\ell}$ corresponds to the partial derivate with respect to the $\ell^{th}$ coordinate in $\R^n$.
\end{enumerate}

In addition, we say that $\cL$ satisfies the small Carleson hypothesis $(\cH^1_\kappa)$, if $\cL$ satisfies $(\cH^1)$ and the constant $C$ in \eqref{Ckappa}--\eqref{Ckappabis} can be chosen to be smaller than $\kappa$.
\end{defi}

\begin{rmk}
If the block matrix $\cB_3$ is not a real matrix, we could include its imaginary part in $\cC$ as long as $\Im \, \cB_3$ satisfies the appropriate Carleson measure condition. The assumptions that $b, \cB_3$ are real-valued are used in the proof of Lemma \ref{lemS<Na}, which is the key estimate in proving $S<N$ and $S<\kappa N + \Tr$. Observe that for our example in (4.6) of \cite{DFM3}, $\cB_3 $ has the form $(\cB_3)_{ij} = b_j(x) t_i/|t|$, and thus $|t|^{n-d}\div_t [|t|^{d+1-n}\cB_3]$ equals $0$ and satisfies the Carleson measure assumption trivially.
\end{rmk}

\medskip

We say that $u \in W^{1,2}_{loc}(\R^n \setminus \R^d)$ is a weak solution to $\cL u = 0$ if we have
$$\iint_{(x,t) \in \R^n} \cA \nabla u \cdot \nabla \bar \Psi \, dx\, \frac{dt}{|t|^{n-d-1}} = 0 \qquad \text{ for any } \Psi \in C^\infty_0(\R^n \setminus \R^d).$$
Our main theorem is the following result.

\begin{thm} \label{MainTh}
Let $a>0$, $M>0$, and $q\in (1,+\infty)$. There exists a constant $\kappa_0 := \kappa_0(n,a,M,q)>0$ such that if $\cL := -\div |t|^{d+1-n} \cA \nabla$ is a $q$-elliptic operator with the following properties
\begin{enumerate}[(i)]
\item it satisfies the hypothesis $(\cH^1_{\kappa_0})$,
\item the $q$-ellipticity constant $\lambda_q$ that appears in \eqref{defpellipt} is bigger that $M^{-1}$,
\item $\|\cA\|_\infty + \|b\|_\infty + \|b^{-1}\|_\infty \leq M$, where $b$ is the one that is defined via \eqref{decompA},
\end{enumerate}
then the Dirichlet problem (D$_q$) is well-posed. That is, for any $g\in L^q(\R^d)$, there exists a unique weak solution $u:=u_g \in W^{1,2}_{loc}(\R^n \setminus \R^q)$ to $\cL u = 0$ such that 
\begin{equation}\label{intro:qntcv2}
	\lim_{(z,r) \in \Gamma_a(x) \atop{r\to 0}} \frac{1}{\W_a(z,r)} \iint_{W_a(z,r)} u(y,s) \, dy\, ds = g(x) \qquad \text{ for almost any } x\in \R^d,
\end{equation}
and
$$\|\wt N_{a,2}(u)\|_q<+\infty.$$

Furthermore, there hold
\[ \|\wt N_{a,q}(u)\|_q \leq C \|g\|_q, \]
and 
\[ \|\wt N_{a,q}(u)\|_q \approx \|S_{a,q}(u)\|_q, \quad \|\wt N_{a,2}(u)\|_q \approx \|\wt N_{a,q}(u)\|_q, \]
with constants that depend only on the dimension $n$, $a$, $q$ and $M$. One also has the stronger convergence
\begin{equation}\label{intro:qntcv}
	\lim_{(z,r) \in \Gamma_a(x) \atop{r\to 0}} \frac{1}{\W_a(z,r)} \iint_{W_a(z,r)} |u(y,s)-g(x)|^q \, dy\, ds = 0 \qquad \text{ for almost any } x\in \R^d.
\end{equation}
\end{thm}

\begin{rmk}
	In particular \eqref{intro:qntcv} implies the pointwise non-tangential convergence of the solution $u$ to $g$:
	\[ \lim_{(z,r)\in \Gamma_a(x) \atop {r\to 0}} \left(\frac{1}{\W_a(z,r)} \iint_{W_a(z,r)} |u(y,s)|^q \, dy\, ds \right)^{1/q} = |g(x)| \qquad \text{ for almost any } x\in \R^d. \]
\end{rmk}

\medskip

Except for the main result in \cite{DFMAinfty}, the result above is the first treatment of the well-posedness of the Dirichlet problem for domains with higher co-dimensional boundary. 
In particular, in the case where the domain is $\Omega = \R^n \setminus \R^d$,  the well-posedness of (D$_q$) for $q$ small has never been considered before, even in the case where the operator $L=-\div A \nabla$ is such that $A$ has real coefficients. 
At some point in the proof, we need to use the saw-tooth domains $\{|t| > h(x)\}$ where $h:\R^{d} \to \R_+$ is some non-negative Lipschitz function. The strategy in codimension 1 is to use a bi-Lipschitz change of variable that sends $\{t > h(x)\}$ to the upper half plane, and to deduce the result on saw-tooth domains from the one in the upper half plane. This method is not available anymore in higher co-dimension (notice that even the boundary of the saw-tooth domain becomes an object of mixed dimensions in this context) and one of the main difficulties here is to establish desired estimates directly on saw-tooth domains. Fortunately, emerging new ideas brought considerably stronger results in the ``classical" co-dimension 1 setting of $\R^{n}_+ = \R^{d+1}_+$ as well.

Indeed, the methods used to prove Theorem \ref{MainTh} (and its result) can be easily adapted to the classical case where the boundary of the domain is of codimension 1, $\Omega = \R^n_+$, with an obvious reformulation. Here are some remarks that are important even in codimension 1:
\begin{itemize}
\item As pointed out above, contrary to \cite{DP}, we do not impose any conditions on the first $d$ rows of the matrix of coefficients.
\item We allow Carleson measure perturbations (addition of a matrix $\mathcal C$). This is vital in our method and formally speaking new even in the classical scenario. Indeed, known results about the Carleson measure perturbation are tied up to real coefficients \cite{FKP} or to perturbing from the $t$-independence matrix \cite{HMaMo, AA}, which is not the setting of the present paper. 
\item This is the second time to the best of our knowledge (after \cite{DP}) that the Dirichlet problem for elliptic operators with complex coefficients whose gradients are Carleson measures is attacked. We keep the remarkable idea of Dindo\v s and Pipher to use a notion of $q$-ellipticity for elliptic operators, and we improve it in several ways (see the first two points above and the last two points below).
\item We offer a new proof of the existence of solutions to the Dirichlet problem  in unbounded domains. A general difficulty  is  to prove that the quantity $\|\wt N_{a,q} u\|_q$ is actually finite for a large class of solutions. Indeed, the formal proof of well-posedness consists of showing that $\|\wt N_{a,q} u\|_q \leq C \|S_{a,q} u\|_q$ and then that $\|S_{a,q} u\|_q \leq \eta \|\wt N_{a,q} u\|_q + \|\Tr u\|_q$. If $\eta$ is small, the conclusion $\|\wt N_{a,q} u\|_q \leq C \|\Tr u\|_q$ holds provided that $\|\wt N_{a,q} u\|_q$ is {\em a priori} finite. In the case of co-dimension 1 one can rely on a plethora of known results for smooth coefficients, on layer potential techniques, and other methods. The proof in this paper is self-contained (which partially explains the length). 

\item We prove local rather than global $S<N$ and $N<S$ estimates. Even in co-dimension one, such results in full generality are only available for real coefficients, typically by localization from the global estimates \cite{KePiDrift,DPP, HKMP}. Unfortunately, localizing the global $S<N$ and $N<S$ estimates requires sufficient decay of solutions and a maximum principle, generally failing for complex coefficients. 
However,  the local bounds are important both in the present argument and for the future use (for instance, in the extrapolation techniques on uniformly rectifiable domains) -- a more detailed discussion will be presented in the next paragraph. 

\item We prove that, when the trace is smooth, the solutions $u$ obtained by the Lax-Milgram theorem match the solution(s) of the Dirichlet problem (see \cite{Ax, BM} for a discussion of importance and possible failure of this property under various circumstances). In order to do it, we prove, in particular, the finiteness of the quantity $\int_\Omega |\nabla u|^2 |u|^{q-2} dm$ (where $dm$ is the Lebesgue measure in the co-dimension 1 case, and $dm(x) = \dist(x,\Omega)^{d+1-n} dx$ when the boundary has higher codimension) for $q\geq 2$. This finiteness is new even in the codimension 1 case, and holds under the only condition that the operator is $q$-elliptic - in particular we do not require $\cL$ to satisfy $(\cH^1)$ - and that the domain is the complement of an Ahlfors regular set. Our control on (the energy of) solutions will allow us to derive the global $S<N$ and $N<S$ estimates from the local ones (contrary to the classical approach). A more detailed discussion can be found in Subsection \ref{ssES}.

\item We improve the reverse H\"older estimates (that can be seen as weak Moser estimates) proven in \cite[Lemmas 2.6 and 2.7]{DP}, and we also give a boundary version. These results hold under the sole assumption of $q$-ellipticity and will hopefully be useful in a wide variety of problems. They will be presented in Subsection \ref{ssES} as well.
\end{itemize}

\medskip

If an operator $L=-\div A \nabla$ has real coefficients, slight modifications of the argument of Theorem \ref{MainTh} gives the result below. Before stating it, let us introduce some additional notation. Thanks to the classical Moser's estimate, we can directly work with the classical non-tangential maximal function $N_{a}$ (rather than in the average sense, see \eqref{defNq} \eqref{defw}): For any function $v\in L^\infty_{loc}(\Omega)$ and any $x\in \R^d$, we define
\begin{equation}\label{defNq22}
N_{a}(v)(x) = \sup_{\wh \Gamma_a(x)}  v.
\end{equation}
We say that $f$ satisfies the {\bf real} Carleson measure condition if 
\[d\mu(x,t):= |f(x,t)|^2 dx\, \frac{dt}{|t|^{n-d}}\]
is a Carleson measure. Note that if the elliptic operator $L$ has real coefficients, the assumption $(\cH^1_\kappa)$ is effectively weakened by requiring only real Carleson measure condition in (iii).
Clearly, the assumption that $b$ and $\cB_3$ are real becomes void.

\begin{thm}[Analogue of Theorem \ref{MainTh} for real coefficients] \label{MainThR}
Let $a>0$, $M>0$, and $q\in (1,+\infty)$. There exists a constant $\kappa_0 := \kappa_0(n,a,M,q)>0$ such that if $\cL := -\div |t|^{d+1-n} \cA \nabla$ is an elliptic operator with real coefficients that satisfies
\begin{enumerate}[(i)]
\item the hypothesis $(\cH^1_{\kappa_0})$,
\item the ellipticity constant is bounded from below by $M^{-1}$, i.e.
\[\inf_{x\in \Omega, \, \xi \in \R^n\atop |\xi| = 1} \cA(x)\xi \cdot \xi  \geq M^{-1},\]
\item $\|\cA\|_\infty + \|b\|_\infty + \|b^{-1}\|_\infty \leq M$, where $b$ is defined as in \eqref{decompA},
\end{enumerate}
then the Dirichlet problem (D$_q$) is well-posed. That is, for any $g\in L^q(\R^d)$, there exists a unique weak solution $u:=u_g \in W^{1,2}_{loc}(\R^n \setminus \R^q)$ to $\cL u = 0$ such that 
\begin{equation}\label{intro:qntcv22}
	\lim_{(z,r) \in \Gamma_a(x) \atop{r\to 0}} \frac{1}{\W_a(z,r)} \iint_{W_a(z,r)} u(y,s) \, dy\, ds = g(x) \qquad \text{ for almost any } x\in \R^d,
\end{equation}
and
$$\|N_{a}(u)\|_q<+\infty.$$
Furthermore,
\[ \|S_{a,q}(u)\|_q \approx \|N_{a}(u)\|_q \leq C \|g\|_q, \]
where the constants depend only on the dimension $n$, $a$, $q$ and $M$. 
\end{thm}

A corollary of Theorem \ref{MainThR} is the well-posedness of the Dirichlet problem (D$_q$) when the domain $\Omega$ is the complement of the graph of a Lipschitz function with small Lipschitz constant. Let us be more precise. Set $\Gamma = \{(x,\varphi(x)), \, x\in \R^d\}$ be the graph of a Lipschitz function $\varphi: \, \R^d \to \R^{n-d}$, and then $\Omega = \R^n \setminus \Gamma$. 
The non-tangentially maximal function $N_{a,\Gamma}$ on $\Gamma$ is the one defined as
$$N_{a,\Gamma}(u)(x') = N_{a}(v)(x) \qquad \text{ for } x' =(x,\varphi(x)) \in \Gamma,$$
where $v(y,s) = u(y, s- \varphi(x))$. Observe that the definition of $N_{a,\Gamma}$ makes sense if either $a$ or the Lipschitz constant of $\varphi$ is small. Indeed, in either case, the set $\wh \Gamma_a(x)$ stays inside the domain of definition of $v$.

\begin{cor} \label{MainCor}
Let $\Gamma := \{(x,\varphi(x)), \, x\in \R^d\}$, where $\varphi: \, \R^d \to \R^{n-d}$ is a Lipschitz function, and $\Omega = \R^n \setminus \Gamma$. Choose $\alpha>0$ and define 
$$D_\alpha(X) := \left( \int_\Gamma |X-y|^{-d-\alpha} d\sigma(y) \right)^{-1/\alpha},$$
where $\sigma$ is the $d$-dimensional Hausdorff measure on $\Gamma$. Consider the elliptic operator 
$$L := -\div D_\alpha (X)^{d+1-n} \nabla.$$

Let $a>0$ and $q\in (1,+\infty)$. There exists a constant $\kappa_0 := \kappa_0(n,\alpha,a,q)>0$ such that if $\|\nabla \varphi\|_\infty \leq \kappa_0$, then the Dirichlet problem (D$_q$) is well-posed. That is, for any $g\in L^q(\Gamma)$, there exists a unique weak solution $u:=u_g \in W_{loc}^{1,2}(\Omega)$ to $L u = 0$ such that
$$\lim_{X \in \Gamma_a(x')} \frac{1}{|B(X,\delta(X)/2)|} \iint_{B(X,\delta(X)/2)} u(Y) \, dY = g(x') \qquad \text{ for almost any } x'\in \Gamma, $$
and
$$\|N_{a}(u)\|_{L^q(\Gamma)} <+\infty.$$
Furthermore, 
$$\|N_{a}(u)\|_{L^q(\Gamma)} \leq C \|g\|_{L^q(\Gamma)},$$
where the constant $C>0$ depends only on the dimension $n$, $\alpha$, $a$, and $q$. 
\end{cor}

\begin{rmk} The operator above is the one for which $A^\infty$ property of elliptic measure and hence, solvability of the Dirichlet problem for {\it some large} $q<\infty$ was proved in \cite{DFMAinfty}. Clearly, the Corollary could be extended to a more general class of elliptic operators which are tied up to the change of variables defined below in the proof. However, it does not seem to be possible to work, e.g., with the Euclidean distance and so the precise choice of the coefficients is rather delicate.  \end{rmk}

\bp Assume that $\|\nabla \varphi\| \ll 1$. We use the same bi-Lipschitz change of variable $\rho$ as the one used in \cite[equation (3.3)]{DFMAinfty}, that is 
\begin{equation} 
\rho(x,t) = (x,\eta_{|t|}*\varphi(x)) + h(x,|t|) \, R_{x,|t|}(0,t) \quad \text{ for } (x,t)\in \R^n,
\end{equation}
where $\eta_r$ is a mollifier, $R_{x,r}$ is a linear isometry of $\R^n$, and $h(x,r) > 0$ is a dilation factor,
that we rapidly discuss now. We construct $R_{x,r}$ (with a convolution formula and projections)
so that it maps $\R^d$ to the $d$-plane $P(x,r)$ tangent to 
$\Gamma_r := \{(x,\eta_r*\varphi(x)), \, x\in \R^d\}$ at the point 
$\Phi_r(x):= (x,\eta_r*\varphi(x))$, and hence also $R_{x,r}$ maps $\R^{n-d} = (\R^d)^\bot$ 
to the orthogonal plane to $P(x,r)$ at $\Phi_r(x)$.
The map $\rho^{-1}$ sends $\Omega$ to $\R^n \setminus \R^d$ and $L$ to a (real-coefficient) elliptic operator $\cL = - \div A \nabla$. 
Lemma 3.40 in \cite{DFMAinfty} and the fact that $D_\alpha(X) \simeq \dist(X,\Gamma)$ prove that $A$ satisfies \eqref{defEllip}--\eqref{defBdd}. 
Moreover, if $\|\nabla \varphi\|_\infty$ is small enough, Lemma 6.22 in \cite{DFMAinfty} establishes that $\cA := |t|^{n-d-1} A$ satisfies assumption (i)--(iii) of Theorem \ref{MainThR} \footnote{The smallness of the Carleson measure is not explicitly written in \cite{DFMAinfty}, but the smallness holds as long as the Ahlfors measure $\sigma$ is close enough to a flat measure, which is the case if $\varphi$ has a small Lipschitz constant and $\sigma$ is the $d$-dimensional Hausdorff measure on $\Gamma$.}. Theorem \ref{MainThR} implies now that the Dirichlet problem (D$_q$) is well-posed when the domain is $\R^n \setminus \R^d$ and the operator is $\cL$, and therefore, by using again the change of variable $\rho$, we conclude that the Dirichlet problem (D$_q$) is well-posed when the domain is $\Omega \setminus \Gamma$ and the elliptic operator is $L$.
\ep

\subsection{Local bounds}

In this paragraph, we stay with $\Omega = \R^n \setminus \R^d$ and $\Gamma = \R^d$. We write $(x,t)$ or $(y,s)$ for a running point of $\R^n = \R^d \times \R^{n-d}$. 

\medskip

Let $a>0$ be fixed. 
We need local versions of the square function and the non-tangential maximal function. The following definitions are in parallel to previous definitions \eqref{defSq}, \eqref{defNq} and \eqref{defw}. If $\Psi$ is a cut-off function -  that is we ask $\Psi$ to be at least a locally Lipschitz function and to satisfy $0\leq \Psi \leq 1$ everywhere - then we define $S_{a,q}(\cdot|\Psi)$ and $\wt N_{a,q}(\cdot|\Psi)$ as
\begin{equation}\label{defSq2}
S_{a,q}(v|\Psi)(x) := \left( \iint_{\wh \Gamma_a(x)} |\nabla v(y,s)|^2 |v(y,s)|^{q-2} \Psi(y,s) dy \frac{ds}{|s|^{n-2}} \right)^\frac1q,
\end{equation}
and
\begin{equation}\label{defNq2}
\wt N_{a,q}(v|\Psi)(x) = \sup_{\Gamma_a(x)}  (v|\Psi)_{\W,a,q}
\end{equation}
where $\vW$ is defined on $\R^{d+1}_+:= \R^d \times (0,+\infty)$ by
\begin{equation} \label{defw2}
	(v|\Psi)_{\W,a,q}(z,r) := \left( \frac1{|W_a(z,r)|} \iint_{W_a(z,r)}|v(y,s)|^q \Psi(y,s) dy \, ds   \right)^\frac1q.
\end{equation}

\medskip

Particularly we are interested in the following cut-off functions.
Choose a function $\phi \in C^\infty_0(\R_+)$ such that $0 \leq \phi \leq 1$, $\phi \equiv 1$ on $(0,1)$ and $\phi \equiv 0$ on $(2,+\infty)$, $\phi$ non-increasing and $|\phi'| \leq 2$. If $e(x)$ is a positive $a^{-1}$-Lipschitz function, we define $\Psi_e$ as
$$\Psi_e(x,t) = \phi\left( \frac{e(x)}{|t|} \right).$$
If $B \subset \R^d$ is a ball with radius  bigger or equal  to $l$, 
$$\Psi_{B,l}(x,t) = \phi\left( \frac{a|t|}{l}\right) \phi\left(1+\frac{\dist(x,B)}{l}\right).$$
We keep in mind that the functions $\Psi_e$ and $\Psi_{B,l}$ depend on $a$, but we don't write the dependence to lighten the notation. Whatever choice we make for $e>0$, $l$, and $B$, observe that $\Psi_e \Psi_{B,l}$ is a smooth cut-off function which is compactly supported in $\R^n \setminus \R^d$. The following results hold.

\begin{thm} \label{ThLocal}
Let $a>0$ and $q\in (1,+\infty)$. Let $\cL = - \div |t|^{d+1-n}\cA \nabla$ be a $q$-elliptic operator that satisfies $(\cH^1$). For any $a^{-1}$-Lipschitz function $e$, any $l>0$ and any ball $B$ whose radius is bigger than $l$, and any weak solution $u\in W^{1,2}_{loc}(\R^n\setminus \R^d)$ to $\cL u = 0$, we have
\begin{enumerate}
\item if $k>2$ and $1<p<\infty$,
$$\|S_{a,q}(u|\Psi_e^k \Psi_{B,l}^k)\|_p \leq C \|\wt N_{a,q}(u|\Psi_e^{k-2} \Psi_{B,l}^{k-2})\|_p,$$
\item if $k>12$,
$$\|\wt N_{a,q}(u|\Psi_e^k \Psi_{B,l}^k)\|_q^q \leq C \|S_{a,q}(u|\Psi_e^{k-12} \Psi_{B,l}^{k-12})\|_q^q + C \iint_{(y,s)\in \R^n}  |u|^q \Psi_{B,l}^{k-3} \dr_r[\Psi_e^{k-3}] \, \frac{ds}{|s|^{n-d-1}} \, dy,$$
\item if $k>12$,
$$\|\wt N_{a,q}(u|\Psi_e^k \Psi_{B,l}^k)\|_q^q \leq C \|S_{a,q}(u|\Psi_e^{k-12} \Psi_{B,l}^{k-12})\|_q^q + C \iint_{(y,s)\in \R^n} |u|^q \Psi_{e}^{k-3} \dr_r[-\Psi_{B,l}^{k-3}] \, \frac{ds}{|s|^{n-d-1}} \, dy,$$
\end{enumerate}
where $\dr_r$ represents the derivative is the radial direction of $t$.
The constants in the above inequalities depend on $a$, $q$, $n$, a lower bound for the value $\lambda_q$ in \eqref{defpellipt}, $\|\cA\|_\infty$, an upper bound for the value $C$ in \eqref{Ckappa}, $\|b\|_\infty + \|b^{-1}\|_\infty$, and $k$. In (1), the constant depends also on $p$.
\end{thm}

The proofs of the three parts of the Theorem are in Lemmas \ref{lemS<Nd}, \ref{lem6.9} and \ref{lem6.10}, respectively.

As we will see in the next paragraph, the fact that $u\in W^{1,2}_{loc}(\R^n \setminus \R^d)$ is enough to ensure that the terms in Theorem \ref{ThLocal} are finite. On the other hand, nothing guarantees for instance the finiteness of the quantities $\|S_{a,q}(u|\Psi_{B,l}^k)\|_p$ or $\|\wt N_{a,q}(u|\Psi_{B,l}^{k})\|_p$ that are obtained by taking $e \to 0$, i.e., the boundary behavior.

\medskip

The above definitions using the functions $\Psi_e$ and $\Psi_{B,l}$ are smooth versions of the standard local  square function and non-tangential maximal function. They are more powerful since they will allow us to hide terms from the right hand side of an estimate.
In order to build the reader's intuition, let us show some examples of the use of the theorem above. 

\medskip

If $\epsilon>0$ and $l>0$, we define the $q$-adapted square function $S_{a,q}^{\epsilon,l}$ and the non-tangential maximal function $\wt N_{a,q}^{\epsilon,l}$ in a similar manner to $S_{a,q}$ and $\wt N_{a,q}$, but with the truncated cones
$$\wh \Gamma_a^{\epsilon,l}(x) := \{(y,s) \in \wh \Gamma_a(x), \, \epsilon < |s| < l/a \} \ \text{ and } \ \Gamma_a^{\epsilon,l}(x) := \{(z,r) \in \Gamma_a(x), \, \epsilon < r < l/a \}$$
respectively. Take a ball $B \subset \R^d$ with radius $l$ and then $\epsilon >0$. We construct the Lipschitz function $e$ as $e \equiv \epsilon$ on $\R^d$. We choose then $B' = 2B$ and we write $\Psi$ for $\Psi_{e} \Psi_{B',l}$. Due to (1) of Theorem \ref{ThLocal}, we have
$$\|S_{a,q}(u|\Psi^3)\|_p \leq C \|\wt N_{a,q}(u|\Psi)\|_p.$$
However, our choice of $\Psi$ is such that 
$$S_{a,q}(u|\Psi^3)(x) \geq S_{a,q}^{\epsilon,l} (u)(x) \qquad \text{ for } x\in B,$$
and in addition, one can check that 
$$\|\wt N_{a,q}(u|\Psi)\|_p \leq C_p \|\wt N_{a,q}^{\epsilon/2,2l}(u)\|_{L^p(4B)}.$$
So Theorem \ref{ThLocal} implies that 
\begin{equation} \label{IntroSNa}
\|S_{a,q}^{\epsilon,l} (u)\|_{L^p(B)} \leq C_p \|\wt N_{a,q}^{\epsilon/2,2l}(u)\|_{L^p(4B)},
\end{equation}
which is a more customary statement of the local $S<N$ bound.

\medskip

By reasoning similar to the one above, (3) of Theorem \ref{ThLocal} gives a bound
\begin{equation} \label{IntroSNb}
\|\wt N_{a,q}^{\epsilon,l}(u)\|_{L^q(B)} \leq C \|S_{a,q}^{\epsilon/2,2l} (u)\|_{L^p(10B)} + C \iint_{(y,s) \in \R^n} |u|^q \Psi_{e}^{k-3} \dr_r[-\Psi_{B',l}^{k-3}] \, \frac{ds}{|s|^{n-d-1}} \, dy.
\end{equation}
Since $\phi$ is non-increasing and non-negative, the term $\Psi_{e}^{k-3} \dr_r[-\Psi_{B',l}^{k-3}]$ is non-negative. Moreover, by the construction of  $\Psi_e$ and $\Psi_{B',l}$, the term $\Psi_{e}^{k-3} \dr_r[-\Psi_{B',l}^{k-3}]$ is bounded by $\frac Cl$ and supported in $\{(y,s)\in \R^n, \, y\in 10B, \, l \leq a|s| \leq 2l\}$.
So the last term in \eqref{IntroSNb} is bounded by $Cl^d$ times the average of the function $|u|^q$ over a Whitney box $W_B$ associated to the ball $B$. If we assume that $\iint_{W_B} u = 0$, the Poincar\'e inequality implies that the average of $|u|^q$ over $W_B$ can be bounded by $\|S_{a,q}^{\epsilon/2,2l} (u)\|_{L^q(10B)}$. We obtain then, if the average of $u$ over the Whitney box $W_B$ is 0, that  
\begin{equation} \label{IntroSNc}
\|\wt N_{a,q}^{\epsilon,l}(u)\|_{L^q(B)} \leq C \|S_{a,q}^{\epsilon/2,2l} (u)\|_{L^q(10B)}.
\end{equation}
The latter is a customary statement of the local $N<S$ bound.

\subsection{Reverse H\"older estimates (weak Moser) and energy solutions}

\label{ssES}

The results in the present paragraph require a lot fewer assumptions than before, either on the domain or on the elliptic operator $L$.  The ball in $\R^n$ with center $X\in \R^n$ and radius $r$ is denoted by $B(X,r)$. Let $d<n-1$ be a (not necessarily integer) positive number, and $\Gamma$ be a $d$-dimensional Ahlfors regular set, that is, there exists $C>0$ and a measure $\sigma$ on $\Gamma$ such that 
\begin{equation} \label{defADR}
C^{-1} r^{d} \leq \sigma(B(X,r)) \leq C r^{d} \qquad \text{ for any } X \in \Gamma, \; r>0.
\end{equation}
It is well known that if $\Gamma$ is Alfors regular, then \eqref{defADR} also holds for $\sigma = \mathcal H^d_{|\Gamma}$, the $d$-dimensional Hausdorff measure on $\Gamma$ (and a different constant $C$), see \cite[Theorem 6.9]{Mattila}.
Set now $\Omega = \R^n \setminus \Gamma$, $\delta(x) = \dist(x,\Gamma)$, and a measure $m$ defined as
$$m(E) = \int_E \delta(X)^{d+1-n} dX,$$
that is $dm = \delta^{d+1-n} dX$. Note that when $\Gamma = \R^d$ and $\Omega = \R^n \setminus \R^d$, for any $X=(x,t) \in \R^n$, we have $\delta(X) = |t|$ and $dm(X) = |t|^{d+1-n} dt dx$, and we recognize the weight used in the previous subsections.

\medskip

We say an operator $L = - \div A \nabla$ on $\Omega$ is elliptic \textit{with weight $\delta(X)^{d+1-n}$} if there exists $C>0$ such that the complex matrix $A$ satisfies
\begin{equation} \label{defEllip2}
 \Re( A(X)\xi \cdot \bar \xi) \geq C^{-1} \delta(X)^{d+1-n} |\xi|^2 \qquad \text{ for } X\in \Omega, \, \xi \in \C^n,
\end{equation}
and
\begin{equation} \label{defBdd2}
|A(X)\xi \cdot \bar \zeta| \leq C \delta(X)^{d+1-n} |\xi| |\zeta| \qquad \text{ for } X\in \Omega, \,  \xi,\zeta \in \C^n.
\end{equation}
These two conditions are the generalization of \eqref{defEllip}--\eqref{defBdd} in the case where $\Gamma$ is Ahlfors regular. The operator $L$ is said to be {\it $q$-elliptic} if $L$ is elliptic and $\cA(X) := \delta^{n-d-1}(X)A(X)$ satisfies the condition \eqref{defpellipt} given in the next section. In addition, we say that $u$ is a weak solution to $Lu=0$ if
$$\int_{\Omega} A\nabla u \cdot \nabla \bar \Psi \, dX = \int_{\Omega} \cA \nabla u \cdot \nabla \bar \Psi \, dm  = 0 \qquad \text{ for any } \Psi \in C^\infty_0(\Omega).$$
In the complex case, the classical Moser's estimates (i.e. $L^\infty$-local bounds) don't necessarily hold but we have the following weaker version.

\begin{prop} \label{lemMoserI}
Let $L = - \div A \nabla$ be a $q$-elliptic operator. For any ball $B \subset \R^n$ of radius $r$ that satisfies $3B \subset \Omega$ and any weak solution $u\in W_{loc}^{1,2}(\Omega)$ to $L u = 0$, we have
$$\int_{B} |\nabla u |^2 |u|^{q-2} dm \leq \frac{C}{r^2} \int_{2B} |u|^q dm.$$

Moreover, the following reverse H\"older estimates hold.
\begin{enumerate}[(i)]
\item If $q>2$,
$$ \left( \frac{1}{m(B)} \int_B |u|^q dm \right)^{\frac1q} \leq C \left( \frac{1}{m(2B)} \int_{2B} |u|^2 dm \right)^{\frac12}.$$
\item If $q<2$,
$$ \left( \frac{1}{m(B)} \int_B |u|^2 dm \right)^{\frac12} \leq C \left( \frac{1}{m(2B)} \int_{2B} |u|^q dm \right)^{\frac1q}.$$
\end{enumerate}

In all three inequalities, the constant $C>0$ depends only on $n$, $q$, a lower bound on constant $\lambda_q$ in \eqref{defpellipt}, and the ellipticity constant in \eqref{defEllip2}--\eqref{defBdd2}.
\end{prop}

The proposition validates the fact that the quantities invoked in Theorem \ref{ThLocal} are indeed finite. The analogue of this result in codimension 1 is written in the next subsection, and we will see that the bounds (in co-dimension 1) when $q>2$ were already stated in \cite[Lemma 2.6]{DP}, but when $q<2$, our proposition is an improvement of \cite[Lemma 2.7]{DP}.

\medskip

Before we state our next result, let us introduce a bit of the theory given in \cite{DFMprelim}. We denote by $W$  the weighted Sobolev space of functions $u\in L^1_{loc}(\Omega)$ whose distribution gradient in $\Omega$ lies in $L^2(\Omega,dm)$:
$$W = \left\{u \in L^1_{loc}(\Omega):\, \|u\|_W := \left(\int_\Omega |\nabla u |^2 dm\right)^\frac12 < +\infty\right\}.$$
Clearly $W$ is contained in $W_{loc}^{1,2}(\Omega)$.
Lemma 3.3 and Lemma 4.13 in \cite{DFMprelim} establishes that 
\[ W = \{ f \in L^2_{loc}(\R^n,dm), \, \nabla u \in L^2(\R^n, dm)\}. \] 
This observation is useful to see that we will not have any problems integrating $u\in W$ across the boundary of $\Omega$. We denote by $\mathcal M(\Gamma)$ the set of measurable functions on $\Gamma$, and we set
$$H:= \{g\in \mathcal M(\Gamma), \, \int_{\Gamma} \int_{\Gamma} \frac{|g(x)-g(y)|^2}{|x-y|^{d+1}} \, d\sigma(x) \, d\sigma(y).\}$$
Now, by \cite[Theorem 3.13]{DFMprelim}, the trace operator  $\Tr: \, W \to H$ defined by
$$\Tr u(x) = \lim_{\epsilon\to 0} \fint_{y \in B_\epsilon(x)} \fint_{|s|\leq \epsilon} u(y,s) \, ds\, dy \qquad x\in \Gamma$$
is linear and bounded. We are ready to state a version of Proposition \ref{lemMoserI} at the boundary.

\begin{prop} \label{lemMoserB3}
Let $L:= -\div A \nabla$ be an $q$-elliptic operator.
For any weak solution $u\in W$ to $Lu = 0$ and for any  $B$ a ball of radius $r$ centered on $\Gamma$ that satisfies $\Tr u = 0$ on $3B$ there holds
\[\int_{B} |u|^{q-2} |\nabla u|^2 \, dm  \leq \frac{C}{r^2} \int_{2B \setminus B} |u|^q \, dm.\]
Furthermore, if $q>2$, one has
\[\left( \frac{1}{m(B)} \int_B  |u|^q \, dm \right)^{\frac{1}{q}} \leq C  \left(\frac{1}{m(2B)} \int_{2B}  |u|^2   \, dm \right)^{\frac{1}{2}},\]
and if $q<2$, we have
\[\left(\frac{1}{m(B)} \int_B  |u|^2 \, dm \right)^{\frac{1}{2}} \leq C  \left(\frac{1}{m(2B)} \int_{2B}  |u|^q   \, dm\right)^{\frac{1}{q}}.\]
The constant $C>0$ depends only on $n$, $q$, a lower bound on constant $\lambda_q$ in \eqref{defpellipt}, and $\|A\|_\infty$.
\end{prop}

The following result (proved as Lemma 9.3 in \cite{DFMprelim} in the case where $\cA$ has real coefficients but valid with the same proof in the present setting) gives the existence of weak solutions in $W$.

\begin{lem} \label{lemDirichlet}
Let $L:= -\div A \nabla$ be an elliptic operator. For any $g\in H$, there exists a unique $u_g\in W$ such that 
$$\int_{\Omega} A \nabla u_g \cdot \nabla \bar \varphi\, dX = 0, \qquad \text{ for any } \varphi \in C^\infty_0(\Omega)$$
and $\Tr u_g = g$ $\sigma$-a.e. on $\Gamma$. Moreover $\|u\|_W \leq C\|g\|_H$.
\end{lem} 

A weak solution $u \in W^{1,2}_{loc}(\Omega)$ is called an {\bf energy solution} to $L u = 0$ if $u\in W$ and $\Tr u \in C^\infty_0(\Gamma)$. Since $\Gamma$ is a closed set in $\R^n$, by extension theorem of Whitney type (see \cite[Chapter VI Section 2.2]{SteinSI}), the assumption $\Tr u\in C^\infty_0(\Gamma)$ implies that there exists $g\in C^\infty_0(\R^n)$ such that $\Tr u = \Tr g = g$ for a.e. $x\in \Gamma$. Since $\Tr g \in H$, Lemma \ref{lemDirichlet} shows that there exists a unique energy solution $u:=u_g$ to $L u = 0$ that satisfies $\Tr u = \Tr g$.

\begin{thm} \label{theoES}
Let $L$ is a $q$-elliptic operator. For any energy solution $u\in W$ to $L u = 0$, we have
\begin{enumerate}[(i)]
\item if $q\geq 2$,
$$\int_\Omega |\nabla u|^2 |u|^{q-2} \, dm < +\infty;$$
\item if $q \in (1,2)$, there exists a ball $B$ such that 
$$\int_{\Omega\sm B} |\nabla u|^2 |u|^{q-2} \, dm < +\infty.$$
\end{enumerate}
\end{thm}

\medskip

As we have previously mentioned, we prove the existence in the Dirichlet problem by proving that the energy solutions satisfy $\|N_{a,q}(u)\|_q <+\infty$ whenever $L$ is $q$-elliptic (the range of such $q$ is, as discussed in the next section, an open subset of $(1,+\infty)$ which is symmetric around $2$). Thanks to the a-priori finiteness proved in Theorem \ref{theoES}, we will be able to pass the local estimates in Theorem \ref{ThLocal} to global estimates by taking $e \to 0$, $B \nearrow \R^d$, and $l\to +\infty$, whenever $u$ is an energy solution (see Section \ref{SExistence}).

\subsection{The analogues of the results from the previous section in domains with co-dimension 1 boundaries.} As mentioned above, all our results have analogues in the upper half space or, respectively, a domain above an $(n-1)$-dimensional Lipschitz graph in $\R^n$. This statement is also valid for results of the previous section, but the geometric conditions become slightly more involved and for that reason we choose to restate the results carefully.

In this subsection, we say that $L = - \div A \nabla$ is a $q$-elliptic operator if the matrix $A$ lies in $L^\infty(\Omega)$ and $A$ satisfies \eqref{defpellipt}.

\begin{prop} \label{lemMoserI2}
Let $\Omega$ be a  domain in $\R^n$ and let $L = - \div A \nabla$ be a $q$-elliptic operator. For any ball $B \subset \R^n$ of radius $r$ that satisfies $2B \subset \Omega$ and any weak solution $u\in W_{loc}^{1,2}(\Omega)$ to $L u = 0$, we have
$$\int_{B} |\nabla u |^2 |u|^{q-2} dX \leq \frac{C}{r^2} \int_{2B} |u|^q dX.$$

Moreover, the following reverse H\"older estimates hold.
\begin{enumerate}[(i)]
\item If $q>2$,
$$ \left( \frac{1}{|B|} \int_B |u|^q dX \right)^{\frac1q} \leq C \left( \frac{1}{|2B|} \int_{2B} |u|^2 dX \right)^{\frac12}.$$
\item If $q<2$,
$$ \left( \frac{1}{|B|} \int_B |u|^2 dX \right)^{\frac12} \leq C \left( \frac{1}{|2B|} \int_{2B} |u|^q dX \right)^{\frac1q}.$$
\end{enumerate}

In all three inequalities, the constant $C>0$ depends only on $n$, $q$, a lower bound on constant $\lambda_q$ in \eqref{defpellipt}, and $\|A\|_\infty$.
\end{prop}

Observe that when $q\geq2$, the above result is the same as \cite[Lemma 2.6]{DP}. However, when $q<2$ and under the assumptions of Proposition \ref{lemMoserI2}, \cite[Lemma 2.7]{DP} states that for any $\epsilon >0$, we can find $C_\epsilon>0$ such that 
$$\int_{B} |\nabla u |^2 |u|^{q-2} dX \leq \frac{C_\epsilon}{r^2} \int_{2B} |u|^q dX + \frac{\epsilon}{r^2} \left(\int_{2B} |u|^2 dX\right)^\frac q2$$
and
$$ \left( \frac{1}{|B|} \int_B |u|^2 dX \right)^{\frac12} \leq C_\epsilon \left( \frac{1}{|2B|} \int_{2B} |u|^q dX \right)^{\frac1q} + \epsilon \left( \frac{1}{|2B|} \int_{2B} |u|^2 dX \right)^\frac12.$$
Our result is stronger than the one in \cite{DP} because we can remove the second terms of the two right-hand sides above.

\medskip

For the sequel, let us introduce three topological conditions on $\Omega$. We say that $\Omega$ satisfies the interior Corkscrew point condition when there exists a constant $C_1 > 0$ such that
for $x\in \dr\Omega$ and $0 < r < \diam(\dr\Omega)$, 
\begin{equation} \label{2.7}
\text{one can find a point $A_{x,r} \in \Omega \cap B(x,r)$ such that 
$B(A_{x,r},C_1^{-1} r) \subset \Omega$.}
\end{equation}
Similarly, we say that $\Omega$ satisfies the exterior Corkscrew point condition if the complement $\Omega^c$ satisfies the interior Corkscrew condition.
We also say that $\Omega$  satisfies the Harnack chain condition if there is a constant $C_2 \geq 1$ and, for each $\Lambda \geq 1$, 
an integer $N \geq 1$ such that, whenever $X,Y\in \Omega$ and $r \in (0, \diam(\dr\Omega))$ 
are such that
\begin{equation} \label{cs1} 
\min\{\dist(X,\dr\Omega),\dist(Y,\dr\Omega)\} \geq r \ \text{ and } \  |X-Y| \leq \Lambda r,
\end{equation}
we can find a chain of $N+1$ points $Z_0=X, Z_1,\dots, Z_{N} = Y$ in $\Omega$ such that
\begin{equation} \label{cs2}
C_2^{-1} r \leq \dist(Z_i,\dr\Omega) \leq C_2 \Lambda r
\text{ and } 
|Z_{i+1}-Z_i| \leq \dist(Z_i,\dr\Omega)/2 
\text{ for } 1 \leq i \leq N.
\end{equation}

Similarly to the higher codimension case, we define the space
$$W := \left\{u \in L^1_{loc}(\Omega):\, \|u\|_W := \left(\int_\Omega |\nabla u |^2 dX\right)^\frac12 < +\infty\right\},$$
which is clearly is contained in $W_{loc}^{1,2}(\Omega)$. If $\Omega$ satisfies the interior Corkscrew point condition, the exterior Corkscrew point condition, and the Harnack chain condition, and if the boundary $\dr \Omega$ is Ahlfors regular - i.e. verifies \eqref{defADR} with $d=n-1$ - then we can define notion of trace on $W$, that is there exists a bounded operator  $\Tr$ from $W$ to $L^2_{loc}(\Gamma,\sigma)$ such that $\Tr u = u$ if $u\in W \cap C^0(\overline \Omega)$.

We are now ready for the analogue of Proposition \ref{lemMoserI2} at the boundary, which is completely new.

\begin{prop} \label{lemMoserB33}
Let $\Omega$ satisfy the Corkscrew point condition and the Harnack chain condition, and be such that its boundary $\dr \Omega$ is Ahlfors regular of dimension $n-1$. Let $L:= -\div A \nabla$ be a $q$-elliptic operator.
Let $u\in W$ be a weak solution to $Lu = 0$ in $\Omega$ and  $B$ be a ball of radius $r$ centered on $\dr \Omega$ such that $\Tr u = 0$ on $2B \cap \dr \Omega$. There holds
\[\int_{B \cap \Omega} |u|^{q-2} |\nabla u|^2 \, dX  \leq \frac{C}{r^2} \int_{(2B \cap \setminus B)\cap \Omega} |u|^q \, dX.\]
Furthermore, if $q>2$, one has
\[\left( \frac{1}{|B \cap \Omega|} \int_{B\cap \Omega}  |u|^q \, dX \right)^{\frac{1}{q}} \leq C  \left(\frac{1}{|2B \cap \Omega|} \int_{2B \cap \Omega}  |u|^2   \, dX \right)^{\frac{1}{2}},\]
and if $q<2$, we have
\[\left(\frac{1}{|B \cap \Omega|} \int_{B\cap \Omega}  |u|^2 \, dX \right)^{\frac{1}{2}} \leq C  \left(\frac{1}{|2B \cap \Omega|} \int_{2B \cap \Omega}  |u|^q   \, dX\right)^{\frac{1}{q}}.\]
The constant $C>0$ depends only on $n$, $q$, a smaller bound on constant $\lambda_q$ in \eqref{defpellipt}, and $\|A\|_\infty$.
\end{prop}

This result, along with a few others that we didn't recall here, can be used to establish the finiteness of the following integrals.

\begin{thm} \label{theoES2}
Retain the assumptions of Proposition~\ref{lemMoserB3}. 
For any solution $u\in W$ to $L u = 0$ whose trace $\Tr u$ is a restriction to $\dr \Omega$ of a function in $C^\infty_0(\R^n)$, we have
\begin{enumerate}[(i)]
\item if $q\geq 2$,
$$\int_\Omega |\nabla u|^2 |u|^{q-2} \, dX < +\infty;$$
\item if $q\in (1,2)$, there exists a ball $B$ centered on $\dr \Omega$ such that 
$$\int_{\Omega\sm B} |\nabla u|^2 |u|^{q-2} \, dX < +\infty.$$
\end{enumerate}
\end{thm}

\subsection{Plan of the article}

Section \ref{spellip} is devoted to the presentation of the $q$-ellipticity.
 In Section \ref{STrace}, we prove the results stated in Subsection \ref{ssES}, which hold in the general context when $\Omega$ is the complement in $\R^n$ of an Ahlfors regular set. Section \ref{SIntroLocal} serves as an introduction to the work with the square function and the non-tangential maximal function, for instance, we establish there the equivalence $\|\wt N_{a,q}(u)\|_p \approx \|\wt N_{1,2}(u)\|_p$ whenever $u$ is a weak solution to $\cL u = 0$, $a>0$,  and $q$ is in the range of ellipticity of $\cL$. 
Sections \ref{sS<N} and \ref{N<S} are devoted to, respectively, the local $S<N$ and $N<S$ bounds, so altogether, these sections contain the proof of Theorem \ref{ThLocal}.
At last, in Section \ref{SExistence} and \ref{sUniqueness}, we prove, respectively, the existence and uniqueness of the solutions to the Dirichlet problem (D$_q$), and the combination of Lemma \ref{lem7.5} (Existence), Lemma \ref{lem8.1} (Equivalence between $S$ and $N$), Lemma \ref{lem8.11} (Improvement \eqref{intro:qntcv}), and Lemma \ref{lem8.4} (Uniqueness) gives Theorem \ref{MainTh}.

\section{The $q$-ellipticity and its consequences}

\label{spellip}
Throughout this section, we assume $\Gamma \subset \R^n$ is an Ahlfors-David regular set of dimension $d<n-1$ and $\Omega := \R^n \setminus \Gamma$. For $x\in \Omega$, define $\delta(x) = \dist(x,\Gamma)$ and we write $dm(x)$ for $\delta^{d+1-n}(x) dx$.

\ms

Consider a matrix $\cA(X)$ with complex coefficients, the usual ellipticity assumption is that there exist constants $\lambda = \lambda_{\cA}>0$ and $\Lambda = \|\cA\|_{\infty} <\infty$ such that for almost every $x\in \Omega$ and every $\xi,\zeta\in \C^n$,
\begin{equation} \label{defellipt}
\lambda |\xi|^2 \leq \Re( \cA(X)\xi \cdot \overline\xi) \quad \text{ and } \quad |\cA(X) \xi \cdot \overline \eta| \leq \Lambda|\xi| |\eta|. 
\end{equation}
A stronger form of ellipticity was introduced in \cite{CD} and \cite{DP}, and see also \cite{CM} where you can find the older - but related - notion of $L^p$ dissipativity. For $q>1$, we say that the matrix $\cA$ is $q$-elliptic if there exists $\lambda_q(\cA)>0$ such that for almost every $X\in \Omega$ and every $\xi,\zeta\in \C^n$,
 \begin{equation} \label{defpellipt}
\lambda_q |\xi|^2 \leq \Re(\cA(X)\xi \cdot \overline{\cJ_q \xi}) \quad \text{ and } \quad |\cA(X) \xi \cdot \overline \eta| \leq \Lambda|\xi||\eta|.
\end{equation}
where $\cJ_q:\C^n \to \C^n$ is defined as
 \begin{equation} \label{defJp}
\cJ_q (\alpha + i\beta ) = \frac{1}q \alpha + \frac i{q'}\beta, \qquad \text{ and } \frac1q + \frac1{q'} = 1.
\end{equation}
Let us make a few simple remarks on $q$-ellipticity. First, a matrix $\cA$ is elliptic in the usual sense (i.e. it satisfies \eqref{defellipt}) if and only if $\cA$ is $2$-elliptic. Second, a matrix $\cA$ can be $q$-elliptic only if $q\in (1,+\infty)$, and moreover $\cA$ is $q$-elliptic if and only if $\cA$ is $q'$ elliptic. In \cite[Proposition 5.17]{CD} (see also the discussions in \cite{DP}) the following nice result can be found.

\begin{prop} \label{prop1.1}
Let $\cA \in L^\infty(\Omega,\C)$ be an elliptic matrix, i.e. a matrix satisfying \eqref{defellipt}. Then $\cA$ is $q$-elliptic if and only if 
 \begin{equation} \label{defmuA}
\mu(\cA) = \essinf_{X\in \Omega} \min_{\xi \in \C^n\setminus\{0\} } \Re \frac{\cA(X) \xi \cdot \overline \xi }{|\cA(X) \xi \cdot  \xi|} > \left| 1- \frac{2}{q} \right|.
\end{equation}
In other words, $\cA$ is $q$-elliptic if and only if $q \in (q_0,q'_0)$, where $q_0 = 2/(1+\mu(\cA))$.
Moreover, we can find $\lambda_q$ satisfying \eqref{defpellipt} such that $C^{-1} \lambda_q \leq \mu(\cA) -|1-2/q| \leq C \lambda_q$, where $C$ depends only on $\mu(\cA)$, $\lambda_\cA$ and $\|\cA\|_\infty$.
\end{prop}

Again, let us make a few comments. The minimum in $\xi$ shall be taken over the $\xi$ such that $|\cA(X) \xi \cdot \xi | \neq 0$. Writing $|\cA(X) \xi \cdot \xi |$ is not a mistake for $|\cA(X) \xi \cdot \overline \xi |$. 
By taking $\xi\in \R^n$, it is easy to check that  $\mu(\cA) \leq 1$ and thus $q_0 \in [1,2]$; and since the ellipticity condition on $\cA$ implies $\mu(\cA) \geq \lambda/\|\cA\|_\infty$, we have $q_0 < 2$. 
If $\cA \in L^\infty(\Omega,\C)$ is elliptic, we also have the following (not completely immediate) equivalence: $\cA$ is real valued if and only if $q_0 = 1$. 
It means that the notion of $q$-ellipticity is not relevant when $\cA$ has real coefficients, and that being complex valued prevents $\cA$ to be $q$-elliptic on the full range of $q\in (1,+\infty)$.

\medskip

The notion of $q$-ellipticity will be used via the following result (whose proof is completely identical to the one of \cite[Theorem 2.4]{DP}).

\begin{prop} \label{proppellip}
Assume that $\cA \in L^\infty(\Omega)$ is a $q$-elliptic matrix. Then there exists $\lambda'_q := \lambda'_q(\lambda,\|\cA\|_\infty,\lambda_q) = \lambda'_q(\lambda_\cA,\|\cA\|_\infty,\mu(\cA),q)>0$ such that for any nonnegative function $\chi \in L^\infty(\Omega)$, and any function $u$ such that $|u|^{q-2}|\nabla u|^2 \chi \in L^2(\Omega,dm)$,
one has
$$\Re \int_\Omega \cA \nabla u \cdot \nabla[|u|^{q-2} \overline u]  \, \chi \, dm \geq \lambda'_q \int_\Omega |u|^{q-2} |\nabla u|^2 \chi \, dm.$$
In particular, the right-hand side above is finite.
\end{prop}

We finish the section with a last observation. We will use repeatedly the following fact (see \cite[Lemma 2.5]{DP}). For any $q>1$, any $u$ such that $v:= |u|^{q/2-1}u \in W^{1,2}_{loc}(\Omega,\C)$, and any $X$ for which $u(X)\neq 0$, we have $\nabla |u| = |u|^{-1} \Re\left( \overline u \nabla u\right)$, and thus
\begin{equation} \label{nablau&v}
C^{-1} |u(X)|^{q-2} |\nabla u(X)|^2 \leq  |\nabla v(X)|^2 \leq C |u(X)|^{q-2} |\nabla u(X)|^2,
\end{equation}
where $C>0$ depends only on $q$.

\section{Moser and Energy estimates}

\label{STrace}

The goal of this section is to prove Moser's estimates and the energy estimates, i.e. the estimates of the gradient of solutions. We will start with interior estimates, then prove boundary estimates for solutions with vanishing or non-vanishing traces. Meanwhile, we will use these estimates to show the a priori finiteness of the square function, that is, we will prove Theorem \ref{theoES2} for energy solutions.
We keep the same assumption on $\Gamma$ as the ones given in Subsection \ref{ssES} and Section \ref{spellip}.

\medskip

\subsection{Interior estimates}

We aim to prove the following result, which easily implies Proposition \ref{lemMoserI}. The notation $\fint_E f \, dm$ is used to denote $m(E)^{-1} \int_E f \, dm$.

\begin{lem} \label{lemMoser}
Let $L= -\div A \nabla$ be an elliptic operator, that is, $L$ satisfies \eqref{defEllip2}--\eqref{defBdd2}. Set $\cA(X) := \delta(X)^{n-d-1} A(X)$ and let $q_0 \in [1,2)$ given by \eqref{defmuA}.
Take $u\in W_{loc}^{1,2}(\Omega)$ a weak solution to $L u = 0$ and a ball $B$ of radius $r$ that satisfies $3B \subset \Omega$. 
\begin{enumerate}[(i)]
\item Let $\Psi$ be a smooth function satisfying $0\leq \Psi \leq 1$ and $|\nabla \Psi(X)| \leq 100/\delta(X)$, and let $k>2$. For $q\in (q_0,q_0')$, there holds
\begin{equation}\label{eq:2gradMoserPhi}
\int_{B} |u|^{q-2} |\nabla u|^2 \Psi^k \, dm  \leq \frac{C}{r^2} \int_{2B} |u|^q \Psi^{k-2} \, dm.
\end{equation}
\item For $q \in (2,\frac n{n-2}q_0')$, we have
\begin{equation}\label{eq:Mosergt2}
\left( \fint_B  |u|^q \, dm \right)^{\frac{1}{q}} \leq C  \left(\fint_{2B}  |u|^2   \, dm \right)^{\frac{1}{2}}.
\end{equation}
\item For $q\in (q_0,2)$, we have
\begin{equation}\label{eq:Moserlt2}
\left( \fint_B  |u|^2 \, dm \right)^{\frac{1}{2}} \leq C  \left(\fint_{2B}  |u|^q   \, dm\right)^{\frac{1}{q}}.
\end{equation}
\end{enumerate}
Each of the above constant $C>0$ depends on $n$, $q$, $\lambda_\cA$, $\|\cA\|_\infty$, $\mu(\cA)$, and depends on $k$ in the case of (i).
\end{lem}

\begin{rmk}
	The above lemma and its proof are inspired by \cite[Lemmas 2.6, 2.7]{DP}. However, our Lemma is stronger than \cite[Lemmas 2.6, 2.7]{DP} in the case $q<2$. In addition, our proof is direct, that is, contrary to the proof in \cite{DP}, we don't approximate $L$ by some elliptic operators $L_j$ with smooth coefficients.
\end{rmk}

\begin{proof}

We set $\Phi = \Psi \eta_B$, where $\eta_B \in C^\infty_0(2 B)$ is a smooth function satisfying $0\leq \eta_B \leq 1$, $\eta_B \equiv 1$ on $B$, and $|\nabla \eta_B| \leq C/r$.

\medskip	

\textbf{Step 1: estimate of the gradient.}
For the purpose of a priori finiteness which will be used later, we also define $u_N = \min\{|u|,N\}$ if $q\geq 2$, and $u_N = \max \{|u|, 1/N\}$ if $q<2$. 
Note that $u_N$ is a real-valued non-negative functions and for all $X\in\Omega$, $u_N$ converges monotonously to $|u|$ as $N\to \infty$. It is also easy to see that $u_N^{\frac{q}{2}-1}|u| \in W_{loc}^{1,2}(\Omega)$, and this guarantees a priori boundedness of the following integrals.

For the case $q\geq 2$, let $E_1 =\{X\in\Omega: |u|\leq N \}$ and $E_2 = \{X\in\Omega: |u| > N\}$. Then 
\[ \int_\Omega u_N^{q-2} |\nabla u|^2 \Phi^k \, dm  = \int_{E_1 } |u|^{q-2}|\nabla u|^2 \Phi^k \, dm + \int_{E_2} N^{q-2} |\nabla u|^2 \Phi^k \, dm. \]
By the $q$-ellipticity, we can apply Proposition \ref{proppellip} to $\chi = \Phi^k \1_{E_1}$ and get
\[ \Re \int_{E_1}  \cA \nabla u \cdot \nabla[|u|^{q-2} \overline u ] \, \Phi^k \, dm \geq \lambda'_q \int_{E_1} |u|^{q-2} |\nabla u|^2 \Phi^k \, dm. \]
Similarly, by the $2$-ellipticity, we have
\[ \Re \int_{E_2} \cA \nabla u \cdot \nabla \overline u \  \Phi^k \, dm \geq \lambda'_2 \int_{E_2}  |\nabla u|^2 \Phi^k \, dm. \]
Therefore
\begin{align}
	\int_\Omega u_N^{q-2} |\nabla u|^2 \Phi^k \, dm & \lesssim \Re \int_\Omega \cA \nabla u \cdot \nabla [u_N^{q-2} \bar u] \Phi^k \, dm \nonumber \\
	& = \Re \int_\Omega \cA \nabla u \cdot \nabla [u_N^{q-2} \bar u \Phi^k] \, dm - \Re \int_\Omega \cA \nabla u \cdot \nabla [\Phi^k] \, u_N^{q-2} \bar u \, dm \nonumber \\
	& := T_1 + T_2. \label{eq:ellpdivg}
\end{align}
A similar argument gives \eqref{eq:ellpdivg} in the case $q<2$.

Observe that $|\nabla u_N| \leq |\nabla u|$, so if $q\geq 2$ we have
\begin{equation} \label{trucmuch5}
 |\nabla [u_N^{q-2} \bar u] | \leq (q-1) u_N^{q-2} |\nabla u| \leq (q-1) N^{q-2} |\nabla u|,
 \end{equation}
and if $q<2$ we have
\begin{equation} \label{trucmuch6}
 |\nabla [u_N^{q-2} \bar u] | \leq (3-q) u_N^{q-2} |\nabla u| \leq (3-q) N^{2-q} |\nabla u|.
  \end{equation}
In any case, $u_N^{q-2} \bar u \in W^{1,2}_{loc}$, hence, since $\Phi$ is compactly supported in $\Omega$, we have that $u_N^{q-2} \bar u \Phi^k$ lies in $W_0$ and is compactly supported in $\Omega$. Lemma 8.16 in \cite{DFMprelim} shows that $u_N^{q-2} \bar u \Phi^k$ is a valid test function for $u \in W^{1,2}_{loc}(\Omega)$ and hence $T_1 = 0$.
As for $T_2$, by H\"older inequality
\begin{align}
|T_2| & \lesssim \int_\Omega |\nabla u| \Phi^{k-1} |\nabla \Phi| \, u_N^{q-2} |u| \, dm \nonumber \\
& \leq \left( \int_\Omega |\nabla u|^2 \Phi^{k} \, u_N^{q-2} \, dm \right)^\frac12 \left( \int_\Omega u_N^{q-2} |u|^2 \Phi^{k-2} |\nabla \Phi|^2 \,  \, dm \right)^\frac12 \label{eq:2estT2}
\end{align}
Combining \eqref{eq:ellpdivg}, $T_1=0$ and \eqref{eq:2estT2}, we conclude 
\begin{equation}\label{eq:2estbfdvd}
	\int_\Omega u_N^{q-2} |\nabla u|^2 \Phi^k \, dm \lesssim \left( \int_\Omega u_N^{q-2} |\nabla u|^2 \Phi^{k}   \, dm \right)^\frac12 \left( \int_\Omega u_N^{q-2} |u|^2 \Phi^{k-2} |\nabla \Phi|^2 \,  \, dm \right)^\frac12. 
\end{equation} 
Note that by the definition of $u_N$,
\[ \int_\Omega u_N^{q-2} |\nabla u|^2 \Phi^k \, dm \leq N^{|q-2|} \int_\Omega |\nabla u|^2 \Phi^k \, dm < \infty, \]
hence we may divide the same term on both sides of \eqref{eq:2estbfdvd} and obtain
\begin{equation}\label{eq:2gradN}
	\int_\Omega u_N^{q-2} |\nabla u|^2 \Phi^k \, dm \lesssim \int_\Omega u_N^{q-2} |u|^2 \Phi^{k-2} |\nabla \Phi|^2 \,  \, dm,
\end{equation}
with a constant depending on $q$. 

Since $u_N \to |u|$ pointwise, by Fatou's lemma and the observation that $u_N^{q-2} \leq |u|^{q-2}$ for both $q\leq 2$ and $q\geq 2$, we have
\begin{align}
	\int_\Omega |u|^{q-2} |\nabla u|^2 \Phi^k \, dm & \leq \liminf_{N\to\infty} \int_\Omega u_N^{q-2} |\nabla u|^2 \Phi^k \, dm \nonumber \\
	 & \lesssim \liminf_{N\to\infty} \int_\Omega u_N^{q-2} |u|^2 \Phi^{k-2} |\nabla \Phi|^2 \, dm \nonumber \\
	& \leq \int_\Omega |u|^q \Phi^{k-2} |\nabla \Phi|^2 \, dm.\label{eq:2gradb}
\end{align}
Now, observe that $\Phi^k \geq \Psi^k \1_B$ and $\Phi^{k-2} \leq \Psi^{k-2} \1_{2B}$. In addition, on $2B$, we have
\begin{equation}\label{est:nablaphi}
	|\nabla \Phi| \leq |\nabla \Psi| + |\nabla \eta_B| \lesssim \delta^{-1} + r^{-1} \lesssim \frac1r
\end{equation}
since $\delta(X) \geq r$ on $2 B$. The bound \eqref{eq:2gradMoserPhi} follows. However, we remark that the estimate \eqref{eq:2gradMoserPhi} could be an empty statement unless we prove that its right hand side is finite. 

\bigskip
\textbf{Step 2: Moser estimate. Case $\boldsymbol{q\geq 2}$.}  

Let $2< k'\in \mathbb{R}$. Since $u_N = \min\{|u|,N\} \leq |u|$, for any $p$ such that $q\leq p \leq qn/(n-2)$, we have
\begin{align}
	\left( \int_{2B} u_N^{p-2} |u|^2 \Phi^{k'} \, dm \right)^{\frac{1}{p}} & \leq \left( \int_{2B} u_N^{p-2\frac{p}{q}} |u|^{2 p/q } \Phi^{k'} \, dm \right)^{\frac{1}{p}} \nonumber \\
	& = \left( \int_{2B} \left( u_N^{\frac{q}{2}-1} |u|\Phi^{k'q/2p} \right)^{2p/q} \, dm \right)^{\frac{1}{p}}.
	\label{eq:2cgpq}
\end{align}
We remark that the first inequality is not true for the case $q<2$.
By the Sobolev-Poincar\'e inequality, see Lemma 4.13 and the following Remark 4.33 of \cite{DFMprelim}, there is a constant $C$ independent of $q, p$ such that
\begin{align}
	\left( \fint_{2B} \left( u_N^{\frac{q}{2}-1} |u|\Phi^{k'q/2p} \right)^{2p/q} \, dm \right)^{\frac{q}{2p}}   \leq C r \left( \fint_{2B} \left|\nabla \left( u_N^{\frac{q}{2}-1} |u|\Phi^{k'q/2p}\right)  \right|^2 dm \right)^{\frac{1}{2}}. \label{eq:2SP}  
\end{align} 
Here we used that the power $2p/q$ is less or equal to the Sobolev exponent $2n/(n-2)$ (and $\infty$ if $n=2$).
Therefore by combining \eqref{eq:2cgpq} and \eqref{eq:2SP}, we get
\begin{align}
	\left( \fint_{2B} u_N^{p-2} |u|^2 \Phi^{k'} \, dm \right)^{\frac{1}{p}} \lesssim  r^{2/q}\left(  \fint_{2B} u_N^{q-2} |\nabla u|^2 \Phi^{k'q/p} \, dm + \fint_{2B} u_N^{q-2} |u|^2 \Phi^{k'q/p-2} |\nabla \Phi|^2 \, dm \right)^{\frac{1}{q}}.\label{eq:2tmp1}
\end{align}
Suppose $k'q/p\geq 2$, we can apply the gradient estimate \eqref{eq:2gradN} and get
\begin{align*}
	\int_{2B} u_N^{q-2} |\nabla u|^2 \Phi^{k'q/p} \, dm \lesssim  \int_{2B} u_N^{q-2} |u|^2 \Phi^{k'q/p-2} |\nabla\Phi|^2 \, dm.
\end{align*}
Combining this with \eqref{eq:2tmp1}, we obtain
\begin{equation}\label{eq:2MoserPhi}
	\left( \fint_{2B} u_N^{p-2} |u|^2 \Phi^{k'} \, dm \right)^{\frac{1}{p}} \lesssim  \left(r^2 \fint_{2B} u_N^{q-2} |u|^2 \Phi^{k'q/p-2} |\nabla \Phi|^2  \, dm \right)^{\frac{1}{q}}. 
\end{equation} 
%
Recall the definition $\Phi= \Psi \eta_B$ and \eqref{est:nablaphi}. We have
\begin{equation}
	\left( \fint_{B} u_N^{p-2} |u|^2 \Psi^{k'} \, dm \right)^{\frac{1}{p}}  \lesssim \left( \fint_{2B} u_N^{p-2} |u|^2 \Phi^{k'} \, dm \right)^{\frac{1}{p}} \lesssim \left( \fint_{2B} u_N^{q-2} |u|^2 \Psi^{k'q/p-2}  \, dm \right)^{\frac{1}{q}}.\label{eq:2MoserpqNPhi}
\end{equation}
In particular, if we choose $\Psi \equiv 1$ and $k'$ big enough (depending only on $n$), the estimate above becomes
\begin{equation}\label{eq:2MoserpqN}
	\left( \fint_{B} u_N^{p-2}  |u|^2 \, dm \right)^{\frac{1}{p}} \lesssim \left( \fint_{2B} u_N^{q-2} |u|^2 \, dm \right)^{\frac{1}{q}}< \infty.
\end{equation} 
whenever $q\in [2,q'_0)$ and $q\leq p \leq qn/(n-2)$. 

\medskip

\textbf{Iteration.} \label{piteration}
By the same argument if we replace $2B$ by $\lambda B$ with $\lambda \in (1,2]$, \eqref{eq:2MoserpqN} also holds, namely
\begin{equation}\label{eq:2MoserpqNl}
	\left( \fint_{B} u_N^{p-2}  |u|^2 \, dm \right)^{\frac{1}{p}} \lesssim \left( \fint_{\lambda B} u_N^{q-2} |u|^2 \, dm \right)^{\frac{1}{q}}.
\end{equation}
with a constant depending on the value of $\lambda$ (as well as $q$, $\|\mathcal{A}\|_\infty$ and the $q$-ellipticity of $\mathcal{A}$).
In particular \eqref{eq:2MoserpqNl} holds for any $q\in [2,q'_0)$ and $p=qn/(n-2)$ (for $n>2$; any $p<\infty$ for $n=2$).
Now let $p \in (2, \frac{n}{n-2}q_0')$ be fixed. Let $\ell\in \mathbb{N}$ be the first integer so that $p\left(\frac{n-2}{n} \right)^\ell \leq 2$, and $\lambda=2^{1/\ell} \in (1,2)$. We can iterate \eqref{eq:2MoserpqNl} $\ell$ times and obtain
\[ \left( \fint_B u_N^{p-2} |u|^2 \, dm \right)^{\frac{1}{p}} \lesssim  \left(\fint_{2 B}  |u|^2   \, dm \right)^{\frac{1}{2}}. \]
Note that the right hand side is finite if $u \in W^{1,2}_{loc}(\Omega)$. Therefore by passing $N\to\infty$ we conclude
\begin{equation}\label{eq:2Moserpgt2}
	\left( \fint_B  |u|^p \, dm \right)^{\frac{1}{p}} \leq C_p  \left(\fint_{2 B}  |u|^2   \, dm \right)^{\frac{1}{2}},
\end{equation} 
for any $ p \in (2, \frac{n}{n-2}q_0')$. Note that in the iterative process, we get estimates with $\ell$ constants depending on the powers $p\frac{n-2}{n}, p \left( \frac{n-2}{n} \right)^2, \cdots, \left( \frac{n-2}{n} \right)^{\ell-1}$, and we combine them into one constant $C_p$ depending on $p$ and $n$. This is Moser's estimate for $p>2$. It also justifies the right-hand side of \eqref{eq:2gradMoserPhi} is finite when $q\geq 2$.

\bigskip
\textbf{Step 3: Moser estimate. Case $\boldsymbol{q\leq 2}$.} 
 By H\"older inequality, if $u\in W^{1,2}_{loc}(\Omega)$, we have $u\in L^q(2B,dm)$ for all $q\leq 2$. Hence we do not need to use $u_N$ to approximate $|u|$ in this case.
Similar to the previous case, for any $q\in (q_0,2]$ and $q\leq p \leq qn/(n-2)$ we have
\begin{align}
	\left( \fint_{B} |u|^p \Psi^{k'} \, dm \right)^{\frac{1}{p}} &  = \left( \fint_{B} \left(  |u|^{q/2} \Psi^{k'q/2p} \right)^{2p/q} \, dm \right)^{\frac{1}{p}} \nonumber \\
	& \leq C  \left( r^2 \fint_{B} \left|\nabla \left( |u|^{q/2} \Psi^{k'q/2p} \right) \right|^2  \, dm \right)^{\frac{1}{q}} \nonumber \\
	& \lesssim  \left( r^2\fint_{B} |u|^{q-2} |\nabla u|^2 \Psi^{k'q/p}  \, dm + r^2 \fint_{B} |u|^{q} \Psi^{k'q/p - 2} |\nabla \Psi|^2 \, dm \right)^{\frac{1}{q}}.
\end{align}
We have $|\nabla \Psi| \lesssim \frac1r$ on $B$. Furthermore, since $u\in L^q(2B, dm)$ for $q\leq 2$, the right-hand side of \eqref{eq:2gradMoserPhi} is finite, thus we can plug \eqref{eq:2gradMoserPhi} in the above estimate. Therefore, we get
\begin{equation} \label{trucmuch2}
\left( \fint_{B} |u|^p \Psi^{k'} \, dm \right)^{\frac{1}{p}} \lesssim  \left( \fint_{2B} |u|^{q} \Psi^{k'q/p - 2} \, dm \right)^{\frac{1}{q}}
\end{equation}
whenever $q\in (q_0,2]$ and $q\leq p \leq qn/(n-2)$. Again, we take $\Psi \equiv 1$ and $k'$ large and we obtain
\begin{equation}
	\left( \fint_{B} |u|^p \, dm \right)^{\frac{1}{p}} \lesssim  \left( \fint_{2B} |u|^{q}  \, dm \right)^{\frac{1}{q}},
\end{equation}
with $q\in (q_0,2]$ and $q\leq p \leq qn/(n-2)$.

\medskip

\textbf{Iteration.}
Apply a similar iterative process as in the previous case, we conclude that
\begin{equation}
	\left( \fint_{B} |u|^2 \, dm \right)^{\frac{1}{2}} \leq C_q  \left( \fint_{2B} |u|^{q}   \, dm \right)^{\frac{1}{q}},
\end{equation}
for any $q\in (q_0,2)$.
\end{proof}

Here is a side product of Lemma \ref{lemMoser}, that will be useful later on.

\begin{lem} \label{lemMoserPsi}
Let $L= -\div A \nabla$ be an elliptic operator, that is, assume that $L$ satisfies \eqref{defEllip2}--\eqref{defBdd2}. Set $\cA(X) := \delta(X)^{n-d-1} A(X)$ and let $q_0 \in [1,2)$ be given by \eqref{defmuA}. Take $k>3$.

For any $q\in (q_0,q_0')$, there exists $\epsilon := \epsilon(n,k,q) >0$ such that if $u\in W_{loc}^{1,2}(\Omega)$ is a weak solution to $L u = 0$, if $B$ is a ball $B$ of radius $r$ that satisfies $3B \subset \Omega$, and if $\Psi$ a smooth function satisfying $0\leq \Psi \leq 1$ and $|\nabla \Psi(X)| \leq 100/\delta(X)$, then we have
\begin{equation}\label{trucmuch3}
\left( \fint_B  |u|^q \Psi^k \, dm \right)^{\frac{1}{q}} \leq C  \left(\fint_{2B}  |u|^{q-\epsilon} \Psi^{k-3}   \, dm \right)^{\frac{1}{q-\epsilon}},
\end{equation}
where the constant $C>0$ depends on $n$, $q$, $\lambda_\cA$, $\|\cA\|_\infty$, $\mu(\cA)$, and $k$.
\end{lem}

\bp Actually, this lemma is almost already proven. We just need to find the constant $\epsilon$.

Indeed, the estimate \eqref{eq:2MoserpqNPhi} gives that (here we switch the roles of $p$ and $q$) for $2\leq p< q'_0$, $p\leq q\leq \frac{n}{n-2}p$,
and $k > 2q/p$,
$$\left( \fint_{B} u_N^{q-2} |u|^2 \Psi^{k} \, dm \right)^{\frac{1}{q}} \lesssim \left( \fint_{2B} u_N^{p-2} |u|^2 \Psi^{kp/q-2} \, dm \right)^{\frac{1}{p}}.$$
Thanks to Lemma \ref{lemMoser}, the right-hand side above is bounded uniformly in $N$, and by taking $N \to +\infty$, we get
$$\left( \fint_{B} |u|^q \Psi^{k} \, dm \right)^{\frac{1}{q}} \lesssim \left( \fint_{2B} |u|^p \Psi^{kp/q-2} \, dm \right)^{\frac{1}{p}}.$$
Given $q\in (2,q'_0)$, we choose $p=q-\epsilon$, where $\epsilon =\frac12 \min\{q-2, 2q/n, q(1-2/k)\} >0$, and we obtain  \eqref{trucmuch3} in the case $q>2$.

As for the case $q\leq 2$, we use \eqref{trucmuch2}, and we have
$$\left( \fint_{B} |u|^q \Psi^{k} \, dm \right)^{\frac{1}{q}} \lesssim  \left( \fint_{2B} |u|^{p} \Psi^{kp/q - 2} \, dm \right)^{\frac{1}{p}}$$
whenever $q\in (q_0,2]$, $p \leq q\leq \frac{n}{n-2}p$ and $k>2q/p$. 
We choose again $p=q-\epsilon$, but where $\epsilon$ is now $\frac12 \min\{q-q_0, 2q/n, q(1-2/k)\} >0$, and we obtain  \eqref{trucmuch3} in the case $q \leq 2$.
\ep

\subsection{Boundary estimates}

We want now to prove Theorem \ref{theoES}. In Theorem \ref{theoES}, we assume a stronger assumption, that is we assume that $u\in W$ is an energy solution. Taking $u\in W$ instead of $u\in W^{1,2}_{loc}(\Omega)$ allows us to define a trace on the boundary, and then to get estimates similar to the ones in Lemma \ref{lemMoser} at the boundary.

We introduce the space $W_0$ defined as
$$W_0 := \{ v \in W, \, \Tr v = 0 \}.$$
Recall that the space $W_0$ is the completion of $C_0^\infty(\Omega)$ under the norm $\|\cdot \|_W$, see Lemma 5.30 of \cite{DFMprelim}.

\begin{lem}[Boundary estimates with vanishing trace] \label{lemMoserB}
Let $L= -\div A \nabla$ be an elliptic operator, that is, assume that $L$ satisfies \eqref{defEllip2}--\eqref{defBdd2}. Set $\cA(X) := \delta(X)^{n-d-1} A(X)$ and let $q_0 \in [1,2)$ be given by \eqref{defmuA}.
Let $u\in W$ be a weak solution to $L u = 0$ such that for a ball $B$ of radius $r$ centered on $\Gamma$, we have $\Tr u = 0$ on $3B$. 
\begin{enumerate}[(i)]
\item For $q\in (q_0,q_0')$, there holds
\begin{equation}\label{eq:2gradMoserPhiB}
\int_{B} |u|^{q-2} |\nabla u|^2 \, dm  \leq \frac{C}{r^2} \int_{2B \setminus B} |u|^q \, dm.
\end{equation}
\item For $q \in (2,\frac n{n-2}q_0')$, we have
\begin{equation}\label{eq:Mosergt2B}
\left( \fint_B  |u|^q \, dm \right)^{\frac{1}{q}} \leq C  \left(\fint_{2B}  |u|^2   \, dm \right)^{\frac{1}{2}}.
\end{equation}
\item For $q\in (q_0,2)$, we have
\begin{equation}\label{eq:Moserlt2B}
\left( \fint_B  |u|^2 \, dm \right)^{\frac{1}{2}} \leq C  \left(\fint_{2B}  |u|^q   \, dm\right)^{\frac{1}{q}}.
\end{equation}
\end{enumerate}
Each of the above constant $C>0$ depends on $n$, $q$, $\lambda_\cA$, $\|\cA\|_\infty$, and $\mu(\cA)$.
\end{lem}

\bp The proof of this lemma is similar to the one of Lemma \ref{lemMoser} and we will only talk about the differences.

\medskip	

\textbf{Step 1: estimate of the gradient.}
Here, we take $\Psi \equiv 1$ and thus $\Phi = \eta_B$, where $\eta_B \in C^\infty_0(2B)$, $0\leq \eta_B \leq 1$, $\eta_B \equiv 1$ on $B$, and $|\nabla \eta_B| \leq C/r$. We take $k=3$ ($k$ has no importance as long as it is bigger than $2$). The function $u_N$ is defined as in Step 1 of the proof of Lemma \ref{lemMoser}.

The proof of \eqref{eq:2gradMoserPhiB} is then done as the one of \eqref{eq:2gradMoserPhi}. The only delicate point is to verify
\begin{equation} \label{trucmuch4}
T_1:= \Re \int_\Omega \cA \nabla u \cdot \nabla [u_N^{q-2} \bar u \Phi^k] \, dm = 0.
\end{equation}
Since $u\in W$ is a solution, according to Lemma 8.16 in \cite{DFMprelim} it suffices to show that $u_N^{q-2} \bar u \Phi^k \in W_0$. The bounds \eqref{trucmuch5}--\eqref{trucmuch6} proves that $u_N^{q-2} \bar u \in W$. The fact that $u_N^{q-2} \bar u \Phi^k$ lies in $W_0$ is then a consequence of $\Tr u \equiv 0$ on $3B$, $\supp \, \Phi \subset 2B$, and Lemmas 5.24 and 6.1 in \cite{DFMprelim}.

\bigskip
\textbf{Step 2: Moser estimate.} The proof of both \eqref{eq:Mosergt2B} and \eqref{eq:Moserlt2B} is identical of the ones of respectively \eqref{eq:Mosergt2} and \eqref{eq:Moserlt2}, and is based on $\eqref{eq:2gradMoserPhiB}$ and the boundary Poincar\'e's inequality given in \cite[Lemma 4.13]{DFMprelim}.
\ep

\medskip

Now we set out to prove Theorem \ref{theoES}.

\begin{lem} \label{lemES1}
Let $L:= -\div A \nabla$ be an elliptic operator. Set $\cA(X) := \delta(X)^{n-d-1} A(X)$, let $q_0 \in [1,2)$ be given by \eqref{defmuA} and $q\in (q_0,q_0')$.

For any energy solution $u\in W$ to $L u = 0$, there exists a ball $B$ centered on $\Gamma$ such that
$$\int_{\Omega \setminus B} |\nabla u|^2 |u|^{q-2} \, dm < +\infty.$$
\end{lem}

\bp If $u\in W$ is an energy solution, there exists $g\in C^\infty_0(\R^n)$ such that $\Tr u = g$ on $\Gamma$. So we can find a ball $B'$ centered on $\Gamma$ such that $\supp \, g \subset B'$, that is, the support of $\Tr u$ is in $ B' \cap \Gamma$. We choose $B = 10B'$, and let $r$ be the radius of $B$.

\medskip

We take $k=3$ and then $\Phi \in C^\infty(\R^n)$ such that $0 \leq \Phi \leq 1$, $\Phi \equiv 1$ outside $B$, $\Phi \equiv 0$ in $\frac12 B$, and $|\nabla \Phi|  \leq C/r$. We can apply the argument in in Step 1 of Lemma \ref{lemMoser} to $\Phi$ and get an estimate similar to \eqref{eq:2gradb}
\begin{equation} \label{trucmuch7}
\int_\Omega |u|^{q-2} |\nabla u|^2 \Phi^k \, dm  \lesssim \int_\Omega |u|^q \Phi^{k-2} |\nabla \Phi|^2 \, dm.
\end{equation}
The only delicate point in the proof of \eqref{trucmuch7} is the proof of the fact that 
\begin{equation} 
\Re \int_\Omega \cA \nabla u \cdot \nabla [u_N^{q-2} \bar u \Phi^k] \, dm = 0
\end{equation}
that can be established with the same reasoning used to prove \eqref{trucmuch4}. By using the properties of $\Phi$, the bound \eqref{trucmuch7} becomes
\begin{equation} \label{trucmuch8}
\int_{\Omega\setminus B} |u|^{q-2} |\nabla u|^2 \, dm  \lesssim \frac1{r^2} \int_{B\setminus \frac12 B} |u|^q \, dm.
\end{equation}
The annulus $B\setminus \frac12 B$ can be covered by a finite number of balls $(D_i)_{i\in I}$ of radius $r/100$ that doesn't intersect $2B'$. Due to the Moser estimates \eqref{eq:Mosergt2} and \eqref{eq:Mosergt2B} (if $q>2$) or simply by H\"older inequality (if $q\leq 2$), we have
$$\left(\fint_{D_i} |u|^q \, dm\right)^\frac1q \lesssim \left(\fint_{D_i} |u|^2 \, dm\right)^\frac12.$$
Therefore, we deduce from \eqref{trucmuch8} that
$$\int_{\Omega\setminus B} |u|^{q-2} |\nabla u|^2 \, dm  \lesssim r^{(1+d)\left(1-\frac{q}{2} \right) - 2}  \left(\int_{2B\setminus 2B'} |u|^2 \, dm\right)^{\frac q2},$$
where $C_B$ depends on the ball $B$ (and in particular its radius $r$), but we don't care. The Poincar\'e inequality implies now, since $\Tr u = 0$ on $2B \setminus 2B'$, 
$$\int_{2B\setminus 2B'} |u|^2 \, dm \lesssim r^2 \int_{2B\setminus 2B'} |\nabla u|^2 \, dm.$$
Therefore
$$\int_{\Omega\setminus B} |u|^{q-2} |\nabla u|^2 \, dm  \lesssim r^{(d-1)\left(1-\frac{q}{2}\right)}  \left(\int_{2B\setminus 2B'} |\nabla u|^2 \, dm\right)^{\frac q2} < +\infty.$$
The lemma follows.
\ep

\medskip

When the trace does not vanish, we can still apply similar argument in Lemma \ref{lemMoser} to $u-g \in W_0$. As long as $q\geq 2$, we can morally speaking bound the integral of $|u|^{q-2}|\nabla u|^2$ by that of $u-g$ and $g$.

\begin{lem}[Boundary estimates with non-vanishing trace: case $q\geq 2$] \label{lemES2}
Let $L:= -\div A \nabla$ be an elliptic operator. Set $\cA(X) := \delta(X)^{n-d-1} A(X)$, let $q_0 \in [1,2)$ given by \eqref{defmuA}, and take $q\in [2,q_0')$, $p\in [2,\frac{n}{n-2}q_0')$.

Choose $g\in C^\infty_0(\R^n)$ and set $u\in W$ to be the (unique) energy solution to $L u = 0$ satisfying $\Tr u = g$.  For any ball $B$ of radius $r$ centered on $\Gamma$, we have
\begin{equation} \label{Moser1}
\int_{B} |\nabla u|^2 |u|^{q-2} dm \lesssim \frac1{r^2} \int_{2B} |u|^q dm + r^{q-2} \int_{2B} |\nabla g|^q + r^{d-1} \|g\|_{L^\infty(2B)}^{q},  
\end{equation} 
and
\begin{equation} \label{Moser2}
\left(\fint_{B} |u|^p dm\right)^{\frac1p} \lesssim \left( \fint_{2B} |u|^2 dm \right)^\frac12 + \left( \fint_{2B} |g|^p dm \right)^\frac1p + r \left( \fint_{2B} |\nabla g|^{p_*} dm \right)^\frac1{p_*},
\end{equation}
where  $p_* = \max\{2, \frac{n-2}{n}p\}$, and where the constants depends only on $n$, $\lambda_\cA$, $\mu(\cA)$, $\|\cA\|_\infty$, and, respectively, $q$ and $p$.

As a consequence, for any energy solution $u\in W$ to $L u = 0$ and any ball $B \subset \R^n$ centered on $\Gamma$, one has
\begin{equation} \label{Moser3}
\int_{B} |\nabla u|^2 |u|^{q-2} \, dm < +\infty.
\end{equation}
\end{lem}

\bp We define $v=u-g\in W_0$. The proof follows the same arguments as the ones given in Lemma \ref{lemMoser}, step 1 and 2 or Lemma \ref{lemMoserB}, the only difference being that here we have $Lv=-Lg$ instead.
\bigskip

\textbf{Step 1: estimate of the gradient.}
We set $k$ large - say $k = 100$ - and then $\Phi \in C^\infty_0(2B)$, $0\leq \Phi \leq 1$, $\Phi \equiv 1$ on $B$, and $|\nabla \Phi | \leq C/r$. We also define  $v_N = \min\{|v|,N\}$ to ensure the a priori finiteness of the integrals we work with.

Using $q$-ellipticity and $2$-ellipticity, we obtain similarly to \eqref{eq:ellpdivg} that
\begin{align}
\int_\Omega v_N^{q-2} |\nabla v|^2 \Phi^k \, dm & 
\lesssim \Re \int_\Omega \cA \nabla v \cdot \nabla [v_N^{q-2} \bar v] \Phi^k \, dm \nonumber \\
& = \Re \int_\Omega \cA \nabla v \cdot \nabla [v_N^{q-2} \bar v \Phi^k] \, dm - \Re \int_\Omega \cA \nabla v \cdot \nabla [\Phi^k] \, v_N^{q-2} \bar v \, dm \nonumber \\
& := T_1 + T_2. \label{MoserB1}
\end{align}
Since $q\geq 2$ we have
$$ |\nabla [v_N^{q-2} \bar v] | \leq (q-1) v_N^{q-2} |\nabla v| \leq (q-1) N^{q-2} |\nabla v|, $$
which ensures that $v_N^{q-2} \bar v \in W$. Moreover, since $\Tr v = 0$, Lemma 6.1 in \cite{DFMprelim} gives that $\Tr [v_N^{q-2} \bar v ] = 0$ and then Lemma 5.24 gives that $\varphi := v_N^{q-2} \bar v \Phi^k \in W_0$. So the term $T_1$ is $v$ tested against the function $\varphi \in W_0$. Since $v = u-g$ and $u$ is a solution to $L u = 0$, we deduce
\[\begin{split}
|T_1| & = \left|\int_\Omega \cA \nabla g \cdot \nabla [v_N^{q-2} \bar v \Phi^k] \, dm\right| \lesssim \int_\Omega |\nabla g| |\nabla [v_N^{q-2} \bar v \Phi^k]| \, dm. \\
& \lesssim \int_\Omega  |\nabla g| |\nabla v| v_N^{q-2} \Phi^k \, dm + \int_\Omega |\nabla g| |\nabla \Phi| v_N^{q-2} |v| \Phi^{k-1}  \, dm \\
& \lesssim \left[ \left( \int_\Omega v_N^{q-2} |\nabla v|^2 \Phi^k \, dm  \right)^\frac12 + \left(\int_\Omega v_N^{q-2} |v|^2 \Phi^{k-2} |\nabla \Phi|^2 \, dm  \right)^\frac12 \right] \left(\int_\Omega |\nabla g|^2 v_N^{q-2} \Phi^k \, dm  \right)^\frac12.
\end{split}\]
The last term in the last inequality can be treated as follows: using the fact that $a^{\theta} b^{1-\theta} \lesssim a + b$, where in our case $\theta = 1-\frac2q$, $a = r^{-2} v_{N}^{q}$ and $b = r^{q-2} |\nabla g|^q$, we have
\begin{equation} \label{interpolationforug}\begin{split}
\int_{\Omega} v_{N}^{q-2} |\nabla g|^2 \Phi^k\, dm & \lesssim \frac1{r^2} \int_{\Omega } v_{N}^{q}  \Phi^k\, dm + r^{q-2} \int_{\Omega} |\nabla g|^q \Phi^k \, dm \\
& \lesssim \frac1{r^2} \int_{\Omega } v_{N}^{q-2} |v|^2  \Phi^k\, dm + r^{q-2} \int_{\Omega} |\nabla g|^q \Phi^k \, dm.
\end{split}\end{equation}
Therefore, together with the fact that $|\nabla \Phi| \leq \frac{C}{r}$, the bound on $T_1$ becomes
\begin{equation} \label{boundT'} \begin{split}
|T_1| & \lesssim \left( \int_\Omega v_N^{q-2} |\nabla v|^2 \Phi^k \, dm  \right)^\frac12 \left[ \frac1{r^2} \int_{\Omega } v_{N}^{q-2} |v|^2  \Phi^{k-2}\, dm + r^{q-2} \int_{\Omega} |\nabla g|^q \Phi^k \, dm \right]^\frac12 \\
& \hspace{5cm} +  \frac1{r^2} \int_{\Omega } v_{N}^{q-2} |v|^2  \Phi^{k-2}\, dm + r^{q-2} \int_{\Omega} |\nabla g|^q \Phi^k \, dm.
\end{split}\end{equation}

We turn to the estimate of $T_2$. One has by Cauchy-Schwarz's inequality,
\begin{equation} \label{boundT2}
\begin{split}
|T_2| & \lesssim \int_\Omega |\nabla v| v_N^{q-2} |v| |\nabla \Phi| \Phi^{k-1} \\
& \lesssim \left( \int_\Omega v_N^{q-2} |\nabla v|^2 \Phi^k \, dm  \right)^\frac12 \left( \int_\Omega v_N^{q-2} |v|^2 \Phi^{k-2} |\nabla \Phi|^2 \, dm  \right)^\frac12 \\
& \lesssim \left( \int_\Omega v_N^{q-2} |\nabla v|^2 \Phi^k \, dm  \right)^\frac12 \left( \frac1{r^2} \int_\Omega v_N^{q-2} |v|^2 \Phi^{k-2} \, dm  \right)^\frac12
\end{split}\end{equation}
where we use the fact that $|\nabla \Phi| \lesssim r^{-1}$ in the last line. 
The combination of \eqref{MoserB1}, \eqref{boundT'}, and \eqref{boundT2} implies that
\[\begin{split}
\int_\Omega v_N^{q-2} |\nabla v|^2 \Phi^k \, dm & \lesssim \left( \int_\Omega v_N^{q-2} |\nabla v|^2 \Phi^k \, dm  \right)^\frac12 \left[ \frac1{r^2} \int_{\Omega } v_{N}^{q-2} |v|^2  \Phi^{k-2}\, dm + r^{q-2} \int_{\Omega} |\nabla g|^q \Phi^k \, dm \right]^\frac12 \\
& \hspace{5cm} +  \frac1{r^2} \int_{\Omega } v_{N}^{q-2} |v|^2  \Phi^{k-2}\, dm + r^{q-2} \int_{\Omega} |\nabla g|^q \Phi^k \, dm,
\end{split}\]
which self-improves, since $\int_\Omega v_N^{q-2} |\nabla v|^2 \Phi^k \, dm $ is finite, to
\begin{equation} \label{MoserB2}
\int_\Omega v_N^{q-2} |\nabla v|^2 \Phi^k \, dm \lesssim \frac1{r^2} \int_{\Omega } v_{N}^{q-2} |v|^2  \Phi^{k-2}\, dm + r^{q-2} \int_{\Omega} |\nabla g|^q \Phi^k \, dm.
\end{equation}
Recall that $\Phi \equiv 1$ on $B$ and $\Phi \equiv 0$ outside $2B$, hence
\begin{equation} \label{MoserB3}
\int_B v_N^{q-2} |\nabla v|^2 \, dm \lesssim \frac1{r^2} \int_{2B} v_{N}^{q-2} |v|^2 \, dm + r^{q-2} \int_{2B} |\nabla g|^q  \, dm.
\end{equation}

\bigskip

\textbf{Step 2: Moser estimate.} Using similar arguments as the ones invoked in Step 2 of the proof of Lemma \ref{lemMoser} (mainly based on the Poincar\'e inequality), we obtain an analogue of \eqref{eq:2tmp1}, that is
$$\left(\fint_{2B} v_N^{p-2} |v|^2 \Phi^{k'} \, dm \right)^{\frac{1}{p}} \lesssim  r^{2/q}\left(  \fint_{2B} v_N^{q-2} |\nabla v|^2 \Phi^{k'q/p} \, dm + \fint_{2B} v_N^{q-2} |v|^2 \Phi^{k'q/p-2} |\nabla \Phi|^2 \, dm \right)^{\frac{1}{q}}$$
whenever $2\leq q\leq p\leq qn/(n-2)$. Assuming that $k'$ is such that $k'q/p-2 >0$ and using \eqref{MoserB2} with $k = k'q/p$ for the first term in the right-hand side, we obtain that
\[\begin{split}
\left(\fint_{2B} v_N^{p-2} |v|^2 \Phi^{k'} \, dm \right)^{\frac{1}{p}} & \lesssim  r^{2/q}\left(  \frac1{r^2} \fint_{2B} v_N^{q-2} |v|^2 \Phi^{k'q/p-2} \, dm + r^{q-2} \fint_{2B} |\nabla g|^q \Phi^{k'q/p} \, dm \right.\\
& \hspace{5cm} \left. + \fint_{2B} v_N^{q-2} |v|^2 \Phi^{k'q/p-2} |\nabla \Phi|^2 \, dm \right)^{\frac{1}{q}},
\end{split}\]
whenever $2\leq q< q'_0$.
Now, use the fact that $ |\nabla \Phi| \leq C/r$, $\Psi \equiv 1$ on $B$, and $\Psi \leq 1$, to get
$$\left(\fint_{B} v_N^{p-2} |v|^2\, dm \right)^{\frac{1}{p}} \lesssim  r^{2/q}\left(  \frac1{r^2} \fint_{2B} v_N^{q-2} |v|^2 \, dm + r^{q-2} \fint_{2B} |\nabla g|^q \, dm \right)^{\frac{1}{q}},$$
or
$$\left(\fint_{B} v_N^{p-2} |v|^2\, dm \right)^{\frac{1}{p}} \lesssim  \left(\fint_{2B} v_N^{q-2} |v|^2 \, dm\right)^{\frac{1}{q}} + r \left( \fint_{2B} |\nabla g|^q \, dm \right)^{\frac{1}{q}},$$
The arguments can be slightly modified to get
$$\left(\fint_{B} v_N^{p-2} |v|^2\, dm \right)^{\frac{1}{p}} \leq C_\lambda  \left(\fint_{\lambda B} v_N^{q-2} |v|^2 \, dm\right)^{\frac{1}{q}} + C_{\lambda} r \left( \fint_{\lambda B} |\nabla g|^q \, dm \right)^{\frac{1}{q}},$$
where $\lambda \in (1,2)$ and $C_\lambda$ depends on $\lambda$. We use then an iterative argument as the one given page \pageref{piteration} and we get, if $p\in [2,\frac{n}{n-2}q_0')$ and $p_* = \max\{2,\frac{n-2}{n} p\}$,
\begin{equation} \label{MoserB4}
\left(\fint_{B} v_N^{p-2} |v|^2\, dm \right)^{\frac{1}{p}} \lesssim  \left(\fint_{2B} |v|^2 \, dm\right)^{\frac{1}{2}} + r \left( \fint_{2B} |\nabla g|^{p_*} \, dm \right)^{\frac{1}{p_*}}.
\end{equation} 

\bigskip

\textbf{Step 3: Conclusion.} The estimate \eqref{MoserB4} gives a uniform (and finite) bound on 
$$\left(\fint_{B} v_N^{p-2} |v|^2\, dm \right)^{\frac{1}{p}}$$
so by taking $N\to +\infty$, we obtain, if $p\in [2,\frac{n}{n-2}q_0')$ and $p_* = \max\{2,\frac{n-2}{n} p\}$,
\begin{equation} \label{MoserB5}
\left(\fint_{B} |v|^p\, dm \right)^{\frac{1}{p}} \lesssim  \left(\fint_{2B} |v|^2 \, dm\right)^{\frac{1}{2}} + r \left( \fint_{2B} |\nabla g|^{p_*} \, dm \right)^{\frac{1}{p_*}}.
\end{equation} 
The inequality above is not \eqref{Moser2}, because the estimate is on $v = u-g$ and not on $u$. But by triangle inequality and H\"older inequality
$$\left(\fint_{B} |u|^p\, dm \right)^{\frac{1}{p}} \leq \left(\fint_{B} |v|^p\, dm \right)^{\frac{1}{p}} + \left(\fint_{B} |g|^p\, dm \right)^{\frac{1}{p}}$$
and
$$\left(\fint_{2B} |v|^2\, dm \right)^{\frac{1}{2}} \leq \left(\fint_{2B} |u|^2\, dm \right)^{\frac{1}{2}} + \left(\fint_{2B} |g|^p\, dm \right)^{\frac{1}{p}},$$
we easily obtain the desired estimate \eqref{Moser2}.

\medskip

Thanks to \eqref{MoserB5}, the right-hand side of \eqref{MoserB3} is uniformly bounded in $N$. Thus we may take $N\to+\infty$ in \eqref{MoserB3} and get
\begin{equation} \label{MoserB6}
\int_B |v|^{q-2} |\nabla v|^2 \, dm \lesssim \frac1{r^2} \int_{2B} |v|^q \, dm + r^{q-2} \int_{2B} |\nabla g|^q  \, dm.
\end{equation}
Again, we want to turn \eqref{MoserB6} into an estimate on $u = v+g$. Yet, observe that
\begin{equation} \label{MoserB7}
\frac1{r^2} \int_{2B} |v|^q\, dm  \lesssim \frac1{r^2} \int_{2B} |u|^q\, dm + r^{d-1} \|g\|_{L^\infty(2B)}^q.
\end{equation}
Besides, we have
\[\begin{split}
\int_B |u|^{q-2} |\nabla u|^2 \, dm & \lesssim \int_B |v|^{q-2} |\nabla v|^2 \, dm + \int_B |v|^{q-2} |\nabla g|^2 \, dm \\
& \hspace{3cm} + \int_B |g|^{q-2} |\nabla v|^2 \, dm + \int_B |g|^{q-2} |\nabla g|^2 \, dm \\
& := I_1 + I_2 + I_3 + I_4.
\end{split}\]
We don't change $I_1$. We use the fact that $a^{q-2} b^2 \lesssim r^{d-1} a^q + r^{(d-1)(1-q/2)} b^q$ and we get
\[ \begin{split}
I_3 & \leq \|g\|_{L^\infty(B)}^{q-2} \int_{B} |\nabla v|^2 dm \lesssim r^{d-1} \|g\|_{L^\infty(B)}^q + r^{(d-1)(1-q/2)} \left( \int_{B} |\nabla v|^2 dm \right)^\frac q2 \\
& \lesssim r^{d-1} \|g\|_{L^\infty(B)}^q + r^{(d-1)(1-q/2)} \left( \frac{1}{r^2} \int_{2B} |v|^2 dm + \int_{2B}|\nabla g|^2 dm \right)^{\frac{q}{2}} \\ 
& \lesssim   r^{d-1} \|g\|_{L^\infty(B)}^q +  \frac1{r^2} \int_{2B} |v|^q dm + r^{q-2} \int_{2B}|\nabla g|^q dm. \\
\end{split}\]
Similarly, since $a^{q-2} b^2 \lesssim r^{-2} a^q + r^{q-2} b^q$, we have
$$I_2 \lesssim \frac1{r^2} \int_B |v|^q \, dm + r^{q-2} \int_B |\nabla g|^q \, dm$$
and
$$I_4 \lesssim \frac1{r^2} \int_B |g|^q \, dm + r^{q-2} \int_B |\nabla g|^q \, dm \lesssim r^{d-1} \|g\|_{L^\infty(B)}^q + r^{q-2} \int_B |\nabla g|^q \, dm.$$
Altogether, we deduce
\begin{equation} \label{MoserB8}
\int_B |u|^{q-2} |\nabla u|^2 \, dm  \lesssim  \int_B |v|^{q-2} |\nabla v|^2 \, dm +  \frac1{r^2} \int_{B} |v|^q dm + r^{d-1} \|g\|_{L^\infty(B)}^q + r^{q-2} \int_B |\nabla g|^q \, dm.
\end{equation}
The combination of \eqref{MoserB6}, \eqref{MoserB7}, and \eqref{MoserB8} gives \eqref{Moser1}.

\medskip

Eventually, by combining \eqref{Moser1} and \eqref{Moser2}, we obtain
\[\begin{split}
\int_{B} |\nabla u|^2 |u|^{q-2} \, dm & \lesssim \frac1{r^2} \int_{2B} |u|^q dm + r^{q-2} \int_{2B} |\nabla g|^q + r^{d-1} \|g\|_{L^\infty(2B)}^{q} \\
& \lesssim r^{(d+1)(1-q/2)-2} \left(\int_{4B} |u|^2 dm\right)^\frac{q}2 + r^{q-2} \int_{4B} |\nabla g|^q + r^{d-1} \|g\|_{L^\infty(2B)}^{q} \\
& < +\infty
\end{split}\]
because $g\in C^\infty_0(\R^n)$ and $u\in W$ (recall that the latter forces $u \in L^2_{loc}(\R^n,dm)$, see Lemma 3.2 in \cite{DFMprelim}).
The lemma follows.
\ep

Combining the Lemmas \ref{lemES1} and \ref{lemES2}, we easily prove Theorem \ref{theoES}. 

\section{Preliminaries to the local estimates.}

\label{SIntroLocal} 

In this section, and for the rest of the article, we take $\Omega = \R^n \sm \R^d = \{(x,t) \in \R^n, \, x\in \R^d \text{ and } t\in \R^{n-d}\sm \{0\}\}$; the boundary $\Gamma$ is assimilated to $\R^d$. We write $X = (x,t)$ or $Y = (y,s)$ for the running points in $\Omega$, and we write $B_l(x)$ for the ball in $\R^d$ with center $x$ and radius $l$. In this particular case, the measure $m$ satisfies 
$$dm(X) = \frac{dt}{|t|^{n-d-1}} \, dx.$$

\subsection{Properties of the non-tangential maximal function $N$}

\begin{lem} \label{lemMoser22}Let $a,b>0$, $p\in [1,+\infty]$, and $\cL$ be an elliptic operator. Set $q_0 \in (1,2)$ given by Proposition \ref{prop1.1} and choose  $q_1,q_2\in (q_0,q_0')$. Then for any weak solution $u\in W^{1,2}_{loc}(\Omega)$ to $\cL u = 0$ that satisfies $\|\wt N_{a,q_1}(u)\|_p < +\infty$, we have
$$\|\wt N_{b,q_2}(u)\|_p \leq C \|\wt N_{a,q_1}(u)\|_p < +\infty,$$
where $C>0$ depends only on $n$, $a$, $b$, $q_1$, $q_2$, $\lambda_\cA$, $\|\cA\|_\infty$ and $\mu(A)$.
\end{lem}

\bp We use the notation introduced in Subsection \ref{SSMR}. The proof of Lemma \ref{lemMoser} can be easily adapted to get \eqref{eq:Mosergt2}--\eqref{eq:Moserlt2} where the ball $B$ is replaced by $W_b(z,r)$ and $2B$ is replaced by 
$$\{(y,s) \in \Omega, \, y \in B_{br}(z), \, r/4 \leq |s| \leq 4r\} \subset W_{4b}(z,r/2) \cup W_{4b}(z,2r).$$
So we obtain, for any $(z,r)\in \R^{d+1}_+$, that
\begin{equation} \label{Nprop1}
\left( \fint_{W_b(z,r)} |u|^{q_2} \, dm \right)^\frac1{q_2} \lesssim \left( \fint_{W_{4b}(z,r/2) \cup W_{4b}(z,2r)} |u|^{q_1} \, dm \right)^\frac1{q_1}.
\end{equation}
In addition, we have 
$$|W_{4b}(z,r/2)| \simeq b^{d}r^{n} \simeq |W_{4b}(z,2r)|$$
and 
$$dm(X) \simeq r^{d+1-n} dx \, dt \qquad \text{ for any $X= (x,t)$ in $W_{4b}(z,r/2) \cup W_{4b}(z,2r)$}.$$
We deduce that \eqref{Nprop1} becomes
$$u_{\W,b,q_2}(z,r) \lesssim u_{\W,4b,q_1}(z,r/2) + u_{\W,4b,q_1}(z,2r).$$
By taking the supremum on $(z,r) \in \Gamma_b(x)$, where $x\in \R^d$, we get the pointwise bound
$$\wt N_{b,q'}(u)(x) \leq C \wt N_{4b,q}(u)(x)$$
with a constant $C$ independent of $x$. We take the $L^p$-norm of the inequality above to obtain
$$\|\wt N_{b,q_2}(u)\|_p \lesssim \|\wt N_{4b,q_1}(u)\|_p \lesssim \|\wt N_{a,q_1}(u)\|_p,$$
where the last inequality is obtained by a classical real variable method argument (See Chapter II, equation (25) in \cite{Stein93}).
\ep

We also have the following corollary of Lemma \ref{lemMoserPsi}, that will be useful in Section \ref{N<S}.

\begin{lem} \label{lemMoserPsi2} Let $a>0$ and $\cL:= -\div |t|^{d+1-n}\cA \nabla$ be an elliptic operator. Set $q_0 \in (1,2)$ given by Proposition \ref{prop1.1} and choose $q\in (q_0,q_0')$. Take also $k>3$.

There exists $\epsilon:= \epsilon(n,k,q,q_0)$ such that for any weak solution $u\in W^{1,2}_{loc}(\Omega)$ to $\cL u = 0$, and for any smooth  function $\Psi \in C^\infty_0(\R^n)$ that satisfies $0\leq \Psi \leq 1$ and $|\nabla \Psi| \leq 100/|t|$, then we have
$$\|\wt N_{a,q}(u|\Psi^k)\|_q \leq C \|\wt N_{a,q-\epsilon}(u|\Psi^{k-3})\|_q,$$
where $C>0$ depends only on $n$, $a$, $q$, $\lambda_\cA$, $\|\cA\|_\infty$, $\mu(A)$ and $k$. In particular, if the right-hand side above is finite, so is the left-hand side.
\end{lem}

\bp As in the previous lemma, the use of Lemma \ref{lemMoserPsi} with sets $W_a(z,r)$ instead of balls give that
$$(u|\Psi^k)_{\W,a,q}(z,r) \lesssim (u|\Psi^{k-3})_{\W,4a,q-\epsilon}(z,r/2) +  (u|\Psi^{k-3})_{\W,4a,q-\epsilon}(z,2r).$$
By taking the supremum on the cones $\Gamma_a(x)$, and then the $L^q$ norm on $x\in \R^d$, we get
$$\|\wt N_{a,q}(u|\Psi^k)\|_q \lesssim \|\wt N_{4a,q-\epsilon}(u|\Psi^{k-3})\|_q.$$
We can upgrade $\wt N_{4a,q-\epsilon}$ into $\wt N_{a,q-\epsilon}$ by using, as in the proof of Lemma \ref{lemMoser22}, a real variable argument. The lemma follows.
\ep

\subsection{The Carleson inequality and good cut-off functions.}

\begin{prop}[Carleson inequality] \label{propCarl} Let $a>0$ and $q\in (1,+\infty)$. Let $f$ be a measurable function (scalar, vector valued, or matrix valued) that satisfies the Carleson measure condition (see Definition \ref{defCMI}). For any function $u \in L^q_{loc}(\Omega)$ and any non-negative function $\Psi \in C_0^\infty(\Omega)$, we have 
$$\int_{(x,t)\in \Omega} |f|^2 |u|^q \Psi \, \frac{dt}{|t|^{n-d}} \, dx \leq C \|f\|_{CM,a} \|\wt N_{a,q}(u|\Psi)\|_q^q,$$
where the constant $C$ depends only on the dimensions $n$ and $d$.
\end{prop}

\bp Observe that
\[\begin{split}
\iint_{(x,t)\in \Omega} |f|^2 |u|^q \Psi \, \frac{dt}{|t|^{n-d}} \, dx 
& \lesssim \iint_{(x,r)\in \R^{d+1}_+} \frac1{|W_a(x,r)|}\iint_{(y,s) \in W_a(x,r)} |f|^2 |u|^q \Psi \, ds \, dy \, \frac{dr}{r} \, dx \\
& \lesssim \iint_{(x,r)\in \R^{d+1}_+} (u|\Psi)_{\W,a,q}^q \left(\sup_{(y,s) \in W_a(x,r)} |f|^2 \right) \, \frac{dr}{r} \, dx \\
& \lesssim \|f\|_{CM,a} \|\wt N_{a,q}(u|\Psi)\|_q^q,
\end{split}\]
by the classical Carleson inequality, see for instance \cite[Chapter II, Section 2.2]{Stein93}.
\ep

\bigskip

Let us define `good' cut-off functions.

\begin{defi} \label{defH2}
Let $a>0$. A function $\Psi$ satisfies $(\cH^2_a)$ if $\Psi$ is locally Lipschitz on $\Omega$ and there exists $M>0$ such that
\begin{enumerate}[(i)]
\item $0\leq \Psi \leq 1$ on $\Omega$,
\item $|\nabla \Psi(x,t)| \leq \frac{C}{|t|}$ for any $(x,t) \in \Omega$,
\item for any $x\in \R^d$, 
\begin{equation}\label{def:Ha2}
\iint_{(z,r) \in \Gamma_b(x)} \left( \sup_{W_b(z,r)} |\nabla \Psi| \right) \frac{dr}{r^{d}} \, dz \leq M,
\end{equation}
where $b = a/100$.
\end{enumerate}
We say that $\Psi$ satisfies $(\cH^2_{a,M})$ if $\Psi$ satisfies (iii) with the constant $M$.
\end{defi}

It appears that the condition $(\cH^2_{a,M})$ depends on the constant $C$ in $(ii)$, but we will cheat and not talk about it. It is fairly reasonable because condition $(iii)$ implies $(ii)$, with a constant $C$ depending on $b, M, d$, and condition $(ii)$ is only here to simplify the reading.
Indeed, a point $(x,t)\in \Omega$ belongs to $W_b(z,r)$ for any $z\in \R^d$, $r>0$ which satisfies $|t|/2\leq r \leq 2|t|$ and $z\in B_{b|t|/4}(x)$, so 
\[\begin{split}
|\nabla \Psi(x,t)| & \lesssim \frac1{|t|} \int_{B_{b|t|/4}} \int_{|t|/2 \leq r \leq 2|t|} \sup_{W_b(z,r)} |\nabla \Psi| \frac{dr}{r^{d}} \, dz \leq  \frac1{|t|} \iint_{(z,r) \in \Gamma_b(x)} \left( \sup_{W_b(z,r)} |\nabla \Psi| \right) \frac{dr}{r^{d}} \, dz \\
&\lesssim \frac{M}{|t|},
\end{split}\]
if $(iii)$ is satisfied.

\ms

The following observations are crucial. 

\begin{lem} \label{exgoodcutoff}
Let $\phi \in C^\infty_0(\R_+)$ be such that $0\leq \phi \leq 1$, $\phi \equiv 1$ on $[0,1]$, $\phi \equiv 0$ on $[2,+\infty)$, and $|\phi'| \leq 2$.
\begin{enumerate}[(1)]
\item If $e(x)$ is a positive $a^{-1}$-Lipschitz function and
$$\Psi_e(x,t) = \phi\left( \frac{e(x)}{|t|} \right),$$
then $\Psi_e$ satisfies $(\cH^2_{a,M})$ for an $M$ that depends only on $a$ and the dimensions $n$, $d$.
\item If $B \subset \R^d$ is a (boundary) ball of radius bigger than or equal to $l$, and
$$\Psi_{B,l}(x,t) = \phi\left( \frac{a|t|}{l}\right) \phi\left( 1+ \frac{\dist(x,B)}{l}\right),$$
then $\Psi_{B,l}$ satisfies $(\cH^2_{a,M})$ for an $M$ that depends only on $a$, $n$.
\item If $\Psi_1$ and $\Psi_2$ are two functions satisfying $(\cH^2_{a,M_1})$ and $(\cH^2_{a,M_2})$ respectively, then $\Psi_1\Psi_2$ satisfies $(\cH^2_{a,M_1+M2})$.
\end{enumerate}
\end{lem}

\bp Part (3) follows easily from the definition \ref{defH2}. We shall only prove (1), and the proof for (2) is similar. The assertion $(i)$ of Definition \ref{defH2} is trivial, and as we said, $(ii)$ follows from $(iii)$, so we only need to verify $(iii)$.

Let $b=a/100$ and take $x\in \R^d$. We want to show that
$$\iint_{(z,r) \in \Gamma_b(x)} \left( \sup_{W_b(z,r)} |\nabla \Psi_e| \right) \frac{dr}{r^{d}} \, dz \leq M.$$

By the properties of $\phi$, we first observe that $|\nabla\Psi_e(y,s)| \leq C/|s|$, and thus $\sup_{W_b(z,r)} |\nabla\Psi_e| \leq C'/r$. Another simple observation is that
\[ \{(y,s)\in \R^n, |\nabla \Psi_e(y,s)| \neq 0\} \subset \left\{(y,s)\in\R^n, e(y)/2 \leq |s| \leq e(y)\right\}. \]

If $(y,s) \in W_b(z,r)$ for some $(z,r) \in \Gamma_b(x)$, then necessarily $(y,s) \in \wh \Gamma_{3b}(x)$. So we want to find for which values of $(y,s)\in \wh \Gamma_{3b}(x)$, we have $ |\nabla \Psi_e(y,s)| \neq 0$. If $y = x$, then we have 
$$S_x:= \{s\in \R^{n-d}, \,  |\nabla \Psi_e(x,s)| \neq 0 \} \subset \{s\in \R^{n-d}, \, e(x)/2 \leq |s| \leq e(x)\}.$$
Besides, if $(y,s) \in \wh \Gamma_{3b}(x) $ is such that $|\nabla \Psi_e(y,s)| \neq 0$, 
since $e$ is $a^{-1}$-Lipschitz, we can find $s_x\in S_x$ such that
$$|s-s_x| \leq \frac{1}{a} |y-x| \leq \frac{3b}{a} |s| \leq \frac{|s|}{2}.$$
We deduce that 
$$\frac{e(x)}{3} \leq \frac23 |s_x| \leq |s| \leq 2s_x \leq 2e(x),$$
and then, the values $(z,r)\in \Gamma_b(x)$ for which $W_b(z,r)$ contains such $(y,s)$ satisfy
\[ \frac{e(x)}{6} \leq \frac{|s|}{2} \leq r\leq 2|s| \leq 4e(x). \]
Therefore,
\[\begin{split}
\iint_{(z,r) \in \Gamma_b(x)} \left( \sup_{W_b(z,r)} |\nabla \Psi_e| \right) \frac{dr}{r^{d}} \, dz 
& \lesssim \int_{e(x)/6 \leq r \leq 4e(x)} \int_{B_{br}(x)} dz\, \frac{dr}{r^{d+1}} \lesssim 1.
\end{split}\]
Part (1) of the lemma follows.
\ep

\ms

Now, let us state another important property of functions which satisfy $(\cH^2_{a})$.

\begin{lem} \label{lemH2}
Let $a>0$ and $q\in (1,+\infty)$. Let $\Psi$  satisfy $(\cH^2_{a,M})$. If $v\in L^q_{loc}(\Omega)$ and $0 \leq \Phi \in C^\infty_0(\Omega)$, then
\begin{equation} \label{H2la}
\iint_\Omega |v|^q \Phi |\nabla \Psi|^2 \, \frac{dt}{|t|^{n-d-2}} \, dx + \iint_\Omega |v|^q \Phi |\nabla \Psi| \, \frac{dt}{|t|^{n-d-1}} \, dx \leq C M \|\wt N_{a,q}(v|\Phi)\|_q^q,
\end{equation}
with a constant $C>0$ that depends only on $n$, $a$.
\end{lem}

\bp The proof is almost straightforward. Let $b = a/100$. First, due to $(ii)$ in Definition \ref{defH2}, we have
$$\iint_\Omega |v|^q \Phi |\nabla \Psi|^2 \, \frac{dt}{|t|^{n-d-2}} \, dx \leq C \iint_\Omega |v|^q \Phi |\nabla \Psi| \, \frac{dt}{|t|^{n-d-1}} \, dx,$$
so we only need to bound the second term in the left-hand side of \eqref{H2la}.
By Fubini's lemma, for any non-negative function $f$,
\[ \begin{split}
	\iint_{(z,r)\in \R^d \times(0,\infty)} & \left(\frac{1}{|W_b(z,r)|} \iint_{W_b(z,r)} f dt\, dx \right) dz\, dr \approx \iint_{(x,t)\in\Omega} f \left( \iint_{(z,r)\in \R^d \times(0,\infty) \atop{\text{s.t. } (x,t)\in W_b(z,r)}}  dz\, \frac{dr}{r^n} \right) \, dt \, dx \\
	& \hspace{2cm} \geq \iint_{\Omega} f  \left( \int_{\frac{|t|}{2}}^{2|t|} \int_{z\in B_{b|t|/2}(x)} dz \, \frac{dr}{r^n} \right) \, dt\, dx \approx \iint_{\Omega} f \frac{dt}{|t|^{n-d-1}}\, dx,
\end{split}  \]
where the constants depend on $b$ and the dimension $d$.
Thus 
\[\begin{split}
\iint_\Omega |v|^q \Phi |\nabla \Psi| \, \frac{dt}{|t|^{n-d-1}} \, dx 
& \lesssim \iint_{(z,r) \in \R^d \times (0,\infty)} \left( \frac1{|W_b(z,r)|} \iint_{W_b(z,r)} |v|^q \Phi |\nabla \Psi| \right) dz \, dr \\
& \lesssim \iint_{(z,r) \in \R^d \times (0,\infty)} |(v|\Phi)_{W,b,q}(z,r)|^q \left(\sup_{W_b(z,r)} |\nabla \Psi| \right) dz \, dr \\
& \lesssim \int_{x\in \R^d} \iint_{(z,r) \in \Gamma_a(x)} |(v|\Phi)_{W,b,q}(z,r)|^q \left(\sup_{W_b(z,r)} |\nabla \Psi| \right) dz \, \frac{dr}{r^d} \, dx \\
& \lesssim \int_{x\in \R^d} |\wt N_{b,q}(v|\Phi)(x)|^q \iint_{(z,r) \in \Gamma_b(x)} \left(\sup_{W_b(z,r)} |\nabla \Psi|\right) dz \, \frac{dr}{r^d} \, dx \\
& \lesssim M \int_{\R^d} |\wt N_{a,q}(v|\Phi)|^q \, dx,
\end{split}\]
by $(iii)$ in Definition \ref{defH2}.
\ep

\section{The local estimate $S<N$}

\label{sS<N}

In this section, we estimate the $q$-adapted functional $S_{a,q}$ by the functional $\wt N_{a,q}$, by using methods similar to the ones in \cite{DPP} and \cite[Section 7]{DFMAinfty}. In fact, the general method was pioneered in \cite{KePiDrift}. To get boundary estimates, the derivative(s) in the transversal $t$-direction clearly plays an essential role. To estimate the $t$-derivative(s), it is convenient to tweak the elliptic matrix $\cA$ so that its lower left corner becomes zero:
\begin{equation} \label{defA'2}
\cA = \begin{pmatrix} \cA_{1} & \cA_2 \\ \cB_3 & bI \end{pmatrix} + \cC \, \longrightarrow \, \cA' : = \begin{pmatrix} \cA_{1} & \cA_2 + (\cB_3)^T \\ 0 & bI \end{pmatrix} + \cC.
\end{equation}
By simple algebra, such an $\cA'$ satisfies:
\begin{itemize}
	\item Let $X=(0,t)$ and $Y=(y,s)\in \R^n$, we have 
		\begin{equation} \label{eq:SNalg}
			\cA'Y\cdot X = b\, s\cdot t + \cC Y \cdot X.
		\end{equation}
		If $Y$ denotes the full derivative, then roughly speaking $\cA'$ highlights the $t$-derivative(s).
	\item The difference between the two elliptic operators $-\div |t|^{d+1-n} \cA \nabla$ and $-\div |t|^{d+1-n} \cA' \nabla$ is of first-order; to be more precise, 
		\begin{equation}\label{eq:SNAto}
			\cL= -\div |t|^{d+1-n} \cA \nabla= -\div |t|^{d+1-n} \cA' \nabla + |t|^{d+1-n} \mathcal D' \cdot \nabla, 
		\end{equation} 
	 where for any $1\leq j \leq d < i \leq n$,
		\begin{equation}\label{defB'2}
			(\cD')_j = - \sum_{\ell>d} |t|^{n-d-1} \dr_{t_\ell}[ (\cB_3)_{\ell j} |t|^{d+1-n}] \quad \text{ and } \quad (\cD')_i = \sum_{\ell\leq d} \dr_{x_\ell} (\cB_3)_{i\ell}.
		\end{equation}
	Observe that the condition ($\cH^1$) imposes in particular that $|t|\cD'$ satisfies the Carleson measure condition.
	\item Moreover, if we assume that the matrix $\cB_3$ is real-valued, then $\cA'\xi \cdot \overline \xi$ has the same real part as $\cA \xi \cdot \overline\xi$ [and we have $\cA'\xi \cdot \xi = \cA \xi \cdot \xi$ without even assuming that $\cB_3$ is real]. Hence, by Proposition \ref{defmuA}, $\cA'$ is $q$-elliptic if and only if $\cA$ is $q$-elliptic; more precisely $\lambda_q(\cA') = \lambda_q(\cA)$.
\end{itemize}

\begin{lem}[Local $S<\kappa N$] \label{lemS<Na}
Let $\cL = - \div |t|^{d+1-n} \cA \nabla$ be an elliptic operator that satisfies $(\cH^1_\kappa)$ for some constant $\kappa \geq 0$ (see Definition \ref{defH1}).
We define $\cA'$ as in \eqref{defA'2}.
Let $a>0$ and $q\in (q_0,q'_0)$ where $q_0$ is given by Proposition \ref{prop1.1}. Take $k >2$ and a  function $\Psi \in C^\infty_0(\Omega)$ which verifies $0\leq \Psi \leq 1$ and $\ds |\nabla \Psi| \leq \frac{100}{|t|}$ everywhere.
Then for any weak solution $u\in W^{1,2}_{loc}(\Omega)$ to $\cL u = 0$, we have
\begin{equation} \label{S<Na} \begin{split}
c \|S_{a,q}(u|\Psi^{k})\|_q^q
& \leq 
\frac12 \ \Re \iint_{\Omega} \cA' \nabla u \cdot \nabla [|u|^{q-2} \bar u] \frac{\Psi^{k}}{b}\, \frac{dt}{|t|^{n-d-2}} dx \\
& \leq  C \kappa \|\wt N_{a,q} (u|\Psi^{k-2})\|_{q}^q  - \,  \Re \iint_{\Omega} \cA'\nabla u \cdot \nabla_x [\Psi^{k}] \, \frac{|u|^{q-2} \bar u}{b} \frac{dt}{|t|^{n-d-2}}\, dx \\
& \qquad \qquad  + \, \frac1q \, \iint_{\Omega} \left(\frac{|u|^q}{|t|}\nabla |t| - \nabla_t [|u|^q] \right)\cdot \nabla_t [\Psi^{k}] \, \frac{dt}{|t|^{n-d-2}}\, dx,
\end{split} \end{equation}
where the constants $c,C$ depend on $a$, $q$, $n$, $k$, $\|\cA\|_\infty$, $\lambda_q(\cA)$ and $\|b\|_{\infty} + \|b^{-1}\|_\infty$. 
\end{lem}

\begin{rmk}
	\begin{enumerate}
		\item In \eqref{S<Na} the differential operators $\nabla_x$ and $\nabla_t$ denote an $n$-dimensional vector, where the missing derivatives are taken to be zero.
		\item  This lemma is a key step in proving the $S<N$ estimates. The corollaries of this lemma are Lemma \ref{lemS<Nb}, and with some additional a priori estimate, Lemma \ref{lemS<Ne}. In fact, the second and last terms on the right hand side can be roughly bounded by $\|\wt N_{a,q}(u|\Psi^{k-1})\|_q^q + \|S_{a,q}(u|\Psi^k)\|_q^{q/2} \|\wt N_{a,q}(u|\Psi^{k-2})\|_q^{q/2}$ (see Lemma \ref{lemS<Nb}); with an appropriate choice of cut-off function they can be bounded more precisely by $\|\Tr u\|_q^q$ (see Lemmas \ref{lem7.1'} and \ref{lem7.4}).
	\end{enumerate}
\end{rmk} 

\bp
First of all, Proposition \ref{lemMoserI} guarantees that $|u|^q$ and $|\nabla u|^2 |u|^{q-2}$ lie in $L^1_{loc}(\Omega)$. Therefore, since $\cA'$ and $b^{-1}$ are bounded and $|\nabla u|^2 |u|^{q-2} \simeq |\nabla u||\nabla [u^{q-2} \bar u]|$ by \eqref{nablau&v}, all the quantities invoked in \eqref{S<Na} are well defined and finite.

\ms

Let $T$ be the quantity
$$T:= \Re \iint_{\Omega} \cA' \nabla u \cdot \nabla [|u|^{q-2} \bar u] \frac{\Psi^{k}}{b}\, \frac{dt}{|t|^{n-d-2}} dx.$$
Since $\cA$ is $q$-elliptic and $\cB_3$ is real valued, we know $\cA'$ is also $q$-elliptic, and more precisely, that $\lambda_q(\cA') = \lambda_q(\cA)$. Together with Proposition \ref{proppellip} and the fact that $b$ is uniformly bounded from above, we obtain
$$\|S_{a,q}(u|\Psi^{k})\|_q^q \simeq \iint_{(x,t)\in \R^n}  |\nabla u|^2 |u|^{q-2} \Psi^{k}\, \frac{dt}{|t|^{n-d-2}} dx \lesssim T,$$
which is exactly the first estimate in \eqref{S<Na}. Here we use the assumption that $b$ is a real-valued scalar.

\ms

We claim that 
\begin{equation} \label{S<Nb} \begin{split}
T & \leq  C \kappa^\frac12 T^\frac12 \|\wt N_{q,a} (u|\Psi^{k-2})\|_{L^q(B)}^{q/2} - \,  \Re \iint_{\Omega} \cA'\nabla u \cdot \nabla_x [\Psi^{k}] \, \frac{|u|^{q-2} \bar u}{b} \frac{dt}{|t|^{n-d-2}}\, dx \\
& \hspace{4cm} + \, \frac1q \,  \iint_{\Omega} \left(\frac{|u|^q}{|t|}\nabla |t| - \nabla_t [|u|^q] \right)\cdot \nabla_t [\Psi^{k}] \, \frac{dt}{|t|^{n-d-2}}\, dx,
\end{split} \end{equation}
which easily implies the second inequality in \eqref{S<Na} since $T$ is finite. So let us prove the claim.
 
First, check that we can write
\[\begin{split}
T & = \Re \iint_{\Omega} \cA' \nabla u \cdot \nabla [|u|^{q-2} \bar u] \frac{|t| \Psi^{k}}{b}\, \frac{dt}{|t|^{n-d-1}} dx \\
& =\Re \iint_{\Omega} \cA'\nabla u \cdot \nabla\left[\frac{|t||u|^{q-2} \bar u \Psi^{k}}{b}\right] \frac{dt}{|t|^{n-d-1}}\, dx \\
& \quad - \,  \Re \iint_{\Omega} \cA'\nabla u \cdot \nabla [\Psi^{k}] \, \frac{|u|^{q-2} \bar u}{b} \frac{dt}{|t|^{n-d-2}}\, dx \\
& \quad + \Re \iint_{\Omega} \cA'\nabla u \cdot \nabla b \, \frac{|u|^{q-2} \bar u \Psi^{k}}{b^2} \frac{dt}{|t|^{n-d-2}}\, dx \\
& \quad - \Re \iint_{\Omega} \cA'\nabla u \cdot \nabla |t| \, \frac{|u|^{q-2} \bar u \Psi^{k}}{b} \frac{dt}{|t|^{n-d-1}}\, dx \\
& \quad := T_1 + T_2 + T_3 + T_4.
\end{split}\]
Let us start by $T_2$. The term $T_2$ is similar to the second term in the right-hand side of \eqref{S<Nb}, except that $\nabla [\Psi^{k}]$ is replaced by $\nabla_x [\Psi^{k}]$. We write
\[\begin{split}
T_2 & = -  \Re \iint_{\Omega} \cA'\nabla u \cdot \nabla_x [\Psi^{k}] \, \frac{|u|^{q-2} \bar u}{b} \frac{dt}{|t|^{n-d-2}}\, dx \\
& \qquad -  \Re \iint_{\Omega} \cA'\nabla u \cdot \nabla_t [\Psi^{k}] \, \frac{|u|^{q-2} \bar u}{b} \frac{dt}{|t|^{n-d-2}}\, dx \\
& \qquad := T_{21} + T_{22}
\end{split}\]
The term $T_{21}$ is the second term in the right-hand side of \eqref{S<Nb}. As for $T_{22}$, by \eqref{eq:SNalg} we know that
\[\begin{split}
T_{22} \leq &  - \Re \iint_{\Omega} \nabla_t u \cdot \nabla_t [\Psi^{k}] \, |u|^{q-2} \bar u \frac{dt}{|t|^{n-d-2}}\, dx \\
& + C \iint_{\Omega} |\cC| |\nabla u| |\nabla_t [\Psi^{k}]| \, |u|^{q-1}\frac{dt}{|t|^{n-d-2}}\, dx \\
& := T_{221} + T_{222},
\end{split}\]
where $T_{222}$ was obtained by using the fact that $b^{-1}$ is uniformly bounded. Notice that $ \nabla |u|^q = q \, \Re (\nabla u \, |u|^{q-2} \bar u)$, which gives that $T_{221}$ is part of the last term in the right-hand side of \eqref{S<Nb}. 
Observe now that by assumption on $\Psi$, we have 
$$|\nabla_t [\Psi^{k}]| \leq \frac{C}{|t|} \Psi^{k-1},$$
hence
\[\begin{split}
T_{222} & \leq C \int_{\Omega} |\cC| |\nabla u| \Psi^{k-1} |u|^{q-1} \frac{dt}{|t|^{n-d-1}}\, dx \\
& \leq C \left( \int_{\Omega} |\nabla u|^2 |u|^{q-2} \Psi^{k} \frac{dt}{|t|^{n-d-2}}\, dx \right)^\frac12 \left( \int_\Omega |\cC|^2 |u|^q \Psi^{k-2} \frac{dt}{|t|^{n-d}}\, dx\right)^{\frac12}
\end{split}\]
Since the matrix $\cA'$ is $q$-elliptic and $b^{-1} \leq C$, the first integral in the right-hand side above is thus bounded by 
$$\Re \int_{\Omega} \cA' \nabla u \cdot \nabla[ |u|^{q-2} \bar u] \frac{\Psi^{k}}{b} \frac{dt}{|t|^{n-d-2}}\, dx = T,$$
while we bound the second integral by using Proposition \ref{propCarl} (the Carleson inequality) and the fact that $\cC$ satisfies the Carleson measure condition. We obtain
\[\begin{split}
T_{222} \lesssim T^{1/2} \kappa^\frac12 \|\wt N_{a,q}(u|\Psi^{k-2})\|_{q}^{q/2}.
\end{split}\]
The bound of $T_2$ follows.

We turn to the estimate of $T_3$. We use the boundedness of $\cA'$ and $b^{-1}$, and then Cauchy-Schwarz's inequality to obtain
\[\begin{split}
|T_3| & \lesssim  \iint_\Omega |\nabla u| |\nabla b| \Psi^{k} |u|^{q-1} \frac{dt}{|t|^{n-d-2}}\, dx \\
&  \lesssim  \left( \int_{\Omega} |\nabla u|^2 |u|^{q-2} \Psi^{k} \frac{dt}{|t|^{n-d-2}}\, dx \right)^\frac12 \left( \iint_\Omega |\nabla b|^2 \Psi^{k} |u|^q \frac{dt}{|t|^{n-d-2}}\, dx \right)^\frac12.
\end{split}\]
As we did for $T_{222}$, the $q$-ellipticity of $\cA'$ and the fact that $|t||\nabla b|$ satisfies the Carleson measure condition give that
\[\begin{split}
|T_{3}| \lesssim T^{1/2} \kappa^\frac12 \|\wt N_{a,q}(u|\Psi^{k})\|_{q}^{q/2} \leq T^{1/2} \kappa^\frac12 \|\wt N_{a,q}(u|\Psi^{k-2})\|_{q}^{q/2},
\end{split}\]
where the last inequality stands because $\Psi \leq 1$.

The term $T_{4}$ is a bit more complicated. By \eqref{defA'2}, the bottom left corner of $\cA'$ is a zero matrix. Thus, similarly to \eqref{eq:SNalg}, we have
\[\begin{split}
T_4 
& = - \Re \iint_{\Omega} b \nabla_t  u \cdot \nabla_t |t| \, \frac{|u|^{q-2} \bar u \Psi^{k}}{b} \frac{dt}{|t|^{n-d-1}}\, dx - \Re \iint_{\Omega} \cC \nabla  u \cdot \nabla_t |t| \, \frac{|u|^{q-2} \bar u \Psi^{k}}{b} \frac{dt}{|t|^{n-d-1}}\, dx \\
 & = - \Re \iint_{\Omega} \cC \nabla  u \cdot \nabla_t |t| \, \frac{|u|^{q-2} \bar u \Psi^{k}}{b} \frac{dt}{|t|^{n-d-1}}\, dx - \Re \iint_{\Omega}  \nabla_t  u \cdot \nabla_t |t| \, |u|^{q-2} \bar u \Psi^{k} \frac{dt}{|t|^{n-d-1}}\, dx\\
& := T_{41} + T_{42}.
\end{split}\]
The term $T_{41}$ can be dealt as $T_{222}$ or $T_3$, by using the Cauchy-Schwarz's inequality and then Carleson inequality (Proposition \ref{propCarl}), and we have
$$|T_{41}| \lesssim T^{1/2} \kappa^\frac12 \|\wt N_{a,q}(u|\Psi^{k})\|_{q}^{q/2} \leq T^{1/2} \kappa^\frac12 \|\wt N_{a,q}(u|\Psi^{k-2})\|_{q}^{q/2}.$$
As for $T_{42}$, recall that $\nabla_t |u|^q = q |u|^{q-2} \Re( \bar u \nabla _t u)$,
hence
\[\begin{split}
T_{42} & = - \frac1q \iint_{\Omega}  \nabla_t [ |u|^q \Psi^{k}] \cdot \nabla_t |t| \frac{dt}{|t|^{n-d-1}}\, dx\\
& \quad + \frac1q \iint_{\Omega}  \nabla_t [\Psi^{k}] \cdot \nabla_t |t| \, |u|^{q} \frac{dt}{|t|^{n-d-1}}\, dx \\
& := T_{421} + T_{422}.
\end{split}\]
The first term $T_{421}$ is 0. Indeed, for any function $v$,  $\nabla_t v \cdot \nabla_t |t|$ equals $\dr_r v$, the derivative in the radial  direction. We switch then to polar coordinate and $T_{421}$ is the integral of a derivative. 
The last term $T_{422}$ is part of the last term in the right-hand side of \eqref{S<Nb}.

\ms

It remains to treat $T_1$. 
As stated in \eqref{eq:SNAto}, we know that if $u$ is a weak solution of $\cL u = 0$, then $u$ is also a weak solution to $\cL' u  = 0$, where $\cL' = -\div |t|^{d+1-n} \cA' \nabla + |t|^{d+1-n} \mathcal D' \cdot \nabla$, the matrix $\cA'$ is defined in \eqref{defA'2}, and $\cD'$ as in \eqref{defB'2}.
Note that $|t|\cD'$ satisfies the Carleson measure condition.

Now, set $v:=\frac{|t|}{b} |u|^{q-2}\bar u \Psi^{k}$. We want to find the value of 
\begin{equation}\label{L'u=0}
\iint_{(x,t)\in \R^n} \cA' \nabla u \cdot \nabla v \, \frac{dt}{|t|^{n-d-1}} dx + \iint_{(x,t)\in \R^n} \cD'\cdot \nabla u \, v \frac{dt}{|t|^{n-d-1}} dx,
\end{equation}
which is formally 0, because \eqref{L'u=0} is only $\cL' u$ tested against a test function $v$. The problem is that the test function $v$ is not smooth, and not even necessary in $W^{1,2}_{loc}(\Omega)$, so the fact that \eqref{L'u=0} equals 0 needs to be proven. 

The function $\Psi$ is compactly supported in $\Omega$, hence is $v$; besides, since $\Psi,b^{-1},\nabla \Psi,\nabla b$ are all in $L^\infty_{loc}(\Omega)$, we have $|v| \lesssim |u|^{q-1}$ and $|\nabla v| \lesssim |\nabla u| |u|^{q-2} + |u|^{q-1}$ uniformly on the support of $\Psi$. Due to Proposition \ref{lemMoserI}, the function $\cA' \nabla u \cdot \nabla v+ \cD'\cdot \nabla u \, v$ is now in $L^1(\Omega, \frac{dt}{|t|^{n-d-1}} dx)$ and we will see that an argument similar to the Lebesgue domination convergence proves that \eqref{L'u=0} is 0. We separate two cases: if $q\geq 2$, then we set $u_N = \min\{N,|u|\}$ and
$$v_N :=\frac{|t|}{b} u_N^{q-2} \bar u \Psi^{k}.$$
The function $v_N$ is in $W^{1,2}_{loc}(\Omega)$, because $|t|\Psi^{2q}$ is Lipschitz, $u_N^{q-2}\bar u \in W^{1,2}_{loc}(\Omega)$ since that  $|\nabla [u_N^{q-2} \bar u]| \leq (q-1) u_N^{q-2} |\nabla \bar u| \in L^2_{loc}(\Omega)$, and $|b|\geq C^{-1}$, $b\in W^{1,\infty}_{loc}(\Omega)$ by assumption. Moreover, $v_N$ is compactly supported in $\Omega$, so by \cite[Lemma 8.16]{DFMprelim}, for any $N\in \N$,
$$\iint_{(x,t)\in \R^n} \cA' \nabla u \cdot \nabla v_N \, \frac{dt}{|t|^{n-d-1}} dx + \iint_{(x,t)\in \R^n} \cD'\cdot \nabla u \, v_N \frac{dt}{|t|^{n-d-1}} dx =0.$$
In addition, we have
$$\cA' \nabla u \cdot \nabla v+ \cD'\cdot \nabla u \, v = \cA' \nabla u \cdot \nabla v_N+ \cD'\cdot \nabla u \, v_N \quad \text{ in } \{|u|  < N\}.$$
It follows that 
\[\begin{split}
& \left| \iint_{(x,t)\in \R^n} [\cA' \nabla u \cdot \nabla v + \cD'\cdot \nabla u \, v] \, \frac{dt}{|t|^{n-d-1}} dx \right| \\
& \hspace{3cm} \leq \int_{\{|u|\geq N\}} (|\cA' \nabla u \cdot \nabla v + \cD'\cdot \nabla u \, v| + |\cA' \nabla u \cdot \nabla v_N + \cD'\cdot \nabla u \, v_N|) \frac{dt}{|t|^{n-d-1}} dx \\
& \hspace{3cm} \longrightarrow 0 \quad \text{ as } N \to +\infty,
\end{split}\]
because $\cA' \nabla u \cdot \nabla v + \cD'\cdot \nabla u \, v$ and $\cA' \nabla u \cdot \nabla v_N + \cD'\cdot \nabla u \, v_N$ are integrable.

We turn now to the case where $q<2$. We define $u_\epsilon := \max\{|u|,\epsilon\}$, $v_\epsilon := \frac{|t|}{b} u_\epsilon^{q-2} \bar u \Psi^{k}.$ Check similarly as before that $v_\epsilon \in W^{1,2}_{loc}(\Omega)$, which gives 
$$\iint_{(x,t)\in \R^n} \cA' \nabla u \cdot \nabla v_\epsilon \, \frac{dt}{|t|^{n-d-1}} dx + \iint_{(x,t)\in \R^n} \cD'\cdot \nabla u \, v_\epsilon \frac{dt}{|t|^{n-d-1}} dx =0$$
and, by the Lebesgue domination theorem, since $\cA' \nabla u \cdot \nabla v + \cD'\cdot \nabla u \, v$ and $\cA' \nabla u \cdot \nabla v_\epsilon + \cD'\cdot \nabla u \, v_\epsilon$ are integrable,
\[\begin{split}
& \left| \iint_{(x,t)\in \R^n} [\cA' \nabla u \cdot \nabla v + \cD'\cdot \nabla u \, v] \, \frac{dt}{|t|^{n-d-1}} dx \right| \\
& \hspace{2cm} \leq \int_{\{|u| \leq \epsilon\}} (|\cA' \nabla u \cdot \nabla v + \cD'\cdot \nabla u \, v| + |\cA' \nabla u \cdot \nabla v_\epsilon + \cD'\cdot \nabla u \, v_\epsilon|) \frac{dt}{|t|^{n-d-1}} dx \\
& \hspace{2cm} \longrightarrow \int_{\{|u|=0\}} (|\cA' \nabla u \cdot \nabla v + \cD'\cdot \nabla u \, v| + |\cA' \nabla u \cdot \nabla v_\epsilon + \cD'\cdot \nabla u \, v_\epsilon|) \frac{dt}{|t|^{n-d-1}} dx
\end{split}\]
as $\epsilon \to 0$. 
However, $\nabla u = 0$ almost everywhere on $\{|u|=0\}$. We deduce
$$ \iint_{(x,t)\in \R^n} [\cA' \nabla u \cdot \nabla v + \cD'\cdot \nabla u \, v] \, \frac{dt}{|t|^{n-d-1}} dx = 0.$$

We just have shown that the term in \eqref{L'u=0} is zero, thus our term $T_{1}$ becomes
$$|T_{1}| \lesssim \left| \iint_{(x,t)\in \R^n} \cD'\cdot \nabla u \, v \frac{dt}{|t|^{n-d-1}} dx \right|\lesssim \iint_{(x,t)\in \R^n} |t| |\cD'| |\nabla u| |u|^{q-1}\Psi^{k}\frac{dt}{|t|^{n-d-1}} dx . $$
Recall that $|t|\cD'$ satisfies the Carleson measure condition. So using Cauchy-Schwarz's inequality and then Carleson's inequality (Proposition \ref{propCarl}), similarly to $T_{222}$, $T_3$ or $T_{41}$,
$$|T_{1}| 
\lesssim T^{1/2} \kappa^{\frac12} \|\wt N_{a,q}(u|\Psi^{k-2})\|_{L^q(B)}^{q/2}.$$
The estimate \eqref{S<Nb} and then the lemma follows.
\ep

From Lemma \ref{lemS<Na}, we can easily deduce the local bounds given in the following lemma.

\begin{lem} \label{lemS<Nb}
Let $\cL = - \div |t|^{d+1-n} \cA \nabla$ be an elliptic operator that satisfies $(\cH^1_\kappa)$ for some constant $\kappa \geq 0$. Let $a>0$ and $q\in (q_0,q'_0)$ where $q_0$ is given by Proposition \ref{prop1.1}.
Choose a  function $\Psi \in C^\infty_0(\Omega)$ that satisfies $(\cH^2_{a,M})$ and $k>2$.
Then, for any weak solution $u\in W^{1,2}_{loc}(\Omega)$, we have
\begin{equation} \label{S<Ncc}
\|S_{a,q}(u|\Psi^k)\|_{q}^q \leq C( \kappa + M) \|\wt N_{a,q}(u|\Psi^{k-2})\|_{q}^q
\end{equation}
where the constant $C$ depends on $a$, $q$, $n$, $\|\cA\|_\infty$, $\lambda_q(\cA)$, $\|b\|_{\infty} + \|b^{-1}\|_\infty$, and $k$.
\end{lem}

\bp
First, recall that
\begin{equation}\label{claim6}
	\|S_{a,q,e}^l(u|\Psi^k) \|_{q}^q \simeq T:= \iint_{\R^n} |\nabla u|^2 |u|^{q-2} \Psi^{k} \frac{ds}{|s|^{n-d-2}} dy.
\end{equation}
By Lemma \ref{lemS<Na}, we have 
\[\begin{split}
\|S_{a,q}(u|\Psi^k) \|_{q}^q & \leq C \kappa \|\wt N_{q,a} (u|\Psi^{k-2})\|_{L^q(4B)}^q  - \,  \Re \iint_{\Omega} \cA'\nabla u \cdot \nabla_x [\Psi^{k}] \, \frac{|u|^{q-2} \bar u}{b} \frac{dt}{|t|^{n-d-2}}\, dx \\
& \qquad \qquad  + \, \frac1q \,  \Re \iint_{\Omega} \left(\frac{|u|^q}{|t|}\nabla |t| - \nabla_t [|u|^q] \right)\cdot \nabla_t [\Psi^{k}] \, \frac{dt}{|t|^{n-d-2}}\, dx \\
& \qquad := T_1 + T_2 + T_3.
\end{split}\]
We claim
\begin{equation} \label{claima}
|T_2 + T_3| \lesssim M \|\wt N_{q,a} (u|\Psi^{k-1})\|_{q}^q + T^{1/2}M^{1/2} \|\wt N_{q,a} (u|\Psi^{k-2})\|_{q}^{q/2}.
\end{equation}
From the claim we deduce that $T \lesssim (\kappa + M) \|\wt N_{q,a} (u|\Psi^{k-2})\|_{L^q(4B)}^q$. Together with \eqref{claim6}, the lemma follows.

So it remains to prove \eqref{claima}. By using the fact that $\nabla_t |u|^q = q|u|^{q-2} \Re (\bar u\nabla_t  u)$ and that $\cA'/b$ is bounded, we get
\[\begin{split}
|T_2+T_3| & \lesssim \iint_\Omega |\nabla u| |\nabla [\Psi^{k}]| |u|^{q-1} \frac{dt}{|t|^{n-d-2}}\, dx + \iint_\Omega |u|^q |\nabla [\Psi^{k}]| \frac{dt}{|t|^{n-d-1}}\, dx \\
& \lesssim \iint_\Omega |\nabla u| \Psi^{k-1} |\nabla \Psi| |u|^{q-1} \frac{dt}{|t|^{n-d-2}}\, dx + \iint_\Omega |u|^q \Psi^{k-1} |\nabla \Psi| \frac{dt}{|t|^{n-d-1}}\, dx \\
& := T_4 + T_5.
\end{split}\]
Applying Lemma \ref{lemH2} to $\Phi = \Psi^{k-1}$, we get 
$$T_5 \lesssim M \|\wt N_{q,a} (u|\Psi^{k-1})\|_{q}^q.$$
As for $T_4$, we use Cauchy-Schwarz's inequality to get
\[\begin{split}
\iint_\Omega |\nabla u| \Psi^{k-1} |\nabla \Psi| |u|^{q-1} \frac{dt}{|t|^{n-d-2}}\, dx & \leq T^{1/2} \left( \iint_\Omega |u|^{q}  \Psi^{k-2} |\nabla \Psi|^2 \frac{dt}{|t|^{n-d-2}}\, dx \right) \\
& \lesssim T^{1/2} M^{1/2} \|\wt N_{q,a} (u|\Psi^{k-2})\|_{q}^q,
\end{split}\]
by applying Lemma \ref{lemH2} again to $\Phi = \Psi^{k-2}$.
The claim \eqref{claima} and then the lemma follows.
\ep

We denote by $\mathcal M$ the Hardy-Littlewood maximal function on $\R^d$, that is, if $f$ is a locally integrable function on $\R^d$,
$$\mathcal M f(x) := \sup_{\text{balls } B \ni x} \fint_{B} |f|.$$
As it will be useful later on, we also introduce here the maximal operator $\mathcal M_q$ defined for any $q\in (1,+\infty)$ on $L^q_{loc}(\R^d)$ as
$$\mathcal M_q[f](x) = [\mathcal M [f^q](x)]^{\frac1q}.$$
It is well known that the operator $\mathcal M_q$ is weak type $(q,q)$ and strong type $(p,p)$ for $p>q$.

\ms

A key result of our paper is Lemma \ref{lemS<Nd}, which states that the $L^q$ bounds for the $q$-adapted square function given in Lemma \ref{lemS<Nb} self-improve to $L^p$ bounds for all $p>0$. We start by proving a preliminary good-$\lambda$ inequality.

\begin{lem} \label{lemS<Nc}
Let $\cL = - \div |t|^{d+1-n} \cA \nabla$ be an elliptic operator that satisfies $(\cH^1_\kappa)$ for some constant $\kappa \geq 0$. Let $a>0$ and $q\in (q_0,q'_0)$ where $q_0$ is given by Proposition \ref{prop1.1}.
Choose a  function $\Psi \in C^\infty_0(\Omega)$ that satisfies $(\cH^2_{a,M})$ and $k>2$.
Then, there exists $\eta \in (0,1)$ that depends only on $d$ and $q$ such that for any weak solution $u\in W^{1,2}_{loc}(\Omega)$, any $\nu >0$ and any $\gamma \in (0,1)$,
\begin{align} \label{S<Nc} 
|\{ x\in \R^d, \ S_{a,q}(u|\Psi^k)(x)>\nu,\, \mathcal M[\wt N_{a,q}(u|\Psi^{k-2})](x) & \leq \gamma \nu \}| \nn\\
& \leq C \gamma^q |\{x\in \R^d, \ \mathcal M[S_{a,q}(u|\Psi^k)](x)>\eta \nu\}|,
\end{align}
where the constant $C$ depends on $a$, $q$, $n$, $\|\cA\|_\infty$, $\lambda_q(\cA)$, $\|b\|_{\infty} + \|b^{-1}\|_\infty$, $\kappa$, $k$ and $M$.
\end{lem}

\bp
Let $\eta$ to be chosen later on. Take $\nu>0$. We define the set
$$\mathcal S := \{x\in \R^d, \ \mathcal M[S_{a,q}(u|\Psi^k)](x)>\eta \nu\},$$
which is open [$\Psi$ is compactly supported in $\Omega$, which makes $S_{a,q}(u|\Psi^k)$ continuous] and bounded. Indeed, $\mathcal S$ is bounded because $\Psi^k$ is compactly supported so $S_{a,q}(u|\Psi^k) \equiv 0$ outside of a big ball.  

We construct a Whitney decomposition of $\mathcal S$ in the following manner. For any $x\in \mathcal S$, we set $B_x$ as the ball of center $x$ and radius $\frac1{10}\dist(x,\mathcal S^c)$. The balls $(B_x)_{x\in \mathcal S}$ clearly cover $\mathcal S$. Moreover, the radii of the balls are uniformly bounded, since $\mathcal S$ is bounded. 
So by Vitali's lemma, there exists a non-overlapping sub-collection of  balls $(B_{x_i})_{i\in I}$ such that $\bigcup_{i\in I} 5B_{x_i} \supset \mathcal S$. We write $B_i = 10B_{x_i}$ and $l_i$ for its radius. By construction, for every $i\in I$,
\begin{equation} \label{defyi}
\text{there exists $y_i \in \R^d$ such that $|x_i - y_i| = l_i$ and $\ds \mathcal M[S_{a,q}(u|\Psi^k)](y_i) \leq \eta\nu$}.
\end{equation}
The balls $B_{x_i} = \frac1{10} B_i$ are mutually disjoint sets contained in $\mathcal S$, so we deduce
\begin{equation} \label{claimb}
\sum_{i\in I} |B_i| = 10^d \sum_{i\in I} |B_{x_i}| \lesssim 10^d |\mathcal S|.
\end{equation}

\ms

Thanks to \eqref{claimb}, the estimate \eqref{S<Nc} will be obtained if we can prove that
\begin{equation} \label{claimc}
F^i_\gamma := \{ x\in B_i, \ S_{a,q}(u|\Psi^k)(x)>\nu,\, \mathcal M[\wt N_{a,q}(u|\Psi^{k-2})](x) \leq \gamma \nu \} \leq C \gamma^q |B_i|, 
\end{equation}
where $C$ is independent of $\gamma \in (0,1)$. If $F^i_\gamma = \emptyset$, there is nothing to prove, so we can assume $F^i_\gamma \neq \emptyset$. 
Set the function $\Phi_i$ as
$$\Phi_i(x,t) = \phi\left( \frac{a|t|}{l_i}\right) \phi\left( \frac{|x-x_i|}{2l_i}\right),$$
where $\phi \in C^\infty_0(\R)$ is such that $0\leq \phi \leq 1$, $\phi \equiv 1$ on $[-1,1]$, $\phi \equiv 0$ outside $[-2,2]$, and $|\phi'| \leq 2$. We claim that we can find $\eta$ small such that 
\begin{equation} \label{claimd}
S_{a,q}(u|\Psi^k \Phi_i^k)(z) \geq\frac\nu2 \qquad \text{ for } z\in F^i_\gamma.
\end{equation}
Indeed, for any $z\in \R^d$, one has by definition
\begin{equation} \label{eq:SNtmp1}
	\begin{split}
S_{a,q}(u|\Psi^k [1-\Phi_i^k])(z) & = \left( \iint_{(y,s) \in \wh \Gamma_a(z)} |\nabla u|^2 |u|^{q-2} \Psi^k [1-\Phi_i^k] \frac{ds}{|s|^{n-2}} \, dy \right)^\frac1q.
\end{split}
\end{equation} 
By definition, $1-\Phi_i^k \equiv 0$ if $|s| \leq l_i/a$, thus in the above integral \eqref{eq:SNtmp1} we only need to consider the integration region $\wh \Gamma_a(z)$ above the level of $l_i/a$. 
Note that we can find $N$ points $(z_j)_{j\leq N} \in B_{l_i}(z)$ with $N$ depending only on $d$ (for example by taking each $z_j$ to be at least $l_i/3a$-distance away from all the other $z_j$'s), such that 
\[ \wh \Gamma_a(z) \cap \{|s|> l_i/a  \} \subset \bigcup_{i=1}^N \wh \Gamma_a(z_j,l_i/2a). \]
Here we define 
\begin{equation} \label{defwtGa}
\wh \Gamma_a(x,\rho) := \left\{(y,s) \in \Omega,\,  |y-x| < a\left( |s|- \rho \right) \right\},
\end{equation}
 which are cones raised to the level $\rho$. Hence
\[ \begin{split}
S_{a,q}(u|\Psi^k [1-\Phi_i^k])(z) & \leq \sum_{j=1}^N \left(  \iint_{(y,s) \in \wh \Gamma_a(z_j,l_i/2a)} |\nabla u|^2 |u|^{q-2} \Psi^k [1-\Phi_i^k] \frac{ds}{|s|^{n-2}} \, dy \right)^\frac1q,
\end{split}\]
By simple geometry we observe that
$$ \wh\Gamma_a(z_j,l_i/2a) \subset \wh \Gamma_a(z') \qquad \text{ for any } z'\in B_{l_i/2}(z_j);$$
if the point $z$ belongs to $F_\gamma^i \subset B_i = B(x_i, l_i)$ and $|x_i - y_i| = l_i$ (see \eqref{defyi}), we have $B_{l_i/2}(z_j) \subset B_{4l_i}(y_i)$. We conclude that
\[\begin{split}
S_{a,q}(u|\Psi^k [1-\Phi_i^k])(z) & \leq \sum_{j=1}^N  \fint_{z'\in B_{l_i/2}(z_j)} \left( \iint_{(y,s) \in \wh \Gamma_a(z')} |\nabla u|^2 |u|^{q-2} \Psi^k [1-\Phi_i^k] \frac{ds}{|s|^{n-2}} \, dy\, \right)^\frac1q dz' \\
& \lesssim \sum_{j=1}^N  \fint_{z'\in B_{4l_i}(y_i)} |S_{a,q}(u|\Psi^k)(z')| dz' \\
& \leq C_d \mathcal M[S_{a,q}(u|\Psi^k)](y_i) \leq C_d \eta \nu
\end{split}\]
where we use \eqref{defyi} in the last inequality. We choose $\eta := \eta(d)$ such that $(C_d \eta)^q \leq (1-2^{-q})$. Consequently, $|S_{a,q}(u|\Psi^k [1-\Phi_i^k])(z)|^q \leq \nu(1-2^{-q})$ and then
$$|S_{a,q}(u|\Psi^k \Phi_i^k)(z)|^q = |S_{a,q}(u|\Psi^k)(z)|^q - |S_{a,q}(u|\Psi^k [1-\Phi_i^k])(z)|^q \geq \frac{\nu^q}{2^q},$$
which is the desired claim \eqref{claimd}. 

We define now the cut-off function $\chi_i$ as
$$\chi_i (x,t) = \phi\left( \frac{\dist(x,F^i_\gamma)}{a|t|} \right).$$
It is easy to check that for any $z\in F^i_\gamma$, the cone $\wh \Gamma_a(z) \subset \{ \dist(x,F^i_\gamma) \leq a|t| \}$, and thus $ \chi_i \equiv 1$ on $\wh \Gamma_a(z)$. In fact $\chi_i$ describes a smooth version of the classic sawtooth domain on top of $F_\gamma^i$. Therefore 
$$S_{a,q}(u|\Psi^k \Phi_i^k \chi_i^k)(z) = S_{a,q}(u|\Psi^k \Phi_i^k) \geq \frac{\nu}{2} \qquad \text{ for } z\in F^i_\gamma,$$
and
$$|F^i_\gamma| \leq \frac{2^q}{\nu^q} \int_{\R^d} |S_{a,q}(u|\Psi^k \Phi_i^k \chi_i^k)(z)|^q \, dz.$$
Since $\Psi$ satisfies $(\cH^2_{a,M})$, by Lemma \ref{exgoodcutoff} $\Psi \Phi_i \chi_i$ satisfies $(\cH^2_{M+M'})$ with $M'$ that depends only on $a$ and $n$. Lemma \ref{lemS<Nb} entails that
\begin{equation} \label{Figamma}
|F^i_\gamma| \lesssim \nu^{-q} \int_{\R^d} |\wt N_{a,q}(u|\Psi^{k-2} \Phi_i^{k-2} \chi_i^{k-2})(z)|^q \, dz.
\end{equation}
The function $\Phi_i$ is supported in $\{(x,t) \in \Omega, \, x\in 4B_i, \, |t| \leq 2l_i/a\}$. It forces $\wt N_{a,q}(u|\Psi^{k-2} \Phi_i^{k-2} \chi_i^{k-2})$ to be supported in - say - $20B_i$. So \eqref{Figamma} becomes
\begin{equation} \label{Figamma2}
|F^i_\gamma| \lesssim \nu^{-q} \int_{20B_i} |\wt N_{a,q}(u|\Psi^{k-2}\chi_i^{k-2})(x)|^q \, dx.
\end{equation}
For $(z,r) \in \R^d \times (0,\infty)$, we want to compute $(u|\Psi^{k-2} \chi_i^{k-2})_{\W,a,q}(z,r)$. On one hand, for any $z'' \in B_{ar}(z)$, we have $(z,r) \in \Gamma_a(z'')$ and thus
\[
(u|\Psi^{k-2} \chi_i^{k-2})_{\W,a,q}(z,r) \leq \wt N_{a,q}(u|\Psi^{k-2}\chi_i^{k-2})(z'') \leq \wt N_{a,q}(u|\Psi^{k-2})(z'') .
\]
Hence
\begin{equation}\label{eq:SNtmp2}
	(u|\Psi^{k-2} \chi_i^{k-2})_{\W,a,q}(z,r) \leq \fint_{B_{ar}(z)} \wt N_{a,q}(u|\Psi^{k-2})(z'') dz''.
\end{equation}
On the other hand, if $(u|\Psi^{k-2} \chi_i^{k-2})_{\W,a,q}(z,r)$ is not 0, it means that $W_a (z,r) \cap \supp \, \chi_i \neq \emptyset$, and so that there exists $z'\in F^i_\gamma$ such that $|z-z'| \leq 6ar$. 
Combining with \eqref{eq:SNtmp2} we get
\begin{equation}\label{eq:SNtmp3}
	(u|\Psi^{k-2} \chi_i^{k-2})_{\W,a,q}(z,r) \lesssim \fint_{B_{7ar}(z')} \wt N_{a,q}(u|\Psi^{k-2})(z'') dz'' \leq  \mathcal M [\wt N_{a,q}(u|\Psi^{k-2})](z') \leq \gamma \nu.
\end{equation} 
The last inequality is by the definition of $F_\gamma^i$.
Since the estimate \eqref{eq:SNtmp3} holds for all $(z,r)\in \R^d \times(0,\infty)$, we deduce that $\wt N_{a,q}(u|\Psi^{k-2}\chi_i^{k-2})(x) \lesssim \gamma \nu$ for all $x\in \R^d$. Therefore one can rewrite \eqref{Figamma2} as
\begin{equation} \label{Figamma3}
|F^i_\gamma| \lesssim \nu^{-q} |20B_i| (\gamma\nu)^q \lesssim \gamma^q |B_i|.
\end{equation}
The claim \eqref{claimc} and then the lemma follows.
 \ep

 \begin{lem} \label{lemS<Nd}
Let $\cL = - \div |t|^{d+1-n} \cA \nabla$ be an elliptic operator that satisfies $(\cH^1_\kappa)$ for some constant $\kappa \geq 0$. Let $a>0$, $1<p<\infty$, and $q\in (q_0,q'_0)$ where $q_0$ is given by Proposition \ref{prop1.1}.
Choose a  function $\Psi \in C^\infty_0(\Omega)$ that satisfies $(\cH^2_{a,M})$ and $k>2$.
Then, for any weak solution $u\in W^{1,2}_{loc}(\Omega)$, 
\begin{align} \label{S<Nz} 
\|S_{a,q}(u|\Psi^k)\|_p \leq C \|\wt N_{a,q}(u|\Psi^{k-2})\|_p
\end{align}
where the constant $C$ depends on $a$, $p$, $q$, $n$, $\|\cA\|_\infty$, $\lambda_q(\cA)$, $\|b\|_{\infty} + \|b^{-1}\|_\infty$, $\kappa$, $k$ and $M$.
\end{lem}

\begin{rmk}
The result is proven for $p>1$, but can be easily extended to $p>0$ if, in Lemma \eqref{S<Nc}, we replace the Hardy-Littlewood maximal function $\mathcal M$ by $\mathcal M_r$ with $0<r<1$.
\end{rmk}

\bp We have
\[\begin{split}
\int_{\R^d} |S_{a,q}(u|\Psi^k)|^p \, dx & = c_p \int_0^\infty \nu^{p-1} |\{x\in \R^d, \ S_{a,q}(u|\Psi^k)(x) > \nu\}| d\nu \\
& \leq C \int_0^\infty \nu^{p-1}  |\{x\in \R^d, \ S_{a,q}(u|\Psi^k)(x) > \nu, \, \mathcal M[\wt N_{a,q}(u|\Psi^{k-2})](x) \leq \gamma \nu \}|  d\nu \\
& \qquad + C \int_0^\infty \nu^{p-1} |\{x\in \R^d, \ \mathcal M[\wt N_{a,q}(u|\Psi^{k-2})](x) > \gamma \nu \}| d\nu
\end{split}\]
Lemma \ref{lemS<Nc} implies the existence of $\eta := \eta(d)$ such that
\[\begin{split}
\int_{\R^d} |S_{a,q}(u|\Psi^k)|^p \, dx
& \leq C \gamma^q \int_0^\infty \nu^{p-1}  |\{x\in \R^d, \ \mathcal M[S_{a,q}(u|\Psi^k)](x) > \eta\nu \}|  d\nu \\
& \qquad \qquad + C \gamma^{1-p} \int_0^\infty (\gamma\nu)^{p-1} |\{x\in \R^d, \ \mathcal M[\wt N_{a,q}(u|\Psi^{k-2})](x) > \gamma \nu \}| d\nu \\
& \leq C \gamma^q \int_{\R^d} |\mathcal M[S_{a,q}(u|\Psi^k)]|^p \, dx + C \gamma^{1-p} \int_{\R^d}  |\mathcal M[\wt N_{a,q}(u|\Psi^{k-2})]|^p\, dx \\
& \leq C \gamma^q \int_{\R^d} |S_{a,q}(u|\Psi^k)|^p \, dx + C \gamma^{1-p} \int_{\R^d}  |\wt N_{a,q}(u|\Psi^{k-2})|^p\, dx,
\end{split}\]
where the last inequality is due to the Hardy-Littlewood maximal inequality. We choose $\gamma$ such that $C\gamma^q = \frac12$. On the other hand, since $\Psi$ is compactly supported in $\Omega$, we know $\|S_{a,q}(u|\Psi^k)\|_p < +\infty$. Therefore
$$\int_{\R^d} |S_{a,q}(u|\Psi^k)|^p \, dx \leq C \int_{\R^d}  |\wt N_{a,q}(u|\Psi^{k-2})|^p\, dx.$$
The lemma follows.
\ep

\section{The local estimate $N < S$}

\label{N<S}

First, for $\rho>0$, we introduce similarly to \eqref{defwtGa} the new set defined as
$$\Gamma_a(x,\rho) := \left\{(z,r) \in \R^{d+1}_+,\,  |z-x| < a(r-\rho) \right\},$$
which is a $(d+1)$-dimensional cone raised to the level $\rho$.
Let $a>0$. For $\nu >0$ and $v$, a continuous and compactly supported function, we define the map $h_{\nu,a}(v):\,\R^d \to \R$ as
$$h_{\nu}(x) := h_{\nu,a}(v)(x) = \inf \left\{ r>0: \, \sup_{Y\in \Gamma_a(x,r)} v(Y) < \nu \right\}.$$

Since $v$ is compactly supported, it is pretty clear that $\left\{ r>0: \, \sup_{Y\in \Gamma_a(x,r)} v(Y) < \nu \right\}$ is non-empty and then that $h_{\nu}(x)$ is well defined for all $x\in \R^d$.
 
 \ms
 
 The function $h_\nu$ will play a very important role in the future, so it is important to understand its purpose. 
The original idea (see \cite{KKPT}, \cite{DP}) is the following: we dream to get rid of the supremum in the quantity $\wt N_{a,q}(v)(x)$, so we would like to replace $\wt N_{a,q}(v)(x)$ by $v_{\W,a,q}(x,h_{\nu}(x))$ for some `good' function $h_\nu$. 
The function $h_{\nu}$ is defined so that $h_{\nu}$ captures, roughly, the level sets of $\wt N_{a,q}(v)$ (see Lemma \ref{lem6.3}). 
In Lemma \ref{lem6.4}, we prove local $L^q$ estimates on $u_{\W,a,q}(.,h_{\nu}(.))$ when $u\in W^{1,2}_{loc}(\Omega)$ is a weak solution to $\cL u = 0$.
Lemma \ref{lem6.5} gathers Lemmas \ref{lem6.3} and \ref{lem6.4} to prove, some good $\lambda$-inequality for the bound $\wt N_{a,q} \lesssim S_{a,q}$; Lemma \ref{lem6.8} transforms this weak estimate into one in $L^p$, $p>q$; Lemma \ref{lem6.9} proves that Lemma \ref{lem6.8} can self-improve, which allows us to get an estimate in $L^q$. Lemma \ref{lem6.10} is the aim of the section, and is an easy consequence of Lemma \ref{lem6.8} and Lemma \ref{lem6.4}.
 
 \ms

 Let us start to give some basic properties of the function $h_{\nu}$.
 
 \begin{lem} \label{lem6.1}
Let $a>0$. Let $v$ be a continuous and compactly supported function. 
Choose a positive number $\nu$. Then the following properties hold.
 \begin{enumerate}[(i)]
\item The function $h_{\nu}=h_{\nu,a}(v)$ is Lipschitz with Lipchitz constant $a^{-1}$, that is for $x,y\in \R^d$,
$$|h_{\nu}(x)-h_{\nu}(y)| \leq  \frac{|x-y|}a.$$
\item For an arbitrary $x\in \{y\in \R^d:\, \sup_{\Gamma_a(y)} v > \nu\}$, we set $r_x:=h_{\nu}(x)>0$. Then there exists a point $(z,r_z) \in \dr \Gamma_a(x,r_x)$ such that $v(z,r_z) = \nu$ and $h_{\nu}(z)=r_z$.
 \end{enumerate}
 \end{lem}
 
 \bp
 Let us prove $(i)$. Pick a pair of points $x,y \in \R^d$ and set $r_x = h_{\nu}(x)$ and $r_y  = h_{\nu}(y)$. Without loss of generality, we can assume that $r_x < r_y$. We want to argue by contradiction, hence we also assume
\begin{equation} \label{contradiction}
|x-y| < a(r_y-r_x),
\end{equation}
which can be rewritten $(y,r_y)\in \Gamma_a(x,r_x)$. As a consequence, $\overline{\Gamma_a(y,r_y)}$ is a subset of $\Gamma_a(x,r_x)$. Now, it is easy to improve the inclusion $\overline{\Gamma_a(y,r_y)} \subset \Gamma_a(x,r_x)$ to 
$$\overline{\Gamma_a(y,r_y-\eta)} \subset \Gamma_a(x,r_x+\eta),$$ with $\eta>0$ small enough. We deduce that
$$\sup_{\Gamma_a(y,r_y - \eta)} v \leq \sup_{\Gamma_a(x,r_x+\eta)} v < \nu,$$
the last inequality being true by the definition of $r_x$. It follows that
$$r_y - \eta \geq h_{\nu}(y) = r_y,$$
which is a contradiction. We deduce that \eqref{contradiction} is false and thus Part $(i)$ of the lemma follows.

\medskip

We turn now to the proof of Part $(ii)$. For any $x\in \R^d$, the function $r \mapsto \ds \sup_{\Gamma_a(x,r)} v$ is continuous (and non-increasing). By definition of $r_x :=h_{\nu}(x)$, for any $m\in \N$ big enough, the quantity $\ds \sup_{\Gamma_a(x,r_x-\frac1m)} v$ is bigger than or equal to $\nu$, and the aformentioned continuity implies then the existence of $Z_m \in \Gamma_a(x,r_x-\frac1m)$ such that $v(Z_m) = \nu$. 
 Besides, again by definition of $r_x$, none of the $Z_m$ lies in $\Gamma_a(x,r_x)$.
Obviously, the sequence $(Z_m)_m$ lies in the compact set $\supp \, v_\W$, so $Z_m$ has at least one accumulation point $Z$. Such accumulation point $Z$ has to lie in $\bigcap_{m} \Gamma_a(x,r_x-\frac1m)  \setminus \Gamma_a(x,r_x)  = \dr \Gamma_a(x,r_x) $, and has to satisfy $v(Z) = \nu$ by the continuity of $v$. Hence, $r_z = h_{\nu,a}(z)$.
\ep

\ms

The function $h_\nu$ has the following interesting property.

\begin{lem} \label{lem6.3}
Let $a>0$ and $q\in (1,+\infty)$. Choose a function $v \in L^q_{loc}(\Omega)$, a smooth function $\Psi$ which satisfies $0\leq \Psi \leq 1$, and $k>2$. For any $\nu >0$ and any point $x$ satisfying $\wt N_{a,q}(v|\Psi^k)(x) >\nu$, one has
\begin{equation} \mathcal M\left[ \left( \fint_{y\in B_{ah_\nu(\cdot)/2}(\cdot)} \int_{s \in \R^{n-d}} |v|^q \Psi^k \dr_r[-\chi_\nu^k] \frac{ds}{|s|^{n-d-1}} \, dy \right)^\frac1q \right](x) \geq c\nu,
\end{equation}
where $c$ depends only on $k$ and $n$, $h_\nu := h_{\nu,a}((v|\Psi^k)_{\W,a,q})$, and $\chi_\nu$ is a cut-off function defined as
\begin{equation}\label{eq:NSlncutoff}
	\chi_{\nu}(y,s) = \phi \left(\frac{|s|}{h_\nu(y)} \right), \quad \text{ where } \phi(r) = \left\{\begin{array}{ll}  
1 & \text{ if } 0\leq r<\frac15 \\
\frac{25}{24} - \frac{5}{24} r & \text{ if } \frac15 \leq r \leq 5 \\
0 & \text{ if } r>5
\end{array}\right.
\end{equation} 
and $\chi_\nu(y,\cdot) \equiv 0$ if $h_\nu(y) = 0$
\end{lem}

\begin{rmk}
	\begin{enumerate}
		\item Since $\dr_r[-\chi_\nu^k] \approx 1/h_\nu(y)$ is supported in $\{(y,s), h_\nu(y)/5\leq |s| \leq 5h_\nu(y)\}$, and for $y\in B_{ah_\nu(x)/2}(x)$ we have $h_\nu(y) \approx h_\nu(x)$, we have that (roughly speaking)
			\[ \left( \fint_{y\in B_{ah_\nu(x)/2}(x)} \int_{s \in \R^{n-d}} |v|^q \Psi^k \dr_r[-\chi_\nu^k] \frac{ds}{|s|^{n-d-1}} \, dy \right)^\frac1q \approx (v|\Psi^k)_{W,a,q}(x,h_\nu(x)). \]
		\item Since we define the non-tangential maximal function by the average of $|v|^q\Psi^k$, instead of the pointwise values, the assumption $\wt N_{a,q}(v|\Psi^k)(x)>\nu$ can only give us information about the average value of $|v|^q\Psi^k$ on the level of $h_\nu(y)$. Therefore we use a function $\chi_\nu$ such that $\dr_r[-\chi_\nu^k]$ is supported in a band of width $\approx h_\nu(y)$. 
	\end{enumerate}
\end{rmk}

\bp
First, observe that $v_{W,a,q}$ is the $L^q$ average of a locally $L^q$ function, so $(v|\Psi^k)_{\W,a,q}$ is continuous and the function $h_\nu:= h_{\nu,a}((v|\Psi^k)_{\W,a,q})$ is well defined.

\ms

Fix $x\in \R^d$ be such that $\wt N_{a,q}(v|\Psi^k)(x) >\nu$. Set $r_x := h_\nu(x)>0$.  
From part $(ii)$ of Lemma \ref{lem6.1}, there exists $Z=(z,r_z) \in \dr \Gamma_a(x,r_x) \subset \Gamma_a(x)$ such that $(v|\Psi^k)_{W,a,q}(z,r_z) = \nu$ and $ h_{\nu}(z) = r_z$.
Define now $B$ as the ball $B_{ar_z/2}(z) \subset \R^d$.
Since $Z$ belongs to the boundary of $\Gamma_a(x,r_x)$, we deduce that $|x-z| = a (r_z-r_x) < a r_z $ and thus 
\begin{equation} \label{xin2R}
x\in 2B.
\end{equation}
Besides, due to Lemma \ref{lem6.1}, the function $h_{\nu}$ is $a^{-1}$-Lipschitz and we obtain that 
\begin{equation} \label{boundhnua2}
0 < \frac{r_z}2 \leq r_y:= h_{\nu}(y) \leq \frac{3r_z}{2} \qquad \text{ for $y\in B$,}
\end{equation}
and in particular
$$ \frac{r_z}{2} \geq \frac{r_y}{3} \quad \text{ and } \quad 2r_z \leq 4r_y \qquad \text{ for $y\in B$}.$$
The inequalities above imply that $W_a(z,r_z) \subset \{(y,s)\in \Omega, \ y\in B, \ \frac13 h_\nu(y) \leq |s| \leq 4h_\nu(y)\}$. Besides, notice that by defintion of $\chi_\nu$,
$$\dr_r[-\chi_\nu^k] \gtrsim \frac{1}{r_z} \quad \text{on } \{(y,s)\in \Omega, \ y\in R, \frac13 h_\nu(y) \leq |s| \leq 4h_\nu(y)\} \supset W(z,r_z).$$
It yields that
\begin{align}
\fint_{y\in B} \int_{s \in \R^{n-d}} |v|^q \Psi^k \dr_r[-\chi_\nu^k] \frac{ds}{|s|^{n-d-1}} \, dy \gtrsim |(v|\Psi^k)_{\W,a,q}(z,r_z)|^q = \nu^q.
\end{align}
(It is also easy to see that the left hand side is finite).
For any $y\in B$, we define $B_y = B_{ar_y/100}(y)$. By Vitali's covering lemma, we can find $N$ points $(y_i)_{i\leq N}$ such that the balls $B_{y_i}$ are non-overlapping are $ \bigcup 5B_{y_i} \supset B$. We deduce first that, since $r_y \simeq r_z$ for any $y\in B$, the value $N$ depends only on the dimension $d$, and second that there exists at least one $i\leq N$ such that 
$$ \fint_{y\in 5B_{y_i}}  \int_{s \in \R^{n-d}} |v|^q \Psi^k \dr_r[-\chi_\nu^k] \frac{ds}{|s|^{n-d-1}} dy \gtrsim \fint_{y\in B} \int_{s \in \R^{n-d}} |v|^q \Psi^k \dr_r[-\chi_\nu^k] \frac{ds}{|s|^{n-d-1}} \, dy \gtrsim \nu^q.$$
For any $y' \in B_{ar_z/10}(y_i) \subset 2B$, we let $r_{y'} = h_\nu(y')$. By Lemma \ref{lem6.1} and \eqref{boundhnua2} we get 
\[ \frac25 r_z \leq r_{y'} \leq \frac85 r_z, \]
and thus $r_{y'} \simeq r_z \simeq r_{y_i}$. And moreover, simple computations show that $5B_{y_i} = B_{ar_{y_i}/20}(y_i) \subset B_{ar_{y'}/2 }(y') $.
Therefore,
\[\begin{split}
\nu & \lesssim \left(\fint_{y\in 5B_{y_i}}  \int_{s \in \R^{n-d}} |v|^q \Psi^k \dr_r[-\chi_\nu^k] \frac{ds}{|s|^{n-d-1}} \, dy \right)^\frac1q \\
& \lesssim \fint_{y'\in B_{ar_z/10}(y_i)} \left(\fint_{y\in B_{ar_{y'}/2}(y')}  \int_{s \in \R^{n-d}} |v|^q \Psi^k \dr_r[-\chi_\nu^k] \frac{ds}{|s|^{n-d-1}} \, dy \right)^\frac1q dy' \\
& \lesssim \fint_{y'\in 2B} \left( \fint_{y\in B_{ar_{y'}/2}(y')} \int_{s \in \R^{n-d}} |v|^q \Psi^k \dr_r[-\chi_\nu^k] \frac{ds}{|s|^{n-d-1}} \, dy \right)^\frac1q \, dy'.
\end{split}\]
Since $x\in 2B$, we conclude that
$$\mathcal M \left[ \left( \fint_{y\in B_{ah_\nu(.)/2}(.)} \int_{s \in \R^{n-d}} |v|^q \Psi^k \dr_r[-\chi_\nu^k] \frac{ds}{|s|^{n-d-1}} \, dy \right)^\frac1q \right](x) \gtrsim \nu.$$
The lemma follows.
\ep

\begin{lem} \label{lem6.4}
Let $\cL=-\div |t|^{d+1-n} \cA \nabla$ be an elliptic operator that satisfies $(\cH^1_\kappa)$. Let $a>0$ and $q\in (q_0,q'_0)$ where $q_0$ is given by Proposition \ref{prop1.1}. Choose $k>2$ and a  function $\Psi$ compactly supported in $\Omega$ and that satisfies ($H^2_{a,M}$). Then for any weak solution $u \in W^{1,2}_{loc}(\Omega)$ to $\cL u = 0$, we have
\begin{equation} \label{lem6.4c}\begin{split}
\left|\iint_\Omega |u|^q \dr_r [\Psi^k] \, \frac{dt}{|t|^{n-d-1}} \, dx  \right| & \lesssim \|S_{a,q}(u|\Psi^k)\|_{q}^q + \|S_{a,q}(u|\Psi^{k-2})\|_q^{q/2} \|\wt N_{a,q}(u|\Psi^{k})\|_q^{q/2},
\end{split}\end{equation}
where the constant depends only on $a$, $q$ and $n$, $\|\cA\|_\infty$, $\|b^{-1}\|_\infty$, $\kappa$ and $M$.
\end{lem}

\begin{rmk}
A careful reader will notice that the constants don't depend on the ellipticity constants, or on an upper bound on $|b|$. In addition, we don't need for the proof of this lemma to assume that $b$ and $\cB_3$ are real valued. 

Observe also that if $u$ is a constant function (and since $\Psi$ has to be compactly supported in $\Omega$), both the right-hand term and the left-hand term in \eqref{lem6.4c} are zero.
\end{rmk}

\bp 
First, Proposition \ref{lemMoserI} ensures that $|u|^q$ and $|\nabla u|^2 |u|^{q-2}$ both lie in $L^1_{loc}(\Omega)$, and so all the quantities in \eqref{lem6.4c} are well defined and finite.

We also use the same trick as in Lemma \ref{lemS<Na}, which says that $u$ is also a weak solution to $\cL' u = 0$ where $\cL' = - \div |t|^{d+1-n} \cA' \nabla + |t|^{d+1-n} \cD'\cdot \nabla$,
\begin{equation}\label{defA'}
\ds \cA' = \begin{pmatrix} \cA_{1} & \cA_2 + (\cB_3)^T \\ 0 & bI \end{pmatrix} + \cC,
\end{equation}
and for any $1\leq j \leq d < i \leq n$ 
\begin{equation}\label{defB'}
(\cD')_j = - \sum_{k>d} |t|^{n-d-1} \dr_{t_k} [(\cB_3)_{kj} |t|^{d+1-n}] \quad \text{ and } \quad (\cD')_i = \sum_{k\leq d} \dr_{x_k} (\cB_3)_{ik}.
\end{equation}
Note in particular that $|t|\cD'$ satisfies the Carleson measure condition. Unless we assume $\cB_3$ is real-valued, nothing guarantees here that the matrix $\cA'$ is elliptic, but we shall not use this assumption.

\ms

Let us denote 
$$T:= \iint_\Omega |u|^q \dr_r [\Psi^k] \, \frac{dt}{|t|^{n-d-1}} \, dx.$$
By using the product rule and the fact that $\dr_{r}[ |u|^q ] = \frac q2 |u|^{q-2} [\bar u \, \dr_{r} u + u \, \dr_{r} \bar u]$, one obtains
\[\begin{split}
T & = \iint_\Omega \dr_r[|u|^q\Psi^k] \, \frac{dt}{|t|^{n-d-1}} \, dx \\
& \quad - \frac q2  \iint_\Omega  \dr_r[u] |u|^{q-2}\bar u   \Psi^k \, \frac{dt}{|t|^{n-d-1}} \, dx - \frac q2  \iint_\Omega  \dr_r[\bar  u] |u|^{q-2} u  \Psi^k \, \frac{dt}{|t|^{n-d-1}} \, dx \\
& \qquad := T_1+T_2 + T_3.
\end{split}\]
The term $T_1$ is 0. Indeed, we switch to cylindrical coordinates, and we have
$$T_1 = \int_{x\in \R^d} \int_{\theta \in \mathbb S_{n-d-1}} \int_0^\infty \dr_r[|u|^q\Psi^k] \, dr\, d\theta\, dx = 0,$$
since $\Psi$ is compactly supported in $\Omega$.

We need to bound $T_2$ and $T_3$. Since, $T_3 = \overline{T_2}$, we only need to treat $T_2$. We remark that the estimate of $T_2$ is in some sense the reverse of the estimate of the term $T_4$ in Lemma \ref{lemS<Na}. Observe that
\[
T_2 = -\frac q2  \iint_{(x,t)\in \Omega} b \nabla_t u \cdot \nabla_t |t| \,  \frac{|u|^{q-2}\bar u \Psi^{k}}{b} \frac{dt}{|t|^{n-d-1}} \, dx,\]
and so with \eqref{defA'}, we have
\[\begin{split}
T_2 & = - \frac q2  \iint_{(x,t)\in \Omega} \cA' \nabla u \cdot \nabla |t| \,  \frac{|u|^{q-2} \bar u \Psi^{k}}{b} \frac{dt}{|t|^{n-d-1}} \, dx \\
& \quad  + \frac q2 \iint_{(x,t)\in \Omega} \cC \nabla u \cdot  \nabla |t| \,  \frac{|u|^{q-2} \bar u \Psi^{k}}{b}  \frac{dt}{|t|^{n-d-1}} \, dx \\
& := T_{21} + T_{22}.
\end{split}\]
We use the fact that $|b| \gtrsim 1$, and then the Cauchy-Schwarz's inequality to bound $T_{22}$ as follows
\[\begin{split}
|T_{22}| & \lesssim  \iint_{(x,t)\in \Omega}  |\cC|  |\nabla u| \, \frac{|u|^{q-1} \Psi^k}{|b|}  \frac{dt}{|t|^{n-d-1}} \, dx \\
& \lesssim  \iint_{(x,t)\in \Omega}  |\cC|  |\nabla u| \, |u|^{q-1} \Psi^k  \frac{dt}{|t|^{n-d-1}} \, dx \\
& \lesssim \left(\iint_{(x,t)\in \Omega} |\nabla u|^2 |u|^{q-2} \Psi^k \, \frac{dt}{|t|^{n-d-2}} \, dx \right)^\frac12 \left( \iint_{(x,t)\in \Omega} |\cC|^2 |u|^q \Psi^k \frac{dt}{|t|^{n-d}} \, dx  \right)^\frac12.
\end{split}\]
The first integral in the right-hand side above can be bounded by $\|S_{a,q}(u|\Psi^k)\|_q^q \leq \|S_{a,q}(u|\Psi^{k-2})\|_q^q$. Since $\cC$ satisfies the Carleson measure condition, we use the Proposition \ref{propCarl} (Carleson inequality) to bound the second integral as follows
\[\begin{split} 
\iint_{(x,t)\in \Omega} |\cC|^2 |u|^q \Psi^k \frac{dt}{|t|^{n-d}} \, dx & \leq \kappa \|\wt N_{a,q}(u|\Psi^k)\|_q^q.
\end{split}\]
As a consequence,
\[\begin{split}
|T_{22}| & \lesssim \|S_{a,q}(u|\Psi^{k-2})\|_q^{q/2} \|\wt N_{a,q}(u|\Psi^k)\|_q^{q/2}.
\end{split}\]

Now, we deal with $T_{21}$. We write
\[\begin{split}
T_{21} & =  - \frac q2  \iint_{(x,t)\in \Omega} \cA' \nabla u \cdot \nabla \left[\frac{|t||u|^{q-2}\bar u \Psi^k}{b}\right] \frac{dt}{|t|^{n-d-1}} \, dx \\
& \quad +  \frac q2  \iint_{(x,t)\in \Omega} \cA' \nabla u \cdot \nabla[ |u|^{q-2} \bar u] \frac{\Psi^k}{b} \frac{dt}{|t|^{n-d-2}} \, dx \\
& \quad +  \frac q2  \iint_{(x,t)\in \Omega} \cA' \nabla u \cdot \nabla [\Psi^k] \frac{|u|^{q-2} \bar u}{b} \frac{dt}{|t|^{n-d-2}} \, dx \\
& \quad -  \frac q2  \iint_{(x,t)\in \Omega} \cA' \nabla u \cdot \nabla b \frac{|u|^{q-2} \bar u \Psi^k}{b^2} \frac{dt}{|t|^{n-d-2}} \, dx \\
& \quad := T_{211} + T_{212} + T_{213} + T_{214}.
\end{split}\]
By using the boundedness of $\cA'/b$, the term $T_{212}$ can be bounded as follows 
\[
|T_{212}| \lesssim \iint_{(x,t) \in \Omega} |\nabla u| |\nabla [|u|^{q-2} \bar u]| \Psi^k \frac{dt}{|t|^{n-d-2}} \, dx.
\] 
But using the product rule and the fact that $|\nabla |u|| \leq |\nabla u|$, one has $|\nabla [|u|^{q-2} \bar u]| \leq (q-1) |u|^{q-2} |\nabla u|$ if $q\leq  2$ and $|\nabla [|u|^{q-2} \bar u]| \leq (3-q) |u|^{q-2} |\nabla u|$  if $q\geq 2$. Hence one can bound
\[
|T_{212}| \lesssim \iint_{(x,t) \in \Omega} |\nabla u|^2 |u|^{q-2} \Psi^k \frac{dt}{|t|^{n-d-2}} \, dx \lesssim \|S_{a,q}(u|\Psi^k)\|_{q}^q.
\]
The terms $T_{213}$ and $T_{214}$ are bounded in a similar manner as $T_{22}$. Let us start with $T_{213}$. Since $\cA'/b$ is uniformly bounded,
\[\begin{split}
|T_{213}| & \lesssim \iint_{(x,t)\in \Omega} |\nabla u| |\nabla [\Psi^k]| |u|^{q-1} \frac{dt}{|t|^{n-d-2}} \, dx \lesssim \iint_{(x,t)\in \Omega} |\nabla u| |\nabla \Psi| \Psi^{k-1} |u|^{q-1} \frac{dt}{|t|^{n-d-2}} \, dx \\
& \lesssim \left( \iint_{(x,t)\in \Omega} |\nabla u|^{2} |u|^{q-2} \Psi^{k-2} \frac{dt}{|t|^{n-d-2}} \, dx \right)^\frac12 \left( \iint_{(x,t)\in \Omega} |u|^{q} \Psi^{k} |\nabla \Psi|^2 \frac{dt}{|t|^{n-d-2}} \, dx \right)^\frac12
\end{split}\] 
by Cauchy-Schwarz's inequality. The first integral in the right-hand side above is bounded by $\|S_{a,q}(u|\Psi^{k-2})\|_q^q$. As for the second one, Lemma \ref{lemH2} gives that
$$\iint_{(x,t)\in \Omega} |u|^{q} \Psi^{k} |\nabla \Psi|^2 \frac{dt}{|t|^{n-d-2}} \, dx \lesssim M \|\wt N_{a,q}(u|\Psi^k)\|_q^q.$$
We conclude that
$$|T_{213}| \lesssim \|S_{a,q}(u|\Psi^{k-2})\|_q^{q/2} \|\wt N_{a,q}(u|\Psi^k)\|_q^{q/2}.$$

Then we deal with $T_{214}$. Using the fact that $\cA'/b$ is uniformly bounded, Cauchy-Schwarz's inequality and that fact that $|t|\nabla b$ satisfies the Carleson measure condition, similarly to the bound on $T_{22}$ or $T_{213}$, we obtain
\begin{equation} \label{T114}\begin{split}
|T_{214}| & \lesssim \iint_{(x,t)\in \Omega} |\nabla u| |\nabla b| |u|^{q-1} \Psi^k \frac{dt}{|t|^{n-d-2}} \, dx \\
& \lesssim \left( \iint_{(x,t)\in \Omega} |\nabla u|^{2} |u|^{q-2} \Psi^{k} \frac{dt}{|t|^{n-d-2}} \, dx \right)^\frac12 \left( \iint_{(x,t)\in \Omega} (|t||\nabla b|)^2 |u|^{q} \Psi^{k} \frac{dt}{|t|^{n-d}} \, dx \right)^\frac12 \\
& \lesssim \kappa^{\frac12} \|S_{a,q}(u|\Psi^{k})\|_q^{q/2} \|\wt N_{a,q}(u|\Psi^k)\|_q^{q/2} \\
& \lesssim \|S_{a,q}(u|\Psi^{k-2})\|_q^{q/2} \|\wt N_{a,q}(u|\Psi^k)\|_q^{q/2} \\
\end{split}\end{equation}

\ms

It remains to bound $T_{211}$. However, $T_{211}$ can be treated as the term $T_1$ in Lemma \ref{lemS<Na}, by using the fact that $u$ is a weak solution to $\cL'u=0$, and we will eventually obtain that
$$|T_{211}| \lesssim \left| \iint_{(x,t)\in \R^n} \cD'\cdot \nabla u \, \frac{|t| |u|^{q-2} \bar u \Psi^k}{b} \frac{dt}{|t|^{n-d-1}} dx \right|\lesssim \iint_{(x,t)\in \R^n} |t||\cD'| |\nabla u| |u|^{q-1}\Psi^{k} \frac{dt}{|t|^{n-d-1}} dx . $$
However, recall that $|\cD'| \lesssim |\nabla \cB_3|$ so $|t|\cD'$ satisfies the Carleson measure condition. By using Cauchy-Schwarz's inequality and then Carleson's inequality (Proposition \ref{propCarl}), similarly to $T_{22}$,
$$|T_{111}| \lesssim \kappa^{\frac12} \|S_{a,q}(u|\Psi^{k-2})\|_q^{q/2} \|\wt N_{a,q}(u|\Psi^k)\|_q^{q/2}.$$
The lemma follows.
\ep

\begin{lem} \label{lem6.5}
Let $\cL$ be an elliptic operator satisfying ($H^1_\kappa$) for some constant $\kappa\geq 0$. Let $a, l>0$ and $q\in (q_0,q'_0)$ where $q_0$ is given by Proposition \ref{prop1.1}. Choose $k>2$, a positive $a^{-1}$-Lipschitz function $e$ and a ball $B:=B_{l'}(x_B) \subset \R^d$ of radius $l'\geq l$. Then there exists $\eta \in (0,1)$ that depends only on $d$ such that for any weak solution $u\in W^{1,2}_{loc}(\Omega)$, any $\nu >0$ and any $\gamma \in (0,1)$,
\begin{equation} \label{lem6.5a} \begin{split}
|\{x\in \R^d, \ \wt N_{a,q}(u|\Psi_e^k\Psi_{B,l}^k)(x) > \nu\} \cap E_{\nu,\gamma}| \leq C\gamma^q |\{x\in \R^d, \, \mathcal M[\wt N_{a,q}(u|\Psi_e^k\Psi_{B,l}^k)](x) > \eta \nu\}|,
\end{split} \end{equation}
where $\Psi_e, \Psi_{B,l}$ are defined as in Lemma \ref{exgoodcutoff}, and
\begin{equation} \label{defEnu} \begin{split}
E_{\nu,\gamma}:= & \Bigg\{x\in \R^d, \, \mathcal M_q\left[\left(\fint_{y\in B_{ae(.)/2}(.)}\int_{s\in \R^{n-d}} |u|^q \Psi_{B,l}^k \dr_r [\Psi_e^k] \frac{ds}{|s|^{n-d-1}} \right)^\frac1q \right](x) \leq \gamma \nu \\
& \hspace{2cm}  \text{ and } \mathcal M_q[\wt N_{a,q}(u|\Psi_e^k\Psi_{B,l}^k)](x)\mathcal M_q[S_{a,q}(u|\Psi_e^{k-2}\Psi_{B,l}^{k-2})](x)  \leq \gamma^2 \nu^2\Bigg\}
\end{split}\end{equation}
and the constant $C>0$ depends on $a$, $q$, $n$, $\|\cA\|_\infty$, $\|b^{-1}\|_\infty$, $\kappa$ and $k$.
\end{lem}

\bp Tak some $\eta$ to be fixed later. Choose $\nu>0$. We write to lighten the notation $\Psi$ for $\Psi_e \Psi_{B,l}$. We define
$$\mathcal S:= \{x\in \R^d, \, \mathcal M[\wt N_{a,q}(u|\Psi^k)](x) > \eta\nu\},$$
which is open and bounded. Indeed, $\mathcal S$ is bounded because $\Psi$ is compactly supported and then $\wt N_{a,q}(u|\Psi^k)] \equiv 0$ outside a big ball, and $\mathcal S$ is open because $(u|\Psi^k)_\W$ is continuous.

We construct a Whitney decomposition as follows (the construction is classical, and we only aim to prove that \eqref{prWhitneyZ} is possible). For any $x\in \mathcal S$, we set $B_x \subset \R^d$ the ball of center $x$ and radius $\dist(x,\mathcal S^c)/10)$. The balls $B_x$ have uniformly bounded radius (because $|B_x| \leq |\mathcal S| <+\infty$) and therefore Vitali's covering lemma entails the existence of a non-overlapping collection of balls $(B_{x_i})_{i\in I}$ such that $\bigcup_{i\in I} 5B_{x_i} = \mathcal S$.
We write $B_i$ for $10B_{x_i}$ and $l_i$ for its radius, so by construction,
\begin{equation} \label{prWhitneyZ}
\bigcup_{i\in I} B_{i} = \mathcal S \quad \text{ and } \quad \sum_{i\in I} |B_i| \leq 10^d |\mathcal S| \quad \text{ and for $i\in I$, there exists $y_i$ in $\overline{B_i} \cap \mathcal S^c$.}
\end{equation}
The point $y_i$ satisfies then
\begin{equation} \label{defyi2}
|x_i - y_i| = l_i \qquad \text{ and } \quad \mathcal M[\wt N_{a,q}(u|\Psi^k)](y_i) \leq \eta \nu.
\end{equation}

\ms

Since $\sum_{i\in I} |B_i| \lesssim |\mathcal S|$, the estimate \eqref{lem6.5a} will be proven if we establish that
\begin{equation} \label{claimy}
F_{\gamma}^i := \{x\in \R^d, \, \wt N_{a,q}(u|\Psi^k)(x) > \nu\} \cap E_{\nu,\gamma} \cap B_i \leq C\gamma^q |B_i|,
\end{equation}
where $C$ is independent of $\gamma \in (0,1)$. If $F^i_\gamma = \emptyset$, there is nothing to prove, so we can assume that $F^i_\gamma$ contains some point $z_i$.  

Similar to Lemma \ref{exgoodcutoff} (2), we denote
\begin{equation}\label{def:Phii}
	\Phi_i = \Psi_{3B_i,2l_i}(x,t) = \phi\left( \frac{a|t|}{2 l_i} \right) \phi \left( 1+ \frac{\dist(x,3B_i)}{2l_i} \right). 
\end{equation} 
We claim that we can find $\eta$ small enough such that
\begin{equation} \label{claimz}
\wt N_{a,q}(u|\Psi^k\Phi_i^k)(z) > \eta \qquad \text{ for } z\in F^i_\gamma.
\end{equation}
Indeed, let $z\in B_i$, for any $(z',r') \in \Gamma_a(z)$ such that $r'\geq l_i/a$,
$$(u|\Psi^k)_{\W,a,q}(z',r') \leq \wt N_{a,q}(u|\Psi^k)(z'') \qquad \text{ for } z'' \in B_{ar'}(z').$$
Now, since 
\[ |z' - y_i| \leq |z'-z| + |z-y_i| < ar' + 2l_i \leq 3ar', \] we deduce
\[\begin{split}
(u|\Psi^k)_{\W,a,q}(z',r') & \leq \fint_{B_{ar'}(z')} \wt N_{a,q}(u|\Psi^k)(z'') dz'' \leq C_d \fint_{B_{3ar'}(z')} \wt N_{a,q}(u|\Psi^k)(z'') dz'' \\
&  \leq C_d \mathcal M[\wt N_{a,q}(u|\Psi^k)](y_i) \leq C_d \eta \nu
\end{split}\]
by \eqref{defyi2}. We choose $\eta$ such that $C_d\eta \leq 1$ and we obtain that 
\begin{equation}\label{eq:NStmp1}
	(u|\Psi^k)_{\W,a,q}(z',r') \leq \nu \qquad \text{ for } (z',r') \in \Gamma_a(z), \, r'\geq l_i/a.
\end{equation} 
We observe that $\Phi_i \equiv 1$ on $W_{a}(z'',r'')$ if $(z'',r'')\in \Gamma_a(z)$ and $r'' \leq l_i/a$, and thus 
\[ (u|\Psi^k)_{\W,a,q}(z'',r'') = (u|\Psi^k\Phi_i^k)_{\W,a,q}(z'',r''). \]
Recall that $\wt N_{a,q}(u|\Psi^k)(z) > \nu$ for $z\in F^i_\gamma$ and \eqref{eq:NStmp1}, so we necessarily have
$$\wt N_{a,q}(u|\Psi^k\Phi^k_i)(z) >\nu.$$
The claim \eqref{claimz} follows.

\ms

We are now ready to use Lemma \ref{lem6.3}. Set $h_\nu := h_{\nu,a}((u|\Psi^k \Phi_i^k)_{\W,a,q})$ and $\chi_\nu$ as in \eqref{eq:NSlncutoff}.
Lemma \ref{lem6.3} gives that for any $z\in F^i_\gamma$
$$
\mathcal M\left[ \left( \fint_{y\in B_{ah_\nu(.)/2}(.)} \int_{s \in \R^{n-d}} |u|^q \Psi^k \Phi_i^k \dr_r[-\chi_\nu^k] \frac{ds}{|s|^{n-d-1}} \, dy \right)^\frac1q \right](z) \gtrsim \nu,
$$
and so if we integrate on $z\in \R^d$,
\begin{equation}
|F^i_\gamma| \lesssim \frac1{\nu^q} \int_{z\in \R^d} \left|\mathcal M\left[ \left( \fint_{y\in B_{ah_\nu(.)/2}(.)} \int_{s \in \R^{n-d}} |u|^q \Psi^k \Phi_i^k \dr_r[-\chi_\nu^k] \frac{ds}{|s|^{n-d-1}} \, dy \right)^\frac1q \right](z) \right|^q dz.
\end{equation}
The Hardy-Littlewood maximal inequality entails that
\begin{equation} \label{5.18}
|F^i_\gamma| \lesssim \frac1{\nu^q} \int_{z\in \R^d} \fint_{y\in B_{ah_\nu(z)/2}(z)} \int_{s \in \R^{n-d}} |u|^q \Psi^k \Phi_i^k \dr_r[-\chi_\nu^k] \, \frac{ds}{|s|^{n-d-1}} \, dy \,  dz.
\end{equation}
If $y\in B_{ah_\nu(z)/2}(z)$, then $|y-z| \leq \frac12 ah_\nu(z)$ and since $h_\nu$ is $a^{-1}$-Lipschitz,
\begin{equation} \label{hnuyz}
\frac{1}{2} h_\nu(z) \leq h_\nu(y) \leq \frac32 h_\nu(z),
\end{equation}
in other words $h_\nu(z) \leq 2 h_\nu(y)$. It follows that $z \in B_{ah_\nu(y)}(y)$ and so by Fubini's theorem,
\begin{equation} \label{5.20}\begin{split}
|F^i_\gamma| & \lesssim \frac1{\nu^q} \int_{y\in \R^d} \int_{s \in \R^{n-d}} |u|^q \Psi^k \Phi_i^k \dr_r[-\chi_\nu^k] \, \frac{ds}{|s|^{n-d-1}} \,  \left(\int_{z\in B_{ah_\nu(y)}(y)} (ah_\nu(z))^{-d} dz\right) dy \\
& \lesssim \frac1{\nu^q} \int_{y\in \R^d} \int_{s \in \R^{n-d}} |u|^q \Psi^k \Phi_i^k \dr_r[-\chi_\nu^k] \, \frac{ds}{|s|^{n-d-1}} \, dy, \\
\end{split} \end{equation}
by \eqref{hnuyz}. 

\bigskip

Now we estimate the right hand side of \eqref{5.20}. By the product rule,
\[\begin{split}
\Psi^k \Phi_i^k \dr_r [-\chi_\nu^k] & = - \dr_r [\Psi^k \Phi_i^k\chi_\nu^k] + \chi_\nu^k \dr_r [\Psi^k \Phi_i^k] = - \dr_r [\Psi^k \Phi_i^k\chi_\nu^k] + \chi_\nu^k (\Psi_e^k \dr_r [\Psi_{B,l}^k \Phi_i^k] + \Psi_{B,l}^k \Phi_i^k \dr_r[\Psi_e^k]) \\
& \leq - \dr_r [\Psi^k\Phi_i^k\chi_\nu^k] + \chi_\nu^k \Psi_{B,l}^k \Phi_i^k \dr_r[\Psi_e^k],
\end{split}\]
where the last line holds because $\Psi_{B,l} \Phi_i$ is decreasing in $r$ (because we build $\Psi_{B,l}$, $\Phi_i$ with the help of $\phi$ which is decreasing). Set 
$$T_1:= - \int_{y\in \R^d} \int_{s \in \R^{n-d}} |u|^q \dr_r[\Psi^k\Phi_i^k\chi_\nu^k] \, \frac{ds}{|s|^{n-d-1}} \, dy$$
and
$$T_2:= \int_{y\in \R^d} \int_{s \in \R^{n-d}} |u|^q \chi_\nu^k \Phi_i^k \Psi_{B,l}^k \dr_r[\Psi^k_e] \, \frac{ds}{|s|^{n-d-1}} \, dy.$$
We have
\begin{equation} \label{FbyT12}
 |F^i_\gamma| \lesssim \nu^{-q} (T_1 + T_2),
\end{equation}
and we want to bound $T_1$ and $T_2$. Let us start with $T_2$. Let $y\in \R^d$ and $z\in B_{ae(y)/4}(y)$. Since $e$ is $a^{-1}$-Lipschitz, similarly to \eqref{hnuyz}, one has
\begin{equation}\label{eq:NStmp2}
\frac{3}{4} e(y) \leq e(z) \leq \frac54 e(y).
\end{equation}
In particular, $y\in B_{ae(z)/3}(z)$ and for any $y\in \R^d$
$$1 \lesssim \int_{z\in B_{ae(y)/4}} (ae(z))^{-d} dz.$$
The bound on $T_2$ becomes then
\[\begin{split}
T_2 & \lesssim \int_{y\in \R^d} \int_{s \in \R^{n-d}} |u|^q \chi_\nu^k \Phi_i^k\Psi_{B,l}^k \dr_r[\Psi^k_e] \, \frac{ds}{|s|^{n-d-1}}  \int_{z\in B_{ae(y)/4}} (ae(z))^{-d} dz \, dy \\
& \lesssim \int_{z\in \R^d} \fint_{y\in B_{ae(z)/3}(z)} \int_{s \in \R^{n-d}} |u|^q \chi_\nu^k \Phi_i^k\Psi_{B,l}^k \dr_r[\Psi^k_e] \, \frac{ds}{|s|^{n-d-1}} \, dy \, dz
\end{split}\]
by Fubini's lemma. We want to see for which $z\in \R^d$ the quantity 
$$\fint_{y\in B_{ae(z)/3}(z)} \int_{s \in \R^{n-d}} |u|^q \chi_\nu^k \Phi_i^k\Psi_{B,l}^k \dr_r[\Psi^k_e] \, \frac{ds}{|s|^{n-d-1}} \, dy$$
is non-zero. First, by the definition \eqref{def:Phii} we know
\[ \supp \, \Phi_i \subset \{(y,s)\in \R^n, |s|\leq 4l_i/a, y\in 5B_i\}. \] 
Thus we need $e(z) \leq 10 l_i/a$, because otherwise $\Phi_i^k \dr_r[\Psi^k_e] \equiv 0$ (we also use \eqref{eq:NStmp2} here).
We also need $B_{ae(z)/3}(z) \cap 5B_i \neq \emptyset$ to guarantee $\Phi_i \neq 0$. Altogether, $z$ needs to lie in - say - $10B_i$. Recall that there exists some $z_i\in F_\gamma^i \subset 10B_i$, we conclude
\begin{equation} \label{T2forF} \begin{split}
T_2 & \lesssim \int_{z\in 10B_i} \fint_{y\in B_{ae(z)/3}(z)} \int_{s \in \R^{n-d}} |u|^q \chi_\nu^k \Phi_i^k\Psi_{B,l}^k \dr_r[\Psi^k_e] \, \frac{ds}{|s|^{n-d-1}} \, dy \, dz \\
& \lesssim \int_{z\in 10B_i} \fint_{y\in B_{ae(z)/3}(z)} \int_{s \in \R^{n-d}} |u|^q \Psi_{B,l}^k \dr_r[\Psi^k_e] \, \frac{ds}{|s|^{n-d-1}} \, dy \, dz \\
& \lesssim |10B_i| \left|\mathcal M_q\left[\left(\fint_{y\in B_{ae(.)/2}(.)}\int_{s\in \R^{n-d}} |u|^q \Psi_{B,l}^k \dr_r [\Psi_e^k] \frac{ds}{|s|^{n-d-1}} \right)^\frac1q \right](z_i) \right|^q\\
& \lesssim |B_i| \gamma^q \nu^q.
\end{split}\end{equation}

We turn now to the treatment of $T_1$. Observe that, by invoking the same argument as the one used to prove Lemma \ref{exgoodcutoff} (1), we can prove that $\chi_\nu$ satisfies ($\cH^2_{a,M}$) with $M$ that depends only on $a$ and $n, d$. (Note, in particular, that the Lipschitz constant for $h_\nu$ is $a^{-1}$, independent of $\nu$.) Therefore, Lemma \ref{lem6.4} entails that
\[\begin{split}
|T_1| & \lesssim \|S_{a,q}(u|\Psi^k\Phi_i^k\chi_\nu^k)\|_{q}^q + \|S_{a,q}(u|\Psi^{k-2}\Phi_i^{k-2}\chi_\nu^{k-2})\|_q^{q/2} \|\wt N_{a,q}(u|\Psi^k\Phi_i^k\chi_i^k)\|_q^{q/2} \\
& \lesssim \|S_{a,q}(u|\Psi^k\Phi_i^k)\|_{q}^q + \|S_{a,q}(u|\Psi^{k-2}\Phi_i^{k-2})\|_q^{q/2} \|\wt N_{a,q}(u|\Psi^k\Phi_i^k)\|_q^{q/2}.
\end{split}\]
Since $\Phi_i$ is supported in $\{(x,t)\in \Omega,\ x\in 5B_i, \, |t| \leq 4l_i/a\}$, the functions $S_{a,q}(u|\Psi^{k-2}\Phi_i^{k-2})$ and $\wt N_{a,q}(u|\Psi^k\Phi_i^k)$ are supported in - say - $10B_i$. The bound on $T_1$ becomes then
\[\begin{split}
|T_1| & \lesssim \|S_{a,q}(u|\Psi^k)\|_{L^q(10B_i)}^q + \|S_{a,q}(u|\Psi^{k-2})\|_{L^q(10B_i)}^{q/2} \|\wt N_{a,q}(u|\Psi^k)\|_{L^q(10B_i)}^{q/2} \\
& \lesssim |B_i| \left(  |\mathcal M_q[S_{a,q}(u|\Psi^k)](z_i)|^q + |\mathcal M_q[S_{a,q}(u|\Psi^k)](z_i)|^{q/2} |\mathcal M_q[\wt N_{a,q}(u|\Psi^k)](z_i)|^{q/2} \right).
\end{split}\]
Since $z_i \in F^i_\gamma$, we have 
$$ |\mathcal M_q[S_{a,q}(u|\Psi^k)](z_i)|^{q/2} |\mathcal M_q[\wt N_{a,q}(u|\Psi^k)](z_i)|^{q/2} \leq \gamma^q \nu^q.$$
In addition,
\[\begin{split}
|\mathcal M_q[S_{a,q}(u|\Psi^k)](z_i)|^q & = \left(\frac{|\mathcal M_q[S_{a,q}(u|\Psi^k)](z_i)| |\mathcal M_q[\wt N_{a,q}(u|\Psi^k)](z_i)|}{|\mathcal M_q[\wt N_{a,q}(u|\Psi^k)](z_i)|} \right)^{q} \\
& \leq \frac{\gamma^{2q} \nu^{2q} }{\nu^q}  \leq \gamma^{2q} \nu^q \leq \gamma^q \nu^q.
\end{split}\]
For the first inequality we also used the fact that $\wt N_{a,q}(u|\Psi^k)(z_i) \geq \nu$ and that $\wt N_{a,q}(u|\Psi^k)$ is continuous.
We deduce
\begin{equation} \label{T1forF}
|T_1| \lesssim |B_i| \gamma^q \nu^q.
\end{equation}
The combination of \eqref{FbyT12}, \eqref{T2forF} and \eqref{T1forF} proves \eqref{claimy}. The lemma follows.
\ep

\begin{lem} \label{lem6.8}
Let $\cL$ be an elliptic operator satisfying $(\cH^1_\kappa)$ for some constant $\kappa\geq 0$. Let $a, l>0$, $q\in (q_0,q'_0)$ where $q_0$ is given by Proposition \ref{prop1.1}, and $p>q$. Choose $k>2$, a positive $a^{-1}$-Lipschitz function $e$ and a ball $B \subset \R^d$ of radius $l'\geq l$. Then for any weak solution $u\in W^{1,2}_{loc}(\Omega)$,
\begin{equation} \label{NbySLpZ}
\|\wt N_{a,q}(u|\Psi_e^k\Psi_{B,l}^k)\|_{p}^p \lesssim \|S_{a,q}(u|\Psi_e^{k-2}\Psi_{B,l}^{k-2})\|_{p}^p + \left\|\left(\fint_{y\in B_{ae(.)/2}(.)} \int_{s\in \R^{n-d}} |u|^q \Psi_{B,l}^k \dr_r [\Psi_e^k] \frac{ds}{|s|^{n-d-1}}\, dy \right)^\frac1q \right\|_p^p
\end{equation}
where the constant depends on $a$, $q$, $n$, $\|\cA\|_\infty$, $\|b^{-1}\|_\infty$, $\kappa$, $k$.
\end{lem}

\begin{rmk}\label{rmk:pq}
	The reader may think of the last term as the average of $u$ in the local region in consideration (determined by $B$ and $e$). It appears on the right side because, roughly speaking, we are estimating $u$ by its gradient.
	
	The aforementioned region is close to boundary, and the reader may be used to see it lying in a Whitney region of the boundary ball $B$. However, we find it easier to get the self-improvement given by Lemma \ref{lem6.9} with the estimate \eqref{NbySLpZ}, and Lemma \ref{lem6.4} will be used (again) in Lemma \ref{lem6.10} to recover a ``classical'' right-hand term.
\end{rmk}

\bp As in the previous proof, we write $\Psi$ for $\Psi_e\Psi_{B,l}$. Besides, $E_{\nu,\gamma}$ denotes the set
\[\begin{split}
E_{\nu,\gamma}:= & \left\{x\in \R^d, \, \mathcal M_q\left[\left(\fint_{y\in B_{ae(.)/2}(.)}\int_{s\in \R^{n-d}} |u|^q \Psi_{B,l}^k \dr_r [\Psi_e^k] \frac{ds}{|s|^{n-d-1}} \right)^\frac1q \right](x) \leq \gamma \nu \right.\\
& \hspace{3cm}  \text{ and } \mathcal M_q[\wt N_{a,q}(u|\Psi_e^k\Psi_{B,l}^k)](x)\mathcal M_q[S_{a,q}(u|\Psi_e^{k-2}\Psi_{B,l}^{k-2})](x)  \leq \gamma^2 \nu^2\Big\}
\end{split}\]
By Lemma \ref{lem6.5}, there exists $\eta = \eta(d)$ such that for any $\gamma \in (0,1)$,
\[\begin{split}
\|\wt N_{a,q} &(u|\Psi^k) \|_p^p = c_p \int_0^\infty \nu^{p-1} |\{x\in \R^d, \, \wt N_{a,q}(u|\Psi^k)(x)>\nu\}|  d\nu \\
& \leq C \int_0^\infty \nu^{p-1} |\{x\in \R^d, \, \wt N_{a,q}(u|\Psi^k)(x)>\nu\} \cap E_{\nu,\gamma}| d\nu \\
& \ + C \int_0^\infty \nu^{p-1} \left|\left\{x\in \R^d, \, \mathcal M_q\left[\left(\fint_{y\in B_{ae(.)/2}(.)}\int_{s\in \R^{n-d}} |u|^q \Psi_{B,l}^k \dr_r [\Psi_e^k] \frac{ds}{|s|^{n-d-1}} \right)^\frac1q \right](x)> \gamma \nu\right\}\right| d\nu \\
& \ + C  \int_0^\infty \nu^{p-1} \left|\left\{x\in \R^d, \, \mathcal M^{1/2}_q[\wt N_{a,q}(u|\Psi^k)](x)\mathcal M^{1/2}_q[S_{a,q}(u|\Psi^{k-2})](x) > \gamma \nu\right\}\right| d\nu. \\
& \leq C \gamma^q \int_0^\infty \nu^{p-1} |\{x\in \R^d, \, \mathcal M_q[\wt N_{a,q}(u|\Psi^k)](x)>\eta\nu\}| d\nu \\
&  \qquad + C\gamma^{1-p} \left\|\mathcal M_q\left[\left(\fint_{y\in B_{ae(.)/2}(.)}\int_{s\in \R^{n-d}} |u|^q \Psi_{B,l}^k \dr_r [\Psi_e^k] \frac{ds}{|s|^{n-d-1}} \right)^\frac1q \right] \right\|_p^p \\
& \qquad + C\gamma^{1-p}  \left\| \mathcal M^{1/2}_q[\wt N_{a,q}(u|\Psi^k)] \ \mathcal M^{1/2}_q[S_{a,q}(u|\Psi^{k-2})]\right\|_p^p. \\
\end{split}\]
Now applying the Hardy-Littlewood maximal theorem with power $p/q>1$ to each term, and using also Cauchy-Schwarz inequality for the last term, we get
\[\begin{split}
\|\wt N_{a,q} (u|\Psi^k) \|_p^p
& \leq C \gamma^q \left\|\mathcal M_q\left[\wt N_{a,q}(u|\Psi^k) \right] \right\|_p^p \\
&  \qquad + C\gamma^{1-p} \left\|\mathcal M_q\left[\left(\fint_{y\in B_{ae(.)/2}(.)}\int_{s\in \R^{n-d}} |u|^q \Psi_{B,l}^k \dr_r [\Psi_e^k] \frac{ds}{|s|^{n-d-1}} \right)^\frac1q \right] \right\|_p^p \\
& \qquad + C\gamma^{1-p}  \left\| \mathcal M_q[\wt N_{a,q}(u|\Psi^k)]\right\|^{1/2}_p  \left\| \mathcal M^{p/2}_q[S_{a,q}(u|\Psi^{k-2})]\right\|_p^{p/2}. \\
& \leq C \gamma^q \|\wt N_{a,q} (u|\Psi^k) \|_p^p \\
& \qquad + C \gamma^{1-p} \left\|\left(\fint_{y\in B_{ae(.)/2}(.)}\int_{s\in \R^{n-d}} |u|^q \Psi_{B,l}^k \dr_r [\Psi_e^k] \frac{ds}{|s|^{n-d-1}} \right)^\frac1q\right\|_p^p\\
& \qquad + C \gamma^{1-p} \|\wt N_{a,q} (u|\Psi^k) \|_p^{p/2} \|S_{a,q} (u|\Psi^{k-2})\|_p^{p/2}.
\end{split}\]
The last term in the above inequality can be bounded by 
$$\frac14 \|\wt N_{a,q} (u|\Psi^k) \|_p^{p} + C \|S_{a,q} (u|\Psi^{k-2})\|_p^{p}.$$
We choose $\gamma$ such that $C\gamma^q = \frac14$. So all the term $\|\wt N_{a,q} (u|\Psi^k) \|_p^{p}$ in the right-hand side can be hidden in the left-hand side. The estimate \eqref{NbySLpZ} and then the lemma follows.
\ep

Combined with Moser's estimate and Lemma \ref{lemS<Nd}, we can improve Lemma \ref{lem6.8}.
\begin{lem}[self-improvement] \label{lem6.9}
Let $\cL$ be an elliptic operator satisfying ($H^1_\kappa$) for some constant $\kappa\geq 0$. Let $a>0$, $l>0$, and $q\in (q_0,q'_0)$ where $q_0$ is given by Proposition \ref{prop1.1}. Choose $k>12$, a positive $a^{-1}$-Lipschitz function $e$ and a ball $B\subset \R^d$ of radius $l'\geq l$. Then for any weak solution $u\in W^{1,2}_{loc}(\Omega)$,
\begin{equation} \label{NbySLp}
\|\wt N_{a,q}(u|\Psi_e^k\Psi_{B,l}^k)\|_{q}^q \lesssim \|S_{a,q}(u|\Psi_e^{k-12}\Psi_{B,l}^{k-12})\|_{q}^q + \int_{y\in \R^d} \int_{s\in \R^{n-d}} |u|^q \Psi_{B,l}^{k-3} \dr_r [\Psi_e^{k-3}] \frac{ds}{|s|^{n-d-1}}\, dy,\end{equation}
where the constant depends on $a$, $q$, $n$, $\lambda_q(\cA)$, $\|\cA\|_\infty$, $\|b^{-1}\|_\infty + \|b\|_\infty$, $\kappa$ and $k$.

\end{lem}

\begin{rmk}
	The last term is bounded by the $L^q$ average on a local region compactly contained in $\Omega$, so is finite. Thus this integral has no singularity at the boundary $\R^d \times\{0\}$, and in fact it is the same as
	\[ \iint_\Omega |u|^q \Psi_{B,l}^{k-3} \dr_r [\Psi_e^{k-3}] \frac{ds}{|s|^{n-d-1}} \, dy. \]
\end{rmk}

\bp As before, we write $\Psi$ for $\Psi_e\Psi_{B,l}$. Lemma \ref{lemMoserPsi2} entails that there exists $\epsilon$ sufficiently small (depending on $n, k, q, q_0$ and without loss of generality, we can assume $q-2\epsilon >q_0$) such that
$$\|\wt N_{a,q}(u|\Psi^k)\|_{q}^q \lesssim \|\wt N_{a,q-\epsilon}(u|\Psi^{k-3})\|_{q}^q.$$
Now, since $q-\epsilon <q$, we can use Lemma \ref{lem6.8} to get
\begin{align}
\|\wt N_{a,q}(u|\Psi^k)\|_{q}^q & \lesssim \|\wt N_{a,q-\epsilon}(u|\Psi^{k-3})\|_{q}^q \nn \\
& \lesssim \|S_{a,q-\epsilon}(u|\Psi^{k-5})\|_{q}^q  + \left\|\left(\fint_{y\in B_{ae(.)/2}(.)} \int_{s\in \R^{n-d}} |u|^{q-\epsilon} \Psi_{B,l}^{k-3} \dr_r [\Psi_e^{k-3}] \frac{ds}{|s|^{n-d-1}}\, dy \right)^\frac1{q-\epsilon} \right\|_q^q \nn \\
& := T_1 + T_2. \label{lem6.9a}
\end{align}
We first estimate $T_2$. By H\"older's inequality, for any $x\in \R^d$
\begin{equation} \label{NbyS:tmp2}
\begin{split}
& \left(\fint_{y\in B_{ae(x)/2}(x)} \int_{s\in \R^{n-d}} |u|^{q-\epsilon} \Psi_{B,l}^{k-3} \dr_r [\Psi_e^{k-3}] \frac{ds}{|s|^{n-d-1}}\, dy \right)^\frac1{q-\epsilon} \\
 & \hspace{3cm} \leq \left(\fint_{y\in B_{ae(x)/2}(x)} \int_{s\in \R^{n-d}} |u|^{q} \Psi_{B,l}^{k-3} \dr_r [\Psi_e^{k-3}] \frac{ds}{|s|^{n-d-1}}\, dy \right)^\frac1{q} \\
 & \hspace{5cm}  \left(\fint_{y\in B_{ae(x)/2}(x)} \int_{s\in \R^{n-d}} \Psi_{B,l}^{k-3} \dr_r [\Psi_e^{k-3}] \frac{ds}{|s|^{n-d-1}}\, dy \right)^\frac{\epsilon}{q(q-\epsilon)} 
\end{split}\end{equation}
Note that $0 \leq \Psi_{B,l}^{k-3} \dr_r [\Psi_e^{k-3}] \leq (k-3) \dr_r \Psi_e$, $\dr_r \Psi_e \leq 2e(x)/|s|^2 $ and it is non-zero only if $e(y)/2 \leq |s| \leq e(y)$. For any $y\in B_{ae(x)/2}(x)$, since $|e(y)-e(x)| \leq |y-x|/a < e(x)/2$, it follows that 
\[ \frac{e(x)}{2} < e(y)< \frac{3e(x)}{2}. \] 
Hence
\begin{align}
	& \fint_{y\in B_{ae(x)/2}(x)} \int_{s\in \R^{n-d}} \Psi_{B,l}^{k-3} \dr_r [\Psi_e^{k-3}] \frac{ds}{|s|^{n-d-1}}\, dy \nn \\
	& \qquad \lesssim \fint_{y\in B_{ae(x)/2}(x)} \int_{e(x)/4 \leq |s| \leq 3e(x)/2} \frac{e(x)}{|s|^2} \frac{ds}{|s|^{n-d-1}}\, dy \lesssim 1.\label{NbyS:tmp1}
\end{align}
Combining \eqref{NbyS:tmp1} with \eqref{NbyS:tmp2}, we deduce that 
\begin{equation} \label{lem6.9b}\begin{split}
T_2 & \lesssim \int_{x\in \R^d} \fint_{y\in B_{ae(x)/2}(x)} \int_{s\in \R^{n-d}} |u|^{q} \Psi_{B,l}^{k-3} \dr_r [\Psi_e^{k-3}] \frac{ds}{|s|^{n-d-1}}\, dy \, dx \\
& \lesssim \int_{y\in \R^d} \int_{s\in \R^{n-d}} |u|^{q} \Psi_{B,l}^{k-3} \dr_r [\Psi_e^{k-3}] \frac{ds}{|s|^{n-d-1}}\, dy,
\end{split}\end{equation}
where, for the last line, we use the same argument as the one used to go from \eqref{5.18} to \eqref{5.20}.

\ms

It remains to bound $T_1$. For any $x\in \R^d$,
\[\begin{split}
|S_{a,q-\epsilon}  (u|\Psi^{k-5})(x)|^q  & = \iint_{(y,s) \in \wh \Gamma_a(x)} |\nabla u|^{2} |u|^{q-2-\epsilon} \Psi^{k-5} \frac{ds}{|s|^{n-2}} \, dy \\
& \leq \left(\iint_{(y,s) \in \wh \Gamma_a(x)} |\nabla u|^{2} |u|^{q-2-2\epsilon} \Psi^{k+2} \frac{ds}{|s|^{n-2}} \, dy\right)^\frac12 \\
& \hspace{5cm}  \left(\iint_{(y,s) \in \wh \Gamma_a(x)} |\nabla u|^{2} |u|^{q-2} \Psi^{k-12} \frac{ds}{|s|^{n-2}} \, dy\right)^\frac12 \\
& \leq |S_{a,q-2\epsilon}  (u|\Psi^{k+2})(x)|^{q/2} |S_{a,q}  (u|\Psi^{k-12})(x)|^{q/2} .
\end{split}\]
So, using Cauchy-Schwarz's inequality, we can bound $T_1$ as follows
\[\begin{split}
T_1 & \leq \|S_{a,q-2\epsilon}  (u|\Psi^{k+2})\|_q^{q/2} \|S_{a,q}  (u|\Psi^{k-12})\|_q^{q/2}.
\end{split}\]
The use of Lemma \ref{lemS<Nd} and then H\"older's inequality gives that
\begin{equation} \label{lem6.9c}\begin{split}
T_1 & \lesssim \|\wt N_{a,q-2\epsilon}  (u|\Psi^{k})\|_q^{q/2} \|S_{a,q}  (u|\Psi^{k-12})\|_q^{q/2} \\
& \lesssim \|\wt N_{a,q}  (u|\Psi^{k})\|_q^{q/2} \|S_{a,q}  (u|\Psi^{k-12})\|_q^{q/2}. \\
\end{split}\end{equation}
The combination of \eqref{lem6.9a}, \eqref{lem6.9b}, and \eqref{lem6.9c} proves that
\[\begin{split}
\|\wt N_{a,q}(u|\Psi^k)\|_{q}^q \leq & C \|\wt N_{a,q}  (u|\Psi^{k})\|_q^{q/2} \|S_{a,q}  (u|\Psi^{k-12})\|_q^{q/2} \\
& + C \int_{y\in \R^d} \int_{s\in \R^{n-d}} |u|^{q} \Psi_{B,l}^{k-3} \dr_r [\Psi_e^{k-3}] \frac{ds}{|s|^{n-d-1}}\, dy,
\end{split}\]
which can be easily improved into \eqref{NbySLp}. The lemma follows.
\ep

\begin{lem} \label{lem6.10}
Let $\cL$ be an elliptic operator satisfying ($H^1_\kappa$) for some constant $\kappa\geq 0$. Let $a, l>0$, $q\in (q_0,q'_0)$ where $q_0$ is given by Proposition \ref{prop1.1}. Choose $k>12$, a positive $a^{-1}$-Lipschitz function $e$ and a ball $B:=B_{l'}(x_B) \subset \R^d$ of radius $l'\geq l$. Then for any weak solution $u\in W^{1,2}_{loc}(\Omega)$,
\begin{equation} \label{NbySLp2}
\|\wt N_{a,q}(u|\Psi_e^k\Psi_{B,l}^k)\|_{q}^q \lesssim \|S_{a,q}(u|\Psi_e^{k-12}\Psi_{B,l}^{k-12})\|_{q}^q + \int_{y\in \R^d} \int_{s\in \R^{n-d}} |u|^q \Psi_e^{k-3}  \dr_r[-\Psi_{B,l}^{k-3}] \frac{ds}{|s|^{n-d-1}}\, dy,\end{equation}
where the constant depends on $a$, $q$, $n$, $\lambda_q(\cA)$, $\|\cA\|_\infty$, $\|b^{-1}\|_\infty + \|b\|_\infty$, $\kappa$ and $k$.
\end{lem}

\bp Let $\Psi$ be the product $\Psi_e\Psi_{B,l}$. By Lemma \ref{lem6.9},
$$\|\wt N_{a,q}(u|\Psi^k)\|_{q}^q \lesssim \|S_{a,q}(u|\Psi^{k-12})\|_{q}^q + \int_{y\in \R^d} \int_{s\in \R^{n-d}} |u|^q \Psi_{B,l}^{k-3} \dr_r [\Psi_e^{k-3}] \frac{ds}{|s|^{n-d-1}}\, dy.$$
It follows by the product rule that
\begin{equation} \label{lem6.10a} \begin{split} 
\|\wt N_{a,q}(u|\Psi^k)\|_{q}^q \lesssim \|S_{a,q}(u|\Psi^{k-12})\|_{q}^q + &  \int_{y\in \R^d} \int_{s\in \R^{n-d}} |u|^q \Psi_e^{k-3} \dr_r [-\Psi_{B,l}^{k-3}]\frac{ds}{|s|^{n-d-1}}\, dy \\
& + \int_{y\in \R^d} \int_{s\in \R^{n-d}} |u|^q  \dr_r [\Psi^{k-3}]\frac{ds}{|s|^{n-d-1}}\, dy.
\end{split}\end{equation}
So it remains to bound the last term in the right-hand side of \eqref{lem6.10a}. Simply by using H\"older inequality with different powers, the proof of Lemma \ref{lem6.4} can be easily adapted to obtain
\begin{equation} \label{lem6.10b}\begin{split}
\left|\iint_\Omega |u|^q \dr_r [\Psi^{k-3}] \, \frac{dt}{|t|^{n-d-1}} \, dx  \right| & \lesssim \|S_{a,q}(u|\Psi^{k-3})\|_{q}^q + \|S_{a,q}(u|\Psi^{k-8})\|_q^{q/2} \|\wt N_{a,q}(u|\Psi^{k})\|_q^{q/2} \\
& \leq C_\eta \|S_{a,q}(u|\Psi^{k-8})\|_q^q + \eta \|\wt N_{a,q}(u|\Psi^{k})\|_q^q \\
& \leq C_\eta \|S_{a,q}(u|\Psi^{k-12})\|_q^q + \eta \|\wt N_{a,q}(u|\Psi^{k})\|_q^q
\end{split}\end{equation}
for all $\eta>0$. By choosing $\eta$ small enough, the combination of \eqref{lem6.10a} and \eqref{lem6.10b} gives \eqref{NbySLp2}. The lemma follows. \ep

\section{From local estimates to global ones and \newline existence of solutions to Dirichlet problem}

\label{SExistence}

By Lemma \ref{lemDirichlet}, for any $g\in C^\infty_0(\R^d) \subset H$, there is a (unique) energy solution $u\in W$ to $\cL u = 0$ such that $\Tr u = g$. The idea of this section is to first prove that if $\cL$ satisfies ($H^1_\kappa$) with $\kappa$ sufficiently small, then any energy solution satisfies
\begin{equation}\label{eq:Nbdbytrace}
	\|\wt N_{a,q}(u)\|_q \leq C \|\Tr u\|_q,
\end{equation}
with a universal constant $C>0$. 
Then, for any $g\in L^q(\R^d)$, the existence of a weak solution $u$ to $\cL u = 0$, whose non-tangential limit on $\R^d$ is given by $g$, can be obtained by a density argument; moreover we can show that this solution $u$ satisfies $\|\wt N_{a,q}(u)\|_q \lesssim \|g\|_q$.

\smallskip

The local inequalities proven in Lemmas \ref{lemS<Na} and \ref{lem6.10} increase the integral regions from left to right. Actually, we use some cut-off functions and we loose power on these cut-off functions, but the idea is the same: by combining Lemma \ref{lem6.10} and Lemma \ref{lemS<Na}, we don't loop back to the same local non-tangential function, so even if $\kappa$ is small, we cannot hide the term on the right-hand side to the left-hand side at the local level, where we are sure that everything is finite.  The idea to prove \eqref{eq:Nbdbytrace} is then to pass the local estimates to infinity.
The main obstacle however is that we don't know {\em a priori} that the energy solutions satisfy that $\|\wt N_{a,q}(u)\|_q$, or even $\|S_{a,q}(u)\|_q$, is finite.

\bigskip

For the sequel we use the following notations for cut-off functions. Choose the same function $\phi \in C^\infty_0([0,\infty))$ such that $0\leq \phi \leq 1$, $\phi \equiv 1$ on $[0,1]$, $\phi \equiv 0$ outside $[0,2]$, $\phi$ decreasing, and $|\phi'| \leq 2$. 
For $\epsilon>0$, we define $\Psi_\epsilon$ as
$$\Psi_\epsilon(x,t)  = \Psi_\epsilon(t) = \phi\left( \frac{\epsilon}{|t|} \right).$$
For $l>0$, we define $\chi_l$ as 
$$\chi_l(x,t) = \chi_l(t) = \phi\left( \frac{a|t|}{l} \right), $$
and if $B \subset \R^d$ is a ball, we define
$$ \Phi_{B,l}(x,t) = \Phi_{B,l}(x) = \phi\left(1+ \frac{\dist(x,B)}{l}\right). $$
The reader may recall Lemma \ref{exgoodcutoff} and recognize that $\Psi_{\epsilon}$ is the function $\Psi_{e}$ there with $e(x) \equiv \epsilon$, and the product $\chi_l \Phi_{B,l}$ is the function $\Psi_{B,l}$ there. These cut-off functions correspond to smooth cut-off away from the boundary, at infinity in the $t$ and $x$ variables, respectively. Also recall that, when the radius of $B$ is bigger than $l$, the function $\Psi_\epsilon \chi_l \Phi_{B,l}$ satisfies $(\cH^2_{a,M_0})$ for some $M_0$ depending only on $a$ and the dimensions $d,n$ (see Lemma \ref{exgoodcutoff}). 

\medskip

\begin{lem}[$N< S + \Tr$ when $q\geq 2$] \label{lem7.1}
Let $\cL$ be an elliptic operator satisfying $(\cH^1_\kappa)$ for some constant $\kappa\geq 0$. Let $a, l>0$, $2\leq q < q'_0$ where $q'_0$ is the conjugate of $q_0$ given by Proposition \ref{prop1.1}. Choose $k>12$. Then for any energy solution $u\in W$ to $\cL u = 0$,
\begin{equation}\label{eq:lem7.1}
\|\wt N_{a,q}(u|\chi_l^k)\|_{q}^q \lesssim \|S_{a,q}(u|\chi_l^{k-12})\|_{q}^q + \|\Tr u\|_{q}^q < +\infty,
\end{equation}
where the constant of the first inequality depends on $a$, $q$, $n$, $d$, $\lambda_q(\cA)$, $\|\cA\|_\infty$, $\|b^{-1}\|_\infty + \|b\|_\infty$, $\kappa$ and $k$.
\end{lem} 

\begin{rmk} \label{rem7.2}
An immediate consequence of the lemma is that if $u\in W$ is an energy solution to $\mathcal L u = 0$, then $\|\wt N_{a,q}(u|\chi_l^k)\|_{q}$ is finite for any $k>0$ (not only for $k>12$). Indeed, since $\chi_l^k \leq \chi_{2l}^{20}$ for all $k>0, l>0$, we have
$$\|\wt N_{a,q}(u|\chi_l^k)\|_{q} \leq \|\wt N_{a,q}(u|\chi_{2l}^{20})\|_{q}< +\infty.$$
The same remark holds for Lemma \ref{lem7.1ter}
\end{rmk}

\bp We start by proving finiteness. First, if $u$ is an energy solution and $q\geq 2$, we have that $\int_\Omega |\nabla u|^2 |u|^{q-2} |t|^{d+1-n} dt \, dx < +\infty$ by Theorem \ref{theoES} (i), and also that $\Tr u \in C_0^\infty(\R^d) \subset L^q(\R^d)$. 
Then, since $\chi_l \leq 1$ is supported in $\{(x,t)\in \Omega, \, |t| \leq 2l/a\}$, we have
\begin{align}
\|S_{a,q}(u|\chi_l^{k-12})\|_{q}^q & \simeq \iint_{(x,t) \in \Omega} |\nabla u|^2 |u|^{q-2} \chi_l^{k-12} \frac{dt}{|t|^{n-d-2}}\, dx \nn \\
& \quad \leq  \frac{2l}{a} \iint_{(x,t) \in \Omega} |\nabla u|^2 |u|^{q-2} \frac{dt}{|t|^{n-d-1}}\, dx < +\infty. \label{lem7.1a}
\end{align}
So we indeed have
\[ \|S_{a,q}(u|\chi_l^{k-12})\|_{q}^q + \|\Tr u\|_{q}^q < +\infty.
\]

\medskip

The proof of the first inequality in \eqref{eq:lem7.1}, in simple words, is by passing to the limit the estimate in Lemma \ref{lem6.9}.

\noindent {\bf Step 1:} 
We claim that for any ball $B$ with radius $l' \geq l$ and for any $k>12$,
\begin{equation} \label{lem7.1b}
\|\wt N_{a,q}(u|\chi_l^{k}\Phi_{B,l}^{k})\|_{q}^q < +\infty.
\end{equation}

\smallskip

We take $\epsilon>0$ sufficiently small ($\epsilon < l/a$). We write $\Psi$ for $\chi_l \Phi_{B,l} \Psi_\epsilon$. According to Lemma \ref{lem6.10},
$$\|\wt N_{a,q}(u|\Psi^{k})\|_{q}^q \lesssim \|S_{a,q}(u|\Psi^{k-12})\|_{q}^q + \int_{y\in \R^d} \int_{s\in \R^{n-d}} |u|^q \Psi_\epsilon^{k-3} \Phi_{B,l}^{k-3} \dr_r[-\chi_{l}^{k-3}] \frac{ds}{|s|^{n-d-1}}\, dy.$$
Note that $\Phi_{B,l}$ is independent of $t$ so we can pull it out of the $r$-derivative.
Since $\Psi \leq \chi_l$, $\Psi_\epsilon \equiv 1$ on $\supp \,  \dr_r \chi_l \subset \{l/a \leq |t| \leq 2l/a\} $, and $|\dr_r \chi_l (s)| \lesssim \frac{1}{|s|}$, we deduce
$$\|\wt N_{a,q}(u|\Psi^k)\|_{q}^q \lesssim \|S_{a,q}(u|\chi_l^{k-12})\|_{q}^q + \iint_{(y,s)\in \supp\, \Phi_{B,l} \dr_r \chi_l} |u|^q \frac{ds}{|s|^{n-d}}\, dy.$$
Notice that the right hand side is independent of $\epsilon$, and is finite. Indeed, $\|S_{a,q}(u|\chi_l^{k-12})\|_{q}^q < +\infty$ by \eqref{lem7.1a}; moreover, $\Phi_{B,l} \dr_r \chi_l$ is compactly supported in $\Omega$ and by Lemma \ref{lemMoser} $u\in L^q_{loc}(\Omega)$. We deduce that 
$$\|\wt N_{a,q}(u|\Psi^{k})\|_{q}^q = \|\wt N_{a,q}(u|\chi_l^{k}\Phi_{B,l}^{k} \Psi_\epsilon^{k})\|_{q}^q$$
is uniformly bounded in $\epsilon$, and so the claim \eqref{lem7.1b} follows by passing $\epsilon \to 0$.

\ms

\noindent {\bf Step 2:} We want to prove that for any $k\geq 1$ and any ball $B$ with radius $l' \geq l$, we have
\begin{equation} \label{lem7.1c}
\limsup_{\epsilon \to 0} \int_{y\in \R^d} \int_{s\in \R^{n-d}} |u|^q \chi_l^{k}\Phi_{B,l}^{k} \dr_r [\Psi_\epsilon^{k}] \frac{ds}{|s|^{n-d-1}}\, dy \lesssim \int_{x\in 4B} |\Tr u|^q \, dx.
\end{equation}

\smallskip

Note that $\supp \, \dr_r \Psi_\epsilon \subset \{(x,t) \in \Omega, \, \epsilon/2 \leq |t| \leq \epsilon\}$, $0\leq \dr_r \Psi_\epsilon \lesssim \frac1{\epsilon}$. 
Hence
\begin{align}
	\int_{y\in \R^d} \int_{s\in \R^{n-d}} |u|^q \chi_l^{k}\Phi_{B,l}^{k} \dr_r [\Psi_\epsilon^{k}] \frac{ds}{|s|^{n-d-1}}\, dy 
& \lesssim \int_{y\in 2B} \fint_{\epsilon/2 \leq |s| \leq \epsilon} |u|^q \chi_l^k \Phi_{B,l}^k  ds\, dy \nonumber \\
& \lesssim \int_{x\in 3B} \fint_{y\in B_{a\epsilon/2}(x)} \fint_{\epsilon/2 \leq |s| \leq \epsilon} |u|^q \chi_l^k \Phi_{B,l}^k ds\, dy\, dx. \label{ex:tmp1}
\end{align} 
For the last inequality we use Fubini's lemma. By \eqref{lem7.1b} of Step 1,
\begin{equation}\label{ex:RFasmp}
	\fint_{y\in B_{a\epsilon/2}(x)} \fint_{\epsilon/2 \leq |s| \leq \epsilon} |u|^q \chi_l^k \Phi_{B,l}^k ds\, dy = \left|(u|\chi_l^k \Phi_{B,l}^k)_{W,a,q}(x,\epsilon) \right|^k \leq |\wt N_{a,q}(u|\chi_l^k \Phi_{B,l}^k)(x)|^q  
\end{equation} 
is integrable uniformly in $\epsilon$.
When $q\geq 2$, Moser's estimate (Lemma \ref{lemMoser} (ii)) gives that
\begin{equation}
	\left(\fint_{y\in B_{a\epsilon}(x)} \fint_{\epsilon/2 \leq |s| \leq \epsilon} |u|^q  ds\, dy\right)^\frac1q \lesssim \left(\fint_{y\in B_{2a\epsilon}(x)} \fint_{\epsilon/4 \leq |s| \leq 2\epsilon} |u|^2 ds\, dy\right)^\frac12. \label{ex:tmp2}
\end{equation} 
(We remark that when $q<2$, the above estimate also holds by H\"older's inequality.)
We claim that
\begin{equation} \label{lem7.1e}
\limsup_{\epsilon\to 0} \left(\fint_{y\in B_{2a\epsilon}(x)} \fint_{\epsilon/4 \leq |s| \leq 2\epsilon} |u|^2  ds\, dy\right)^{\frac12} \lesssim |\Tr u(x)|,
\end{equation}
for $\sigma$-almost every $x\in\Gamma$. Then by reverse Fatou's lemma and the pointwise domination \eqref{ex:RFasmp}, we get
\[ \limsup_{\epsilon \to 0} \int_{x\in 4B} \left(\fint_{y\in B_{2a\epsilon}(x)} \fint_{\epsilon/4 \leq |s| \leq 2\epsilon} |u|^2  ds\, dy\right)^{\frac{q}{2}} dx \lesssim \int_{x\in 4B} |\Tr u|^q dx. \]
This estimate, combined with \eqref{ex:tmp1} and \eqref{ex:tmp2}, proves \eqref{lem7.1c}.

By triangle inequality 
\begin{align}
	\fint_{y\in B_{2a\epsilon}(x)} \fint_{\epsilon/4 \leq |s| \leq 2\epsilon} |u|^2  ds\, dy & \lesssim \fint_{y\in B_{2a\epsilon}(x)} \fint_{\epsilon/4 \leq |s| \leq 2\epsilon} |u - \Tr u(x) |^2  ds\, dy + |\Tr u(x)|^2 \nonumber \\
	& \lesssim \fint_{y\in B_{2a\epsilon}(x)} \fint_{|s| \leq 2\epsilon} |u - \Tr u(x) |^2  ds\, dy + |\Tr u(x)|^2.\label{ex:trg}
\end{align}
Since $u\in W$, by the Lebesgue density property in Theorem 3.13 of \cite{DFMprelim}\footnote{The Lebesgue density property in \cite{DFMprelim} is a $L^1$ version, while we need a stronger $L^2$ version. We pretend that the argument of \cite{DFMprelim} can be modified without much difficulty to get the desired $L^2$ Lebesgue density property, and a similar $L^2$ Lebesgue density property will be properly written by the authors in an article under preparation.}, we know that the first term in the right-hand side above tends to 0 as $\epsilon \to 0$ for $\sigma$-almost every $x\in \Gamma$. The claim \eqref{lem7.1e} follows.

\ms

\noindent {\bf Step 3:} Conclusion.

\smallskip

Let $B=B_{l'}$ be a ball with center $0$ and radius $l' \geq l$ and $\epsilon >0$. Apply Lemma \ref{lem6.9} to the function $\Psi = \chi_l \Phi_{B,l} \Psi_\epsilon$, we get 
$$\|\wt N_{a,q}(u|\Psi^k)\|_{q}^q \lesssim \|S_{a,q}(u|\Psi^{k-12})\|_{q}^q + \int_{y\in \R^d} \int_{s\in \R^{n-d}} |u|^q \chi_l^{k-3} \Phi_{B,l}^{k-3} \dr_r [\Psi_\epsilon^{k-3}] \frac{ds}{|s|^{n-d-1}}\, dy$$
By Step 2, taking the limit as $\epsilon$ goes to 0 gives that
$$\|\wt N_{a,q}(u|\chi_l^k\Phi_{B,l}^k)\|_{q}^q \lesssim \|S_{a,q}(u|\chi_l^{k-12}\Phi_{B,l}^{k-12})\|_{q}^q + \int_{x\in 4B} |\Tr u|^q \, dx.$$
We take now the limit as the radius $l'$ goes to $+\infty$ and we obtain \eqref{eq:lem7.1}. \ep

\begin{rmk} \label{rmklem1}
We remark that the above proof does not use $q\geq 2$ per se: We only need this assumption to guarantee that $\|S_{a,q}(u|\chi_l^{k-12})\|_q$ is finite (see Step 1). That is to say, we can use the same argument for the case $q<2$, if we know a priori the corresponding square function is integrable.
\end{rmk}

\begin{lem}[$N< S+\Tr$ when $q<2$] \label{lem7.1ter}
Let $\cL$ be an elliptic operator satisfying $(\cH^1_\kappa)$ for some constant $\kappa\geq 0$. Let $a, l>0$, $q\in (q_0,q'_0)$ where $q_0$ is given by Proposition \ref{prop1.1}. Choose $k>12$. Then for any energy solution $u\in W$ to $\cL u = 0$,
\begin{equation}\label{eq:lem7.1ter}
	\|\wt N_{a,q}(u|\chi_l^k)\|_{q}^q \lesssim \|S_{a,q}(u|\chi_l^{k-12})\|_{q}^q + \|\Tr u\|_{q}^q < +\infty,
\end{equation}
where the constant of the first inequality depends on $a$, $q$, $n$, $d$, $\lambda_q(\cA)$, $\|\cA\|_\infty$, $\|b^{-1}\|_\infty + \|b\|_\infty$, $\kappa$ and $k$.
\end{lem}

\bp
	We start by proving a priori finiteness.
	\medskip
	
	\noindent \textbf{Step 1:} We claim that if $l$ is sufficiently large (depending only on $a$ and the support of $\Tr u$) and $B\subset \R^d$ is any ball centered at zero with radius greater than $l$, then
		\begin{equation}\label{eq:ql2tail}
			\|\wt N_{a,q}(u|\chi_l^k(1-\Phi_{B,l})^k)\|_{q}^q \lesssim \|S_{a,q}(u|\chi_l^{k-12}(1-\Phi_{B,l})^{k-12})\|_{q}^q + \|\Tr u\|_{q}^q < +\infty.
		\end{equation}
		
		We choose $l$ so that the ball in $\R^n$ centered at zero with radius $l/a$ contains two times the ball $B_0 \subset \R^n$ given by Theorem \ref{theoES} (ii). (Recall that $B_0$ depends on the support of $\Tr u$.)
	 Recall that
\begin{align}
\|S_{a,q}(u|\chi_l^{k-12}(1-\Phi_{B,l})^{k-12})\|_{q}^q & \simeq \iint_{(x,t) \in \Omega} |\nabla u|^2 |u|^{q-2} \chi_l^{k-12}(1-\Phi_{B,l})^{k-12} \frac{dt}{|t|^{n-d-2}}\, dx \nn \\
 & \lesssim l \iint_{(x,t) \in \Omega} |\nabla u|^2 |u|^{q-2} (1-\Phi_{B,l})^{k-12} \frac{dt}{|t|^{n-d-1}}\, dx \nn
\end{align}
because $|t| \leq 2l$ in the support of $\chi_l$. In addition, the support of $(1-\Phi_{B,l})^{k-12}$ is contained in the complement of a cylinder of radius $\sim l/a$, thus by the choice of $l$ it is contained in $\R^n \setminus B_0$. 
Lemma \ref{lemES1} then allows us to conclude that the right-hand side above - and so the left-hand side - is finite. 
Again by Remark \ref{rmklem1}, once we know $\|S_{a,q}(u|\chi_l^{k-12}(1-\Phi_{B,l})^{k-2})\|_{q}<+\infty$, \eqref{eq:ql2tail} follows by the same argument as in the proof of Lemma \ref{lem7.1}.
\medskip

\noindent \textbf{Step 2:} We claim that for any $k>2$,
\begin{equation} \label{lem7.1''a}
\|S_{a,q}(u|\chi_l^{k})\|_{q} \lesssim \|\wt N_{a,q}(u|\chi_l^{k-2})\|_{q} < +\infty,\end{equation}

On one hand we observe that for any $l>0$ and any ball $B'$ of radius at least $l$, one has 
\begin{equation} \label{Naqfini2}
\|\wt N_{a,q}(u|\chi_l^{k-2}\Phi_{B',l}^{k-2})\|_{q} < +\infty.
\end{equation}
Indeed, the above finiteness is an immediate consequence of H\"older's inequality, since 
\[ \wt N_{a,q}(u|\chi_l^{k-2}\Phi_{B',l}^{k-2}) \leq \wt N_{a,2}(u|\chi_l^{k-2}), \] 
$\wt N_{a,q}(u|\chi_l^{k-2}\Phi_{B',l}^{k-2})$ has compact support, and $\| \wt N_{a,2}(u|\chi_l^{k-2})\|_2 <+\infty$ by Lemma \ref{lem7.1}.
On the other hand, in Step 1 we show that $\|\wt N_{a,q}(u|\chi_l^{k-2}(1-\Phi_{B,l})^{k-2})\|_{q}< +\infty$ if $l$ is sufficiently large. (A priori we only have finiteness when $k-2>12$, but since $0\leq \chi_l\leq 1$, the finitness clearly holds for any $k>2$.) A simple computation shows that $1\leq \Phi_{2B,l}^{k-2} + (1-\Phi_{B,l})^{k-2}$, and hence
\[ \| \wt N_{a,q}(u|\chi_l^{k-2})\|_{q}  \lesssim  \|\wt N_{a,q}(u|\chi_{l}^{k-2}\Phi_{2B,l}^{k-2})\|_{q} + \|\wt N_{a,q}(u|\chi_l^{k-2}(1-\Phi_{B,l})^{k-2})\|_{q} < +\infty. \]
Now take an increasing sequence of balls $(B_i)_{i\geq 1}$ such that $B_1$ is of radius $l$, $B_i \subset B_{i+1}$ and $\bigcup B_i = \R^d$. We apply Lemma \ref{lemS<Nb} to the functions $\Psi = \chi_l \Phi_{B_i,l} \Psi_{1/i}$, which gives that
$$\|S_{a,q}(u|\Psi^{k})\|_{q} \lesssim \|\wt N_{a,q}(u|\Psi^{k-2})\|_{q}.$$
By the finiteness of $\|\wt N_{a,q}(u|\chi_l^{k-2})\|_{q}$, we may take the limit as $i \to +\infty$ and obtain \eqref{lem7.1''a}.

\medskip

\noindent \textbf{Step 3:} Conclusion.
Having shown the finiteness of $\|S_{a,q}(u|\chi_l^{k-12})\|_q$, by Remark \ref{rmklem1}, the same argument in the proof of Lemma \ref{lem7.1} gives \eqref{eq:lem7.1ter}.
\ep

The following estimate is a byproduct of Lemma \ref{lem7.1}, in particular, of \eqref{lem7.1c}. It will be used later in the proof of Lemma \ref{lemS<Ne}. 
\begin{cor}\label{lem7.1'}
	Let $\cL$ be an elliptic operator satisfying $(\cH^1_\kappa)$ for some constant $\kappa\geq 0$. Let $a>0$, $q\in (q_0,q'_0)$ where $q_0$ is given by Proposition \ref{prop1.1}. Choose $k\geq 1$. Then for any energy solution $u\in W$ to $\cL u = 0$,
\begin{equation} \label{lem7.1'a}
\limsup_{\epsilon \to 0} \iint_{\Omega} |u|^q \dr_r [\Psi_\epsilon^{k}] \frac{ds}{|s|^{n-d-1}}\, dy \lesssim \|\Tr u\|_q^q.
\end{equation}
where the constant depends on $a$, $q$, $n$, $d$, $\lambda_q(\cA)$, $\|\cA\|_\infty$, $\|b^{-1}\|_\infty + \|b\|_\infty$, $\kappa$ and $k$.
\end{cor} 

\bp Let $l>0$ be fixed. Since $\left| \wt N_{a,q}(u|\chi_l^k) \right|^q$ is integrable, the same argument in Step 2 of Lemma \ref{lem7.1} allows us to conclude
\[ \limsup_{\epsilon \to 0} \int_{y\in \R^d} \int_{s\in \R^{n-d}} |u|^q \chi_l^{k} \dr_r [\Psi_\epsilon^{k}] \frac{ds}{|s|^{n-d-1}}\, dy \lesssim \int_{x\in \R^d } |\Tr u|^q \, dx.
\]
Since $\chi_l \equiv 1$ on $\supp \, \dr[\Psi^{k}_\epsilon]$ if $\epsilon$ is small enough ($\epsilon<l/a$), it follows that
\begin{align*}
	\limsup_{\epsilon \to 0} \int_{y\in \R^d} \int_{s\in \R^{n-d}} |u|^q \dr_r [\Psi_\epsilon^{k}] \frac{ds}{|s|^{n-d-1}}\, dy & = \limsup_{\epsilon \to 0} \int_{y\in \R^d} \int_{s\in \R^{n-d}} |u|^q \chi_l^{k} \dr_r [\Psi_\epsilon^{k}] \frac{ds}{|s|^{n-d-1}}\, dy \\
	& \lesssim \int_{x\in \R^d } |\Tr u|^q \, dx.
\end{align*}
\ep

The next lemma can be seen as the analogue of Lemma \ref{lemS<Na} for energy solutions.

\begin{lem}[$S< \kappa N + \Tr$] \label{lemS<Ne}
Let $\cL$ be an elliptic operator satisfying ($H^1_\kappa$) for some constant $\kappa\geq 0$. Let $a, l>0$, $q\in (q_0,q'_0)$ where $q_0$ is given by Proposition \ref{prop1.1}. Choose $k>2$. Then for any energy solution $u\in W$ to $\cL u = 0$,
\begin{equation} \label{lem7.2a} \begin{split}
c \|S_{a,q}(u|\chi_{l}^{k})\|_q^q
& \leq  C \kappa \|\wt N_{a,q} (u|\chi_{l}^{k-2})\|_{q}^q  + C \|\Tr u\|_q^q \\
& \hspace{3cm} + \, \frac1q \, \iint_{\Omega} \left(\frac{|u|^q}{|t|} - \dr_r [|u|^q] \right) \dr_r[\chi_l^k] \, \frac{dt}{|t|^{n-d-2}}\, dx,
\end{split} \end{equation}
where the constants $c,C>0$ depend on $a$, $q$, $n$, $\lambda_q(\cA)$, $\|\cA\|_\infty$, $\|b\|_{\infty} + \|b^{-1}\|_\infty$, $k$ and (the upper bound of) $\kappa$.
\end{lem}

\bp
We take $\epsilon>0$ and a sequence of balls $(B_i)_{i\geq 1}$ such that $B_1$ is of radius $100l$, $B_i \subset B_{i+1}$, $\bigcup B_i = \R^d$. We apply Lemma \ref{lemS<Na} with the function $\Psi_i = \chi_l \Phi_{B_i,l} \Psi_{\epsilon}$, which gives that
\begin{equation} \label{lem7.2b} \begin{split}
c \|S_{a,q}(u|\Psi_i^{k})\|_q^q
& \leq 
\frac12 \ \Re \iint_{\Omega} \cA' \nabla u \cdot \nabla [|u|^{q-2} \bar u] \frac{\Psi_i^{k}}{b}\, \frac{dt}{|t|^{n-d-2}} dx \\
& \leq  C \kappa \|\wt N_{q,a} (u|\Psi_i^{k-2})\|_{q}^q  - \,  \Re \iint_{\Omega} \cA'\nabla u \cdot \nabla_x [\Psi_i^{k}] \, \frac{|u|^{q-2} \bar u}{b} \frac{dt}{|t|^{n-d-2}}\, dx \\
& \qquad \qquad  + \, \frac1q \, \iint_{\Omega} \left(\frac{|u|^q}{|t|}\nabla |t| - \nabla_t [|u|^q] \right)\cdot \nabla_t [\Psi_i^{k}] \, \frac{dt}{|t|^{n-d-2}}\, dx,
\end{split} \end{equation}
Passing $i\to +\infty$ and then passing $\epsilon \to 0$, the left hand side converges to $c \|S_{a,q}(u|\chi_{l}^{k})\|_q^q$, which is finite by Lemmas \ref{lem7.1} and \ref{lem7.1ter}.
Clearly $ \|\wt N_{q,a} (u|\Psi_i^{k-2})\|_{q}\leq  \|\wt N_{q,a} (u|\chi_l^{k-2})\|_{q}$. So to prove \eqref{lem7.2a} it suffices to establish
\begin{equation} \label{claimu}
\lim_{i\to \infty} \Re \iint_{\Omega} \cA'\nabla u \cdot \nabla_x [\Psi_i^{k}] \, \frac{|u|^{q-2} \bar u}{b} \frac{dt}{|t|^{n-d-2}}\, dx = 0
\end{equation}
and
\begin{align} 
& \limsup_{\epsilon \to 0} \lim_{i\to \infty} \frac1q \, \iint_{\Omega}  \left(\frac{|u|^q}{|t|}\nabla |t| - \nabla_t [|u|^q] \right)\cdot \nabla_t [\Psi_i^{k}] \, \frac{dt}{|t|^{n-d-2}}\, dx  \nn \\
&  \qquad \leq C \|\Tr u\|_q^{q/2} \|S_{a,q}(u|\chi_l^{k})\|_q^{q/2} + C\|\Tr u\|_q^q + \frac1q \, \iint_{\Omega} \left(\frac{|u|^q}{|t|} - \dr_r [|u|^q] \right) \dr_r[\chi_l^k] \, \frac{dt}{|t|^{n-d-2}}\, dx. \label{claimv}
\end{align}
In the last inequality, we have the additional term $C\|\Tr u\|_q^{q/2} \|S_{a,q}(u|\chi_l^{k})\|_q^{q/2}$ which doesn't appear in the right hand side of \eqref{lem7.2a}, but this term can be bounded by $C_\eta \|\Tr u\|_q^q + \eta \|S_{a,q}(u|\chi_l^{k})\|_q^{q}$ and then we can hide $\eta \|S_{a,q}(u|\chi_l^{k})\|_q^{q}$ in the left-hand side by taking $\eta$ small enough.

\ms

Let us start with \eqref{claimu}. Since $\chi_l$ and $\Psi_{\epsilon}$ are $x$-independent, we have 
\[
\nabla_x [\Psi_i^k] = k\Psi_i^{k-1} \cdot \chi_l \Psi_\epsilon \nabla_x \Phi_{B_i,l}.
\] 
But $\nabla_x \Phi_{B_i,l}$ is supported outside the ball $B_i$, and $\Phi_{B_i/2,l} \equiv 0$ on $\supp\, \nabla_x \Phi_{B_i,l}$ if $l_i \geq 2l$. So $\nabla_x \Phi_{B_i,l} = (1-\Phi_{B_i/2,l})^k \nabla_x \Phi_{B_i,l}$ and
\begin{align}
T_1 & := \left| \Re \iint_{\Omega} \cA'\nabla u \cdot \nabla_x [\Psi_i^{k}] \, \frac{|u|^{q-2} \bar u}{b} \frac{dt}{|t|^{n-d-2}}\, dx \right| \nn \\
& \lesssim \iint_{\Omega} |\nabla u| (1-\Phi_{B_i/2,l})^k \Psi_{i}^{k-1} \chi_l \Psi_{\epsilon} |\nabla_x \Phi_{B_i,l}| |u|^{q-1} \frac{dt}{|t|^{n-d-2}}\, dx \nn \\
& \lesssim \left( \iint_{\Omega} |\nabla u|^2 |u|^{q-2} (1-\Phi_{B_i/2,l})^k \Psi_i^k \frac{dt}{|t|^{n-d-2}}\, dx \right)^\frac12 \nn \\
& \hspace{4cm}  \left( \iint_{\Omega} |u|^{q} (1-\Phi_{B_i/2,l})^k \Psi_i^{k-2} |\nabla_x [\Phi_{B_i,l}\chi_l]|^2 \frac{dt}{|t|^{n-d-2}}\, dx \right)^\frac12 \nn \\
& \lesssim \|S_{a,q}(u|\Psi_i^k(1-\Phi_{B_i/2,l})^k)\|_q^{q/2} \|\wt N_{a,q}(u|\Psi_i^{k-2}(1-\Phi_{B_i/2,l})^{k})\|_q^{q/2},\label{ex:T1}
\end{align}
by Lemma \ref{lemH2} and the fact that $\Phi_{B_i,l}\chi_l$ satisfies $(\cH^2_{a})$.
It then follows from Lemma \ref{lemS<Nb} that
\begin{equation}\label{ex:tmp4}
	T_1 \lesssim \|\wt N_{a,q}(u|\Psi_i^{k-2}(1-\Phi_{B_i/2,l})^{k-2})\|_q^{q}.
\end{equation}
By the construction of $\Psi_i$, $\wt N_{a,q}(u|\Psi_i^{k-2}(1-\Phi_{B_i/2,l})^{k-2})$ is supported in $\R^d \setminus (B_i/4)$, so the right-hand side of \eqref{ex:tmp4} is bounded by $\|\wt N_{a,q}(u|\chi_l^{k-2})\|_{L^q(\R^d \setminus B_i/4)}^q$. Lemmas \ref{lem7.1} and \ref{lem7.1ter} imply that $\wt N_{a,q}(u|\chi_l^{k-2})$ is integrable in $L^q(\R^d)$. Therefore since the balls $B_i$ increase to $\R^d$, we have 
$$\lim_{i\to \infty} T_1 = 0.$$
The claim \eqref{claimu} follows. 

We turn to the proof of \eqref{claimv}. First, since $\Psi_i$ depends only on $x, |t|$ and not on $t/|t|$, by simple algebra we have
\begin{align*}
	\left|\left(\frac{|u|^q}{|t|}\nabla |t| - \nabla_t [|u|^q] \right)\cdot \nabla_t [\Psi_i^{k}]\right| & = \left|\left(\frac{|u|^q}{|t|} - \dr_r [|u|^q] \right) \, \dr_r [\Psi_i^{k}]\right| \\
	& \leq \left|\left(\frac{|u|^q}{|t|} - \dr_r [|u|^q] \right) \, \dr_r [\chi_l^k \Psi_{\epsilon}^{k}]\right|.
\end{align*} 
That is, the integrand is bounded by a function independent of $i$. Moreover, we let the reader check that the latter is integrable:
\begin{align*}
	& \iint_{\Omega} \left|\left(\frac{|u|^q}{|t|} - \dr_r [|u|^q] \right) \, \dr_r [\chi_l^{k} \Psi_{\epsilon}^k]\right| \, \frac{dt}{|t|^{n-d-2}}\, dx \\
	& \qquad \qquad \lesssim  \|\wt N_{a,q}(u|\Psi_i^{k-1})\|_q^q + \|\wt N_{a,q}(u|\Psi_i^{k-2})\|_q^{q/2} \|S_{a,q}(u|\Psi_i^k)\|_q^{q/2} \\
& \qquad \qquad \lesssim \|\wt N_{a,q}(u|\Psi_i^{k-2})\|_q^q  \lesssim \|\wt N_{a,q}(u|\chi_l^{k-2})\|_q^q < +\infty.
\end{align*}
Therefore by Lebesgue dominated convergence theorem,
\begin{align*}
	& \lim_{i\to \infty} \frac1q \, \iint_{\Omega}  \left(\frac{|u|^q}{|t|}\nabla |t| - \nabla_t [|u|^q] \right)\cdot \nabla_t [\Psi_i^{k}] \, \frac{dt}{|t|^{n-d-2}}\, dx \\
	& \qquad \qquad = \frac1q \, \iint_{\Omega} \left(\frac{|u|^q}{|t|} - \dr_r [|u|^q] \right) \, \dr_r [\chi_l^k\Psi_\epsilon^{k}] \, \frac{dt}{|t|^{n-d-2}}\, dx < +\infty 
\end{align*}
We denote the second integral as $T_2$.
If $\epsilon$ is small ($\epsilon< l/a$), then $\Psi_{\epsilon} \equiv 1$ on $\supp\, \dr_r \chi_l$ and $\chi_l \equiv 1$ on $ \supp \, \dr_r \Psi_{\epsilon}$. Hence, 
$$\dr_r[\chi_l^k\Psi_\epsilon^{k}] = \dr_r[\Psi_{\epsilon}^k] + \dr_r[\chi_l^k],$$
which implies that
\[\begin{split}
T_2 & = \frac1q \, \iint_{\Omega} \left(\frac{|u|^q}{|t|}- \dr_r [|u|^q] \right) \dr_r [\Psi_{\epsilon}^{k}] \, \frac{dt}{|t|^{n-d-2}}\, dx \\
& \qquad + \frac1q \, \iint_{\Omega} \left(\frac{|u|^q}{|t|}- \dr_r [|u|^q] \right) \dr_r [\chi_l^{k}] \, \frac{dt}{|t|^{n-d-2}}\, dx \\
& =: T_3 + T_4.
\end{split}\]
Both $T_3$ and $T_4$ are finite when $\epsilon$ is small, because $T_3$ and $T_4$ are just the decomposition of the integral $T_2$ into two integrals on disjoint sets $\supp\, \dr_r \Psi_\epsilon$ and $\supp\, \dr_r \chi_l$ (and the function under the integral of $T_2$ is integrable). The term $T_4$ is independent of $\epsilon$, so we don't need to touch it anymore. 
Now, we consider 
$$U_{\epsilon}:= \iint_{\Omega} |u|^q\dr_r [\Psi_{\epsilon}^{k}] \, \frac{dt}{|t|^{n-d-1}}\, dx \geq 0.$$
Since $|\dr_r [|u|^q]| \lesssim |u|^{q-1} |\nabla u|$ and $\dr_r [\Psi_{\epsilon}^{k}] \geq 0$ by construction, 
by H\"older inequality we have 
\[\begin{split}
T_3 & \lesssim U_\epsilon + U^{1/2}_\epsilon \left( \iint_\Omega |\nabla u|^2 |u|^{q-2} \dr_r [\Psi_{\epsilon}^{k}] \, \frac{dt}{|t|^{n-d-3}}\, dx \right)^\frac12.
\end{split}\]
Observe now that $ \dr_r [\Psi_{\epsilon}] \lesssim 1/|t|$ and if $\epsilon$ is small, $\chi_l \equiv 1$ on $\supp\, \dr_r[\Psi_\epsilon^k]$. We deduce
\[\begin{split}
T_3 & \lesssim U_\epsilon + U^{1/2}_\epsilon \left( \iint_\Omega |\nabla u|^2 |u|^{q-2} \chi_l^k \, \frac{dt}{|t|^{n-d-2}}\, dx \right)^\frac12 \lesssim U_\epsilon + U^{1/2}_\epsilon \|S_{a,q}(u|\chi_l^k)\|_q^{q/2}.
\end{split}\]
By taking the limit as $\epsilon \to 0$, and recall that Lemma \ref{lem7.1'} shows $\limsup_{\epsilon\to 0} U_\epsilon \lesssim \|\Tr u\|_q^q$. It follows that entails that
$$\limsup_{\epsilon \to 0} T_3 \lesssim \|\Tr u\|_q^q + \|\Tr u\|_q^{q/2} \|S_{a,q}(u|\chi_l^k)\|_q^{q/2}.$$
 The claim \eqref{claimv} and then the lemma follow.
\ep

The next step is to prove a bound on the growth of the energy solution $u$. It will be used to get a global estimate in the proof of Lemma \ref{lem7.4}.

\begin{lem} \label{lem7.3}
Let $\cL$ be an elliptic operator. Let $q\in (q_0,q'_0)$ where $q_0$ is given by Proposition \ref{prop1.1}. Choose $k\geq 1$. Then for any energy solution $u\in W$ to $\cL u = 0$,
\begin{equation} \label{lem7.3a}
\lim_{l\to +\infty} \frac1l \iint_\Omega |u|^q \dr_r[-\chi_l^k] \frac{dt}{|t|^{n-d-1}}\, dx = 0.
\end{equation}
\end{lem}

\bp First, by the definition of $\chi_l$ and $\phi$, we have $0 \leq \dr_r[-\chi_l^k] \lesssim 1/l$ and 
$\supp\, \dr_r[-\chi_l^k] \subset \{(x,t)\in \Omega, \, l\leq a|t| \leq 2l\}$. Therefore,
\begin{equation} \label{lem7.3b}
\frac1l \iint_\Omega |u|^q \dr_r[-\chi_l^k] \frac{dt}{|t|^{n-d-1}}\, dx  \lesssim \frac1{l^2} \int_{x\in \R^d} \int_{l\leq a|t| \leq 2l} |u|^q \frac{dt}{|t|^{n-d-1}}\, dx.
\end{equation}
Let $(B_i)_{i\in I}$ be balls of radius $l$ that form a finitely overlapping covering of $\R^d \setminus B_{4l}(0)$. Without loss of generality, we only consider balls that intersect $\R^d \setminus B_{4l}(0)$, thus 
$$\R^d \setminus B_{4l}(0) \subset \bigcup_{i\in I} B_i \subset \R^d \setminus B_{2l}(0).$$
Let $D_0$ be a cylindrical annulus
$$D_0:= \{(x,t)\in \R^n, \, x\in B_{4l} \setminus B_{2l}, \, a|t| \leq 2l\} \cup \{(x,t)\in \R^n ,\,  x\in B_{2l}, \, l \leq a|t| \leq 2l\}.$$
The bound \eqref{lem7.3b} becomes then
\begin{align} \label{lem7.3c}
\frac1l \iint_\Omega & |u|^q \dr_r[-\chi_l^k] \frac{dt}{|t|^{n-d-1}}\, dx \nn \\ 
& \lesssim \frac1{l^2} \iint_{(x,t) \in D_0}  |u|^q \frac{dt}{|t|^{n-d-1}}\, dx + \frac1{l^2} \sum_{i\in I} \int_{x\in B_i} \int_{a|t| \leq 2l} |u|^q \frac{dt}{|t|^{n-d-1}}\, dx \nn \\
& \quad = \frac1{l^2} \iint_{(x,t) \in D_0}  \left||u|^{q/2-1}u\right|^2 \frac{dt}{|t|^{n-d-1}}\, dx + \frac1{l^2} \sum_{i\in I} \int_{x\in B_i} \int_{a|t| \leq 2l} \left||u|^{q/2-1}u\right|^2 \frac{dt}{|t|^{n-d-1}}\, dx.
\end{align}

Now, since $u$ is an energy solution, there exists $l_0$ such that $\supp \, \Tr u \subset B_{2l_0}(0)$. If $l\geq l_0$, all the cubes $\{(x,t)\in \R^n, \, x\in B_i,\, a|t| \leq 2l\}$ and the domain $D_0$ intersect the boundary $\dr \Omega = \R^d$ where $\Tr u$ is 0.Thus $\Tr |u|^{q/2-1}u = 0$ by boundary Moser's estimate (when $q> 2$) or Holder's inequality (when $q<2$). So by Poincar\'e's inequality (see Lemma 4.13 in \cite{DFMprelim} (the proof there is written when the domains are balls but it goes through for our domains), we have
\begin{align} \label{lem7.3d}
\frac1l \iint_\Omega & |u|^q \dr_r[-\chi_l^k] \frac{dt}{|t|^{n-d-1}}\, dx \nn \\ 
& \lesssim  \iint_{(x,t) \in D_0}  |\nabla [|u|^{q/2-1}u]|^2 \frac{dt}{|t|^{n-d-1}}\, dx + \sum_{i\in I} \int_{x\in B_i} \int_{a|t| \leq 2l} |\nabla [|u|^{q/2-1}u]|^2 \frac{dt}{|t|^{n-d-1}}\, dx \\
& \lesssim  \iint_{(x,t) \in D_0}  |\nabla u|^2 |u|^{q-2} \frac{dt}{|t|^{n-d-1}}\, dx + \sum_{i\in I} \int_{x\in B_i} \int_{a|t| \leq 2l} |\nabla u|^2 |u|^{q-2} \frac{dt}{|t|^{n-d-1}}\, dx. \nn
\end{align}
where the last line is due to \eqref{nablau&v}. Due to the finite overlapping of the covering $(B_i)_{i\in I}$ and the fact that we always avoid $\{(x,t) \in \Omega, \, x\in B_{2l}(0), \, a|t| \leq l\}$ in the integrations, we have, when $l\geq l_0$,
$$\frac1l \iint_\Omega |u|^q \dr_r[-\chi_l^k] \frac{dt}{|t|^{n-d-1}}\, dx \lesssim \iint_{\Omega \setminus \{(x,t) \in \Omega, \, x\in B_{2l}(0), \, a|t| \leq l\}} |\nabla u|^2 |u|^{q-2} \frac{dt}{|t|^{n-d-1}}\, dx.$$
Thanks to Lemma \ref{lemES1}, the right-hand side above converges to 0 as $l$ goes to $+\infty$. The lemma follows.
\ep

The following lemma is the key point in proving the existence of solutions to the Dirichlet problem, when the boundary function (trace) is smooth. 

\begin{lem}[Global estimate $N< \Tr$ for energy solutions] \label{lem7.4}
Let $\cL$ be an elliptic operator satisfying $(\cH^1_\kappa)$ for some constant $\kappa\geq 0$. Let $a>0$, $q\in (q_0,q'_0)$ where $q_0$ is given by Proposition \ref{prop1.1}. There exist two values $\kappa_0>0$ and $C>0$, both depending only on $a$, $q$, $n$, $\lambda_q(\cA)$, $\|\cA\|_\infty$, $\|b\|_{\infty} + \|b^{-1}\|_\infty$ such that if $\kappa \leq \kappa_0$, then for any energy solution $u\in W$ to $\cL u = 0$,
\begin{equation} \label{lem7.4a} 
\|\wt N_{a,q}(u)\|_q \leq C \|\Tr u\|_q. \end{equation}
\end{lem}

\bp Let us fix for the proof some $k$ large, say $k=20$. In addition, we always consider $\kappa_0$ smaller than 1 and $\kappa \leq \kappa_0$. Assume that $u$ is not trivially zero, otherwise there is nothing to prove.

The combination of Lemma \ref{lem7.1} and Lemma \ref{lemS<Ne} gives that for $l>0$,
\begin{equation} \label{lem7.4b} \begin{split}
\|\wt N_{a,q}(u|\chi_{l}^{k+12})\|_q^q
& \leq  C \kappa \|\wt N_{a,q} (u|\chi_{l}^{k-2})\|_{q}^q + C \|\Tr u\|_q^q   \\
& \hspace{3cm} + \, \frac{C}{q} \, \iint_{\Omega} \left(\frac{|u|^q}{|t|} - \dr_r [|u|^q] \right) \dr_r[\chi_l^{k}] \, \frac{dt}{|t|^{n-d-2}}\, dx,
\end{split} \end{equation}
where the constant $C>0$ depends only on $a$, $q$, $n$, $\lambda_q(\cA)$, $\|\cA\|_\infty$ and $\|b\|_{\infty} + \|b^{-1}\|_\infty$. In particular note that $k=20$ is fixed and $\kappa \leq \kappa_0 \leq 1$, the constant $C$ does not depend $\kappa$ and $l$.

\medskip

In order to control $\|\wt N_{a,q} (u|\chi_{l}^{k-2})\|_{q}$ by $\|\wt N_{a,q}(u|\chi_{l}^{k+12})\|_q$, we make an additional assumption on $\phi$, the smooth function from which $\chi_l$ is defined. Recall that $\phi \in C^\infty_0([0,\infty))$, $0\leq \phi \leq 1$, $\phi \equiv 1$ on $[0,1]$, $\phi \equiv 0$ outside $[0,2]$, $\phi$ is decreasing and $|\phi'| \leq 2$. We assume in addition that 
\begin{equation}\label{ex:asmpphi}
	\phi(x) = 2-x \text{ when } x\in \left(\frac98,\frac{15}8\right).
\end{equation} 
With this additional assumption, it is not hard to verify that 
$$\chi_l^{k-2} \lesssim \chi_l^{k+12} + l \dr_r [-\chi_{3l/2}^k],$$
where the constant is universal (recall $k=20$ is fixed). Indeed, observe that the assumption \eqref{ex:asmpphi} in particular implies that $\dr_r[-\chi_{3l/2}^k]$ has a strictly positive lower bound on the interval $\frac{3}{2} \frac{l}{a} \left( \frac{9}{8}, \frac{15}{8} \right) \supset \left[ \frac{15}{8} \frac{l}{a}, 2 \frac{l}{a} \right]$, where the value of $\chi_l$ is very small.
It then follows that
\begin{equation} \label{lem7.4c}
\|\wt N_{a,q}(u|\chi_{l}^{k-2})\|_q^q \lesssim \|\wt N_{a,q}(u|\chi_{l}^{k+12})\|_q^q + l \int_{x\in \R^d} \left|\sup_{\Gamma_a(x)} (u|\dr_r [-\chi_{3l/2}^k])_{\W,a,q}\right|^q \, dx.
\end{equation}
Since $\dr_r [-\chi_{3l/2}^k]$ is supported in a domain where $|s| \approx l/a$, we deduce 
\[\begin{split}
\left|\sup_{\Gamma_a(x)} (u| \dr_r [-\chi_{3l/2}^k])_{\W,a,q}\right|^q & \lesssim \left|(u| \dr_r [-\chi_{3l/2}^k])_{\W,10a,q}(x,3l/2a)\right|^q \\
& \lesssim l^{-n} \int_{y\in B_{8l}(x)} \int_{3l/4 \leq a|s| \leq 3l} |u|^q \dr_r [-\chi_{3l/2}^k] \, ds\, dy \\
& \lesssim l^{-d-1} \int_{y\in B_{8l}(x)} \int_{s\in \R^{n-d}} |u|^q \dr_r [-\chi_{3l/2}^k] \, \frac{ds}{|s|^{n-d-1}}\, dy.
\end{split}\]
And then by Fubini's lemma,
$$\int_{x\in \R^d} \left|\sup_{\Gamma_a(x)} (u| \dr_r [-\chi_{3l/2}^k])_{\W,a,q}\right|^q \, dx \lesssim l^{-1} \iint_{(y,s) \in \Omega}  |u|^q \dr_r [-\chi_{3l/2}^k] \, \frac{ds}{|s|^{n-d-1}}\, dy.$$
The estimate \eqref{lem7.4c} becomes
\begin{equation} \label{lem7.4d}
\|\wt N_{a,q}(u|\chi_{l}^{k-2})\|_q^q \lesssim \|\wt N_{a,q}(u|\chi_{l}^{k+12})\|_q^q +\iint_{(x,t) \in \Omega}  |u|^q \dr_r [-\chi_{3l/2}^k] \, \frac{dt}{|t|^{n-d-1}}\, dx.
\end{equation}
We eventually want to estimate the last term by $\|\wt N_{a,q}(u|\chi_{l}^{k+12})\|_q^q$. Notice that for any $l' \leq \left(\frac 23 \right)^2 l<\frac{l}{2}$, we have
\begin{equation} \label{lem7.4e} \begin{split}
\iint_{(x,t) \in \Omega}  |u|^q \dr_r [-\chi_{l'}^k] \, \frac{dt}{|t|^{n-d-1}} \, dx & \lesssim (l')^{d-n} \int_{x\in \R^d} \int_{l'/2 \leq a|t| \leq 2l'} |u|^q \, dt \, dx \\
& \lesssim  \int_{x\in \R^d} \left( (l')^{-n} \int_{y\in B_{l'/2}(x)} \int_{l'/2 \leq a|t| \leq 2l'} |u|^q \, dt \, dy \right) dx \\
& \lesssim \int_{x\in \R^d} |\wt N_{a,q}(u|\chi_l^{k+12})|^q \, dx = \|\wt N_{a,q}(u|\chi_{l}^{k+12})\|_q^q.
\end{split} \end{equation}
The last inequality is because $\chi_l \equiv 1$ when $a|t| \leq l$, so in particular when $a|t| \leq 2l' < l$.
Let $w:\, (0,\infty) \to \R$ denote the continuous function
$$w(l) = \iint_{(x,t) \in \Omega}  |u|^q \dr_r [-\chi_{l}^k] \, \frac{dt}{|t|^{n-d-1}}\, dx,$$
then the combination of \eqref{lem7.4b}, \eqref{lem7.4d}, and \eqref{lem7.4e} gives that
\begin{align} 
 \sup_{l' \leq 4l/9} w(l') \lesssim & \|\wt N_{a,q}(u|\chi_{l}^{k+12})\|_q^q \nn \\
& \qquad  \leq  C \kappa \|\wt N_{a,q}(u|\chi_{l}^{k+12})\|_q^q+  C \kappa \, w (3l/2) + C \|\Tr u\|_q^q 
\nn \\
& \hspace{4cm} + \, \frac{C}{q} \, \iint_{\Omega} \left(\frac{|u|^q}{|t|} - \dr_r [|u|^q] \right) \dr_r[\chi_l^{k}] \, \frac{dt}{|t|^{n-d-2}}\, dx, \label{ex:tmp3} 
\end{align}
where the constant $C>0$ is independent of $\kappa$.

The last term of the above estimate \eqref{ex:tmp3} can also be written in terms of the function $w(l)$. In fact, we define 
$$v(r) = \int_{x\in \R^d} \int_{\theta \in \mathbb S^{n-d-1}} |u(x,r\theta)|^q d\theta \, dx,$$
which is finite for almost every $r>0$ since $u\in L_{loc}^q(\Omega)$. On one hand, the last term of \eqref{ex:tmp3} can be rewritten by polar coordinates
\begin{equation} \label{cond4}
\iint_{\Omega} \left(\frac{|u|^q}{|t|} - \dr_r [|u|^q] \right) \dr_r[\chi_l^{k}] \, \frac{dt}{|t|^{n-d-2}}\, dx =
\int_{0}^\infty [rv'(r) - v(r)] \, \dr_r[- \chi_l^{k}] \, dr.
\end{equation}
On the other hand, we have 
\begin{equation}\label{eq:wl}
	w(l) = \iint_{(x,t) \in \Omega}  |u|^q \dr_r [-\chi_{l}^k] \, \frac{dt}{|t|^{n-d-1}}\, dx = \int_{0}^\infty v(r) \dr_r[-\chi_l^k] \, dr. 
\end{equation}
By the construction of $\chi_l$, we can write $\dr_r [-\chi_l^k]$ as $\frac al\xi\left(\frac{ar}{l}\right)$ with a non-negative function $\xi = -\phi^k \phi' \in C^\infty_0([0,\infty))$, and hence
\begin{align}
w'(l) & = \int_{0}^{+\infty} v(r) \dr_l \left[ \frac{a}{l} \xi\left(\frac{ar}{l}\right)\right] \, dr \nn \\
& = -\frac{1}{l} \int_{0}^{+\infty} v(r) \frac{a}{l} \xi \left(\frac{ar}{l}\right) \, dr -\frac{1}{l} \int_{0}^{+\infty} r v(r) \frac{a^2}{l^2} \xi' \left(\frac{ar}{l}\right) \, dr \nn \\
& = -\frac{1}{l} \int_{0}^{+\infty} v(r) \frac{a}{l} \xi \left(\frac{ar}{l}\right) \, dr -\frac{1}{l} \int_{0}^{+\infty} r v(r) \frac{a}{l} \dr_r \left[ \xi\left(\frac{ar}{l}\right) \right] \, dr \nn \\
& = -\frac{1}{l} \int_{0}^{+\infty} v(r) \frac{a}{l} \xi \left(\frac{ar}{l}\right) \, dr +  \frac{1}{l} \int_{0}^{+\infty}  \dr_r[rv(r)] \frac{a}{l} \xi\left(\frac{ar}{l}\right) \, dr \nn \\
& = \frac{1}{l} \int_{0}^{+\infty}  r v'(r) \, \frac{a}{l} \xi\left(\frac{ar}{l}\right) \, dr \nn \\
& = \frac{1}{l} \int_{0}^{+\infty}  r v'(r) \dr_r [-\chi_l^k] \, dr \label{eq:wlder}
\end{align}
where we use the integration by part and the fact that $\xi$ has compact support. By combining \eqref{cond4}, \eqref{eq:wl} and \eqref{eq:wlder}, we can rewrite \eqref{ex:tmp3} as
\begin{align} 
 & \sup_{l' \leq 4l/9} w(l')  \nn \\
 & \qquad \lesssim  \|\wt N_{a,q}(u|\chi_{l}^{k+12})\|_q^q \nn \\
& \qquad \qquad \leq  C \kappa \|\wt N_{a,q}(u|\chi_{l}^{k+12})\|_q^q+  C \kappa \, w (3l/2) + \frac{C}{q} \left( l\, w'(l)-w(l) \right) + C \|\Tr u\|_q^q. \label{ex:tmp3new} 
\end{align}
We claim that for any $l_0>0$, there exists $l\geq l_0$ such that
\begin{equation}\label{claim:taildecay}
	w\left( \frac{3}{2}l \right) \leq \left(\frac32 \right)^5 \sup_{l' \leq \frac{4}{9}l} w(l') \quad \text{and} \quad l\, w'(l)-w(l) \leq 0. 
\end{equation} 
Assume the claim holds, then there is a sequence $l_i \to \infty$ such that
\[ C\kappa w\left(\frac 32 l_i \right) \leq C\kappa \left(\frac32\right)^5 \sup_{l'\leq \frac 49 l_i} w(l') \leq \frac 13 \|\wt N_{a,q}(u|\chi_{l_i}^{k+12})\|_q^q. \]
Thus by choosing $\kappa$ sufficiently small satisfying $C\kappa(3/2)^5 \leq 1/3$, we obtain \eqref{lem7.4a} by combining \eqref{ex:tmp3new} and \eqref{claim:taildecay}.

Recall in Lemma \ref{lem7.3} we proved \eqref{lem7.3a}, which says $w(l)/l\to 0$ as $l\to +\infty$. Thus it is not difficult to show each estimate of the claim \eqref{claim:taildecay} holds for infinitely many $l$'s individually. But we need to find a sequence $l_i\to \infty$ such that both estimates hold at the same time. We prove the claim by contradiction: Assume not, then there exists $l_0>0$ such that for all $l\geq l_0$, either
\begin{equation}\label{claim:contra1}
	\quad w\left( \frac{3}{2}l \right) > \left(\frac32 \right)^5 \sup_{l' \leq \frac{4}{9}l} w(l') 
\end{equation} 
or
\begin{equation}\label{claim:contra2}
	l\, w'(l)-w(l) > 0. 
\end{equation}

We define $M$ as the positive quantity
$$M:= \inf_{\left[\left( \frac 23 \right)^4 l_0,l_0 \right]} w >0.$$
Indeed, note that $w(l)$ is continuous, if $M$ is zero, then there exists $l\in [16l_0/81, l_0]$ such that $w(l) = 0$. Hence $u(x,t) \equiv 0$ almost everywhere whenever $t\in \supp \, \dr_r(-\chi_l^k)$, which is a non-trivial band contained in $\{l \leq a|t| \leq 2l\}$. Extending $u$ by zero outside of this band, we get another weak solution in $W$ to $\cL u=0$ with the same boundary value. By the uniqueness of weak solutions (see Lemma \ref{lemDirichlet}), we conclude that $u(x,t) \equiv 0$ almost everywhere beyond the band, in particular whenever $|t| \geq 2l_0/a$. Then \eqref{claim:contra1} and \eqref{claim:contra2} clearly fail, and thus $M$ is strictly positive. 

Let $I_k$ denote the interval $\left[\left( \frac 32 \right)^{k-4} l_0, \left(\frac 32 \right)^{k-3} l_0 \right)$. We will show by induction that the following property holds for all $k\in \N$
$$\mathcal P(k): \quad \text{ there exists } l \in I_k \text{ such that } w(l) \geq \frac{l}{l_0} M.$$
The base steps when $k=0$, $1$, $2$, $3$ are immediate by the definition of $M$. Assume that $\mathcal P(k)$ and $\mathcal P(k+3)$ hold, that is,
$$\text{ there exists } l_1 \in I_k  \text{ such that } w(l_1) \geq \frac{l_1}{l_0} M,$$
and
$$\text{ there exists } l_2 \in I_{k+3} \text{ such that } w(l_2) \geq \frac{l_2}{l_0} M.$$
We want to prove $\mathcal P(k+4)$ holds. Let $l_3 = \left( \frac 32 \right)^k l_0$ denote the upper endpoint for the interval $I_{k+3}$. We distinguish two cases. Either there exists $l \in [l_2,l_3) \subset I_{k+3}$ such that \eqref{claim:contra1} is satisfied, that is,
$$w\left( \frac 32 l \right) > \left( \frac 32 \right)^5 \sup_{l' \leq \frac 49 l} w(l').$$
Since $\frac 49 l \geq \sup I_{k} \geq l_1$, $\frac 32 l\in I_{k+4}$ and $\left( \frac 32 \right)^5 l_1 \geq \sup I_{k+4} > \frac 32 l$, it follows that
\[ w\left( \frac 32 l \right) > \left( \frac 32 \right)^5  w(l_1) \geq \left( \frac 32 \right)^5 \frac{l_1}{l_0} M \geq \frac{\frac 32 l}{l_0} M.  \]
Therefore $\mathcal{P}(k+4)$ holds for $\frac 32 l \in I_{k+4}$.
Or alternatively in the second case, \eqref{claim:contra1} is never satisfied on the interval $[l_2, l_3)$ and so by assumption \eqref{claim:contra2} is always satisfied on the same interval, that is,
$$l\, w'(l) - w(l) > 0 \quad \text{ for all } l \in [l_2,l_3).$$
By Gr\"onwall's inequality and the continuity of $w$, this implies
\[ w(l_3) \geq \frac{l_3}{l_2}\, w(l_2) \geq \frac{l_3}{l_0} M, \]
i.e. $\mathcal P(k+4)$ holds for $l_3 \in I_{k+4}$. The induction step follows, and we conclude that $\mathcal{P}(k)$ holds for all $k\in \mathbb{N}$. Thus in particular, we have
\[ \limsup_{l\to\infty} \frac{w(l)}{l} \geq \frac{M}{l_0} >0, \]
contradiction Lemma \ref{lem7.3}. Therefore the claim \eqref{claim:taildecay} holds.
\ep

The next result proves the existence of a solution to the Dirichlet problem with boundary value in $L^q(\R^d)$, by approximating using energy solutions.

\begin{lem} \label{lem7.5}
Let $\cL$ be an elliptic operator satisfying $(\cH^1_\kappa)$ for some constant $\kappa\geq 0$. Let $a>0$, $q\in (q_0,q'_0)$ where $q_0$ is given by Proposition \ref{prop1.1}. There exists two values $\kappa_0>0$ and $C>0$, both depending only on $a$, $q$, $n$, $\lambda_q(\cA)$, $\|\cA\|_\infty$, $\|b\|_{\infty} + \|b^{-1}\|_\infty$ such that if $\kappa \leq \kappa_0$, then for any $g\in L^q(\R^d)$, there exists a weak solution $u\in W^{1,2}_{loc}(\Omega)$ to $\cL u = 0$ such that
\begin{equation} \label{lem7.5a} 
   \lim_{(z,r)\in \Gamma_a(x) \atop {r\to 0}} \frac{1}{|\W_a(z,r)|} \iint_{W_a(z,r)} |u(y,s) - g(x)|^q \, dy\, ds  = 0
\end{equation}
for almost every $x\in\R^d$,
and
\begin{equation} \label{lem7.5b} 
\|\wt N_{a,q}(u)\|_q \leq C \|g\|_q, 
\end{equation}
where the constants depend only on $a$, $q$, $n$, $\lambda_q(\cA)$, $\|\cA\|_\infty$, $\|b\|_{\infty} + \|b^{-1}\|_\infty$.
\end{lem}
\begin{rmk}
	\begin{enumerate}
		\item The trace operator is defined for functions in $W$. In general for functions in $W^{1,2}_{loc}(\Omega)$, we use the equality \eqref{lem7.5a} as a weaker version or an adaptation of $\Tr u = g$ almost everywhere on $\R^d$.
		\item By Moser's estimate or H\"older inequality, the Lebesgue density property \eqref{lem7.5a} with power $q$ is equivalent to that of power $2$.
	\end{enumerate}
\end{rmk}

\bp
Since $C^\infty_0(\R^d)$ is dense in $L^q(\R^d)$, we can find a collection $(g_i)_{i}$ of functions in $C^\infty_0(\R^d)$ such that $g_i \to g$ in $L^q(\R^d)$. In particular $(g_i)_i$ is a Cauchy sequence in $L^q(\R^d)$. By passing to a subsequence if necessary, we may assume that $\|g_i-g_{i+1}\|_q \leq 2^{-2i}$ for all $i$. 
Since $C^\infty_0(\R^d) \subset H$, Lemma \ref{lemDirichlet} guarantees the existence of an energy solution $u_i \in W$ such that $\Tr u_i = g_i$. 

Lemma \ref{lemDirichlet} also guarantees that for boundary value $g_i - g_j \in H$, the corresponding solution in $W$ is unique, so by the linearity of the operator $\cL$ the solution is $u_i - u_j$. Thus by Lemma \ref{lem7.4}, we have for all $i, j$,
\begin{equation} \label{lem7.5d}
\|\wt N_{a,q}(u_i)\|_q \lesssim \|g_i\|_q,
\end{equation}
and 
\begin{equation}\label{ex:Cseq}
	\|\wt N_{a,q}(u_i-u_j)\|_q \lesssim \|g_i-g_j\|_q.
\end{equation}
We claim that $\wt N_{a,q}(u_i-u_{i+1})(x)\to 0$ as $i\to\infty$, for almost every $x\in \R^d$. In fact, for any $\lambda >0$, denote
\[ E_\lambda = \{x\in\R^d, \wt N_{a,q}(u_i-u_{i+1})(x) < \frac{\lambda}{2^i} \text{ for all } i \}. \]
Then by \eqref{ex:Cseq} and the assumption on $(g_i)_i$, we have
\begin{align*}
	|E_\lambda^c| \leq \sum_i \left|\left\{x\in \R^d, \wt N_{a,q}(u_i-u_{i+1})(x) \geq \frac{\lambda}{2^i} \right\} \right| & \leq \sum_i \left(\frac{2^{i}}{\lambda} \right)^q \int_{x\in\R^d} |\wt N_{a,q}(u_i-u_{i+1})(x)|^q dx \\
	& \leq \sum_i \left( \frac{2^i } {\lambda} 2^{-2i} \right)^q \\
	& \leq \frac{C}{\lambda^q}.
\end{align*}
Consider an arbitrary sequence $\lambda_j \to \infty$, the above estimate implies $\bigcap_{\lambda_j} E_{\lambda_j}^c$ has measure zero, and thus its complement $\bigcup_{\lambda_j} E_{\lambda_j}$ has full measure in $\R^d$. Therefore for almost every $x\in \R^d$, there exists some $\lambda>0$ (depending on $x$) such that
\[ \wt N_{a,q}(u_i - u_{i+1})(x) < \frac{\lambda}{2^i} \quad \text{ for all }i. \]
Thus by the triangle inequality for the $L^q$ norm, we have that for any $j>i$,
\[  \wt N_{a,q}(u_i - u_{j})(x)\leq \frac{2\lambda}{2^{i}}, \]
and by the definition of $\wt N_{a,q}$, this means that for any $(z,r)\in \Gamma_a(x)$,
\begin{equation}\label{ex:CseqW}
	\|u_i - u_j\|_{L^q(W_a(z,r))} = \left(\iint_{W_a(z,r)} |u_i - u_j|^q dy\, ds \right)^{1/q} \leq \frac{2\lambda}{2^{i}} |W_a(z,r)|^{1/q}. 
\end{equation} 
That is to say $(u_i)_i$ is a Cauchy sequence in $L^q(W_a(z,r))$, and thus it converges to a function $u$. Since this holds for almost every $x\in\R^d$ and any $(z,r)\in \Gamma_a(x)$, this function $u$ is well defined on all of $\Omega$ and $u_i \to u$ in $L_{loc}^q$. Moreover, by passing $j\to \infty$ in \eqref{ex:CseqW}, we have
\begin{equation}\label{ex:Ntgcv}
	\wt N_{a,q}(u_i - u)(x) = \sup_{(z,r)\in \Gamma_a(x)} \left( \frac{1}{|W_a(z,r)|} \iint_{W_a(z,r)} |u_i - u|^q dy\, ds \right)^{1/q} \leq \frac{2\lambda}{2^i},
\end{equation} 
which converges to zero as $i\to \infty$. In particular for any $(z,r)\in \Gamma_a(x)$,
\[ \left| (u_i)_{W,a,q}(z,r) - (u)_{W,a,q}(z,r) \right| \leq \wt N_{a,q}(u_i-u)(x) \leq \frac{2\lambda}{2^i}. \] We deduce that
\begin{equation}\label{ex:Ntgcv2}
	\lim_{i\to \infty} \wt N_{a,q}(u_i)(x) = \wt N_{a,q}(u)(x)
\end{equation} 
for almost every $x\in \R^d$. The estimate \eqref{lem7.5b} follows by \eqref{ex:Ntgcv2}, Fatou's lemma and \eqref{lem7.5d}.

We recall the interior Cacciopoli inequality (see Lemma 8.26 in \cite{DFMprelim}): let $B$ be a fixed ball with radius $r>0$, such that the distance from $4B$ to the boundary is roughly $r$, we know
\[ \iint_{B} |\nabla (u_i - u_j)|^2  dm \lesssim  \frac{1}{r^2} \iint_{2B} |u_i-u_j|^2 dm < + \infty. \]
By assumption the distances from $B$ and $2B$ to the boundary are roughly $r$, so the above estimate is equivalent to
\begin{equation}\label{ex:cvgrad}
	\iint_B |\nabla u_i - \nabla u_j |^2 dx \, dt \lesssim \frac{1}{r^2} \iint_{2B} |u_i - u_j|^2 dx \, dt.
\end{equation} 
We have shown that $(u_i)_i$ is a Cauchy sequence in $L^q(4B)$. Either by H\"older's inequality (in the case of $q\geq 2$) or by Moser's estimate (in the case of $q\leq 2$), see \eqref{eq:Moserlt2B}, it follows that $(u_i)_i$ is a Cauchy sequence in $L^2(2B)$. Therefore \eqref{ex:cvgrad} implies that $(\nabla u_i)_i$ is also a Cauchy sequence in $L^2(B)$, and thus it converges in $L^2(B)$. By the uniqueness of the limit $\nabla u_i \to \nabla u$ in $L_{loc}^2(\Omega)$.
The convergence forces $u$ to be, like the $u_i$'s, a weak solution to $\cL u = 0$.

We now turn to the proof of \eqref{lem7.5a}. Since $u_i\in W$, by Theorem 3.13 of \cite{DFMprelim} (see also step 2 of Lemma \ref{lem7.1}), we have
\begin{equation}\label{ex:Ntglm2}
	\sup_{(z,r)\in \Gamma_a(x) \atop {r<\delta}} \frac{1}{|W_a(z,r)|} \iint_{W_a(z,r)} |u_i(y,s) - g_i(x)|^2 dy\, ds \to 0 \quad \text{ as } \delta \to 0, 
\end{equation} 
for almost every $x\in \R^d$. The set such that \eqref{ex:Ntglm2} holds depends {\em a priori } on $i$, but since a countable union of sets of zero measure is still a set of measure zero, we have \eqref{ex:Ntglm2} for every $i\in \N$ and almost every $x\in \R^d$.
Either by H\"older's inequality (when $q\leq 2$) or by Moser estimate (when $q\geq 2$, see \eqref{eq:Mosergt2}) applied to the weak solution $u_i-g_i(x)$, it follows that
\begin{equation}\label{ex:Ntglm3}
	\sup_{(z,r)\in \Gamma_a(x) \atop {r<\delta}} \frac{1}{|W_a(z,r)|} \iint_{W_a(z,r)} |u_i(y,s) - g_i(x)|^q dy\, ds \to 0 \quad \text{ as } \delta \to 0, 
\end{equation} 
for any $i\in \N$ and almost every $x\in \R^d$.
Since $g_i \to g$ in $L^q(\R^d)$, by passing to a subsequence $g_i(x) \to g(x)$ for almost every $x\in\R^d$. Recall also that \eqref{ex:Ntgcv2} holds for almost every $x\in \R^d$. As a consequence, \eqref{lem7.5a} is obtained by taking the limit as $i\to \infty$ in \eqref{ex:Ntglm3}.
\ep

\section{Uniqueness in the Dirichlet problem}

\label{sUniqueness}

In this section, we assume that $u \in W_{loc}^{1,2}(\Omega)$ is a weak solution to $\cL u = 0$ - where $\cL$ satisfies $(\cH^1_\kappa)$ - and that $u$ satisfies $\|\wt N_{a,q}(u)\|_q < +\infty$. We want to prove that if $\kappa$ is sufficiently small and that
\begin{equation} 
   \lim_{(z,r)\in \Gamma_a(x) \atop {r\to 0}} \frac{1}{|\W_a(z,r)|} \iint_{W_a(z,r)} u(y,s) \, dy\, ds  = 0 \qquad \text{ for almost every } x\in \R^d, 
\end{equation}
then $u$ has to be $0$. This, in turn, proves the uniqueness of solution.

\begin{lem} \label{lem8.1}
Let $\cL$ be an elliptic operator satisfying $(\cH^1_\kappa)$ for some constant $\kappa\geq 0$. Let $a>0$, $q\in (q_0,q'_0)$ where $q_0$ is given by Proposition \ref{prop1.1}. For any weak solution $u\in W^{1,2}_{loc}(\Omega)$ to $\cL u =0$ that satisfies $\|\wt N_{a,q}(u)\|_q < +\infty$, we have
\begin{equation} \label{lem8.1a} 
\|\wt N_{a,q}(u)\|_q \approx  \|S_{a,q}(u)\|_q, 
\end{equation}
where the constants depend only on $a$, $q$, $n$, $\lambda_q(\cA)$, $\|\cA\|_\infty$, $\|b\|_{\infty} + \|b^{-1}\|_\infty$ and $\kappa$.
\end{lem}

\bp
Since $\|\wt N_{a,q}(u) \|_q$ is finite, the proof of
$$\|S_{a,q}(u)\|_q \leq C \|\wt N_{a,q}(u)\|_q < +\infty $$
can be achieved by using Lemma \ref{lemS<Nd} with an increasing sequence of cut-off functions $\Psi_i(x)$ such that $\Psi_i (x) \to 1$ for all $x\in \Omega$.

With the notation of Section \ref{N<S}, for a fixed $k=15$, for any $\epsilon>0$, $l>0$, and any ball $B$ of radius $l'\geq l$, Lemma \ref{lem6.10} shows that
\begin{equation} \label{lem8.1b}
\|\wt N_{a,q}(u|\Psi_\epsilon^k\Psi_{B,l}^k)\|_{q}^q \lesssim \|S_{a,q}(u|\Psi_\epsilon^{k-12}\Psi_{B,l}^{k-12})\|_{q}^q + \int_{y\in \R^d} \int_{s\in \R^{n-d}} |u|^q \Psi_\epsilon^{k-3}  \dr_r[-\Psi_{B,l}^{k-3}] \frac{ds}{|s|^{n-d-1}}\, dy.\end{equation}
By Fubini's lemma and the fact that $B$ has radius $l'\geq l$, the second term in the right-hand side is bounded by
\begin{align*}
	\int_{y\in \R^d} \fint_{l \leq a|s| \leq 2l} |u|^q ds\, dy & \lesssim \int_{x\in \R^d } \left( \frac{1}{|W_a(x, l/a )|} \iint_{(y,s) \in W_a(x, l/a )} |u|^q ds \, dy \right) dx \\
	& \leq \|N_{a,q}(u)\|_{L^q(\R^d \setminus B/2)}^q \rightarrow 0,
\end{align*}
which tends to zero as the ball $B \nearrow \R^d$.
By taking the limit as $l\to \infty$ and $\epsilon \to 0$, we obtain then
$$\|\wt N_{a,q}(u)\|_q \leq C \|S_{a,q}(u)\|_q.$$
The lemma follows.
\ep

\begin{lem} \label{lem8.11}
Let $\cL$ be an elliptic operator satisfying $(\cH^1_\kappa)$ for some constant $\kappa\geq 0$. Let $a>0$, $q\in (q_0,q'_0)$ where $q_0$ is given by Proposition \ref{prop1.1}. Suppose that $g$ lies in $L^q(\R^d)$ and $u\in W^{1,2}_{loc}(\Omega)$ is weak solution to $\cL u = 0$ that satisfy both $\|\wt N_{a,q}(u)\|_q < \infty$ and
\begin{equation} \label{lem8.11a} 
   \lim_{(z,r)\in \Gamma_a(x) \atop {r\to 0}} \frac{1}{|\W_a(z,r)|} \iint_{W_a(z,r)} u(y,s) \, dy\, ds  = g(x) \qquad \text{ for a.e. } x\in \R^d.
\end{equation}
Then \eqref{lem8.11a} self-improves itself into the $q$-Lebesgue property
\begin{equation} \label{lem8.11b} 
   \lim_{(z,r)\in \Gamma_a(x) \atop {r\to 0}} \frac{1}{|\W_a(z,r)|} \iint_{W_a(z,r)} |u(y,s) - g(x)|^q \, dy\, ds  = 0 \qquad \text{ for a.e. } x\in \R^d.
\end{equation}
\end{lem}

\bp
Let us write $\bar u(z,r)$ for $\frac{1}{|\W_a(z,r)|} \iint_{W_a(z,r)} u(y,s) \, dy\, ds$. Observe that 
\[\begin{split}
& \left(\frac{1}{|\W_a(z,r)|} \iint_{W_a(z,r)} |u(y,s) - g(x)|^q \, dy\, ds \right)^\frac1q \\
& \hspace{3cm} \leq \left(\frac{1}{|\W_a(z,r)|} \iint_{W_a(z,r)} |u(y,s) - \bar u(z,r)|^q \, dy\, ds \right)^\frac1q + |\bar u(z,r) - g(x)| \\
& \hspace{3cm} := T_1 + T_2.
\end{split}\]
Due to \eqref{lem8.11a}, the quantity $T_2$ tends to 0 as $(z,r)\in \Gamma_a(x)$, $r\to 0$ for almost every $x\in \R^d$. It remains to prove that 
\begin{equation}\label{claim811}
 \lim_{(z,r)\in \Gamma_a(x) \atop {r\to 0}} \left(\frac{1}{|\W_a(z,r)|} \iint_{W_a(z,r)} |u(y,s) - \bar u(z,r)|^q \, dy\, ds \right)^\frac1q = 0 \qquad \text{ for a.e. } x\in \R^d.
 \end{equation}
We split the proof of the claim into two cases: $q\geq 2$ and $q\leq 2$. 

\medskip

First, let us treat the case $q\geq 2$. Thanks to Proposition \ref{lemMoserI} and then the Poincar\'e inequality, we obtain that
\[\begin{split}
T_1
& \lesssim \left(\frac{1}{|\W_{4a}(z,r/2) \cup \W_{4a}(z,2r)|} \iint_{\W_{4a}(z,r/2) \cup \W_{4a}(z,2r)} |u(y,s) - \bar u(z,r)|^2 \, dy\, ds \right)^\frac12 \\
& \lesssim r \left(\frac{1}{|\W_{4a}(z,r/2) \cup \W_{4a}(z,2r)|} \iint_{\W_{4a}(z,r/2) \cup \W_{4a}(z,2r)} |\nabla u(y,s)|^2 \, dy\, ds \right)^\frac12. \\
\end{split}\]
Notice that $\W_{4a}(z,r/2) \simeq a^dr^n \simeq \W_{4a}(z,2r)$, from which we deduce that
\[\begin{split}
T_1 & \lesssim \left( \iint_{\W_{4a}(z,r/2)} |\nabla u(y,s)|^2 \, dy\, \frac{ds}{|s|^{n-2}} + \iint_{\W_{4a}(z,2r)} |\nabla u(y,s)|^2 \, dy\, \frac{ds}{|s|^{n-2}}  \right)^\frac12 \\
& \lesssim \left( \iint_{\W_{8a}(x,r/2)} |\nabla u(y,s)|^2 \, dy\, \frac{ds}{|s|^{n-2}} \right)^\frac12 + \left(\iint_{\W_{8a}(x,2r)} |\nabla u(y,s)|^2 \, dy\, \frac{ds}{|s|^{n-2}}  \right)^\frac12
\end{split}\]
Therefore, the convergence of $T_1$ to 0 will be established if we can show that 
\[
\lim_{r\to 0}  \iint_{\W_{16a}(x,r)} |\nabla u(y,s)|^2 \, dy\, \frac{ds}{|s|^{n-2}} = 0 \qquad \text{ for a.e. } x\in \R^d,
\]
which is a consequence of the fact that for any ball $B \subset \R^d$,
\begin{equation} \label{lem8.11c} 
   \lim_{{r\to 0}} \int_{x\in B} \iint_{\W_{16a}(x,r)} |\nabla u(y,s)|^2 \, dy\, \frac{ds}{|s|^{n-2}} \, dx = 0.
\end{equation}
Since, for $r$ small enough,
$$\int_{x\in B} \iint_{\W_{16a}(x,r)} |\nabla u(y,s)|^2 \, dy\, \frac{ds}{|s|^{n-2}} \, dx \leq \int_{2B} \int_{|t| \leq r} |\nabla u(y,s)|^2 \, dy\, \frac{ds}{|s|^{n-d-2}} \, dx,$$
the convergence \eqref{lem8.11c} is an immediate biproduct of the fact that
$$\int_{2B} \int_{|t| \leq r} |\nabla u(y,s)|^2 \, dy\, \frac{ds}{|s|^{n-d-2}} \, dx \leq \int_{2B} |S_{a,2}(u)|^2 dx$$
is finite. Thanks to Lemma \ref{lemS<Nd} and an increasing sequence of compactly supported cut-off functions $\Psi_i \uparrow 1$, it is easy to obtain that
\begin{equation}
\|S_{a,2}(u)\|_{L^2(2B)} \leq C_B \|S_{a,2}(u)\|_q \lesssim \|\wt N_{a,2}(u)\|_q \leq  \|\wt N_{a,q}(u)\|_q < +\infty.
\end{equation}
The claim \eqref{claim811}, in the case $q\geq 2$, follows.

\medskip

Let us turn to the proof of the claim in the case $q\leq 2$. By Poincar\'e's inequality,
$$T_1 \lesssim r \left(\frac{1}{|\W_a(z,r)|} \iint_{W_a(z,r)} |\nabla u(y,s)|^q \, dy\, ds \right)^\frac1q \lesssim \left(\iint_{W_{2a}(x,r)} |\nabla u(y,s)|^q \, dy\, \frac{ds}{|s|^{n-2}} \right)^\frac1q.$$
Then, the use of H\"older's inequality with $\alpha = \frac2q \geq 1$ gives that
\[\begin{split}
T_1 & \lesssim \left( \iint_{W_{2a}(x,r)} |\nabla u|^2 |u|^{q-2} \, dy\, \frac{ds}{|s|^{n-2}} \right)^\frac12 \left(\iint_{W_{2a}(x,r)} |u|^{q} \, dy\, \frac{ds}{|s|^{n-2}} \right)^{\frac1q(1-\frac q2)} \\
& \quad \leq \left( \iint_{W_{2a}(x,r)} |\nabla u|^2 |u|^{q-2} \, dy\, \frac{ds}{|s|^{n-2}} \right)^\frac12 |\wt N_{2a,q}(u)|^{1-\frac q2}.
\end{split}\]
In order to prove the claim \eqref{claim811}, we will show that $T_1$ converges to 0 in $L^q(\R^d)$. We have by H\"older's inequality again,
\[\begin{split}
\int_{x\in \R^d} T_1 \, dx & \lesssim \int_{x\in \R^d} \left( \iint_{W_{2a}(x,r)} |\nabla u|^2 |u|^{q-2} \, dy\, \frac{ds}{|s|^{n-2}} \right)^\frac q2 |\wt N_{2a,q}(u)|^{q(1-\frac q2)} \, dx \\
& \lesssim \left\| \left(\iint_{W_{2a}(\cdot,r)} |\nabla u|^2 |u|^{q-2} \, dy\, \frac{ds}{|s|^{n-2}}\right)^\frac1q \right\|^{q/2}_{L^q} \left\| \wt N_{2a,q}(u) \right\|_{q}^{1-q/2} \\
& \lesssim \left(\int_{x\in \R^d} \int_{|t| \leq 2r} |\nabla u|^2 |u|^{q-2} dx \, \frac{dt}{|t|^{n-d-2}} \right)^\frac12  \left\| \wt N_{2a,q}(u) \right\|_{q}^{1-q/2}.
\end{split}\]
The right-hand term above converges to 0 as $r$ tends to 0. Indeed, due to Lemma \ref{lemMoser22},
$$\left\| \wt N_{2a,q}(u) \right\|_{q} \lesssim \left\| \wt N_{a,q}(u) \right\|_{q} < +\infty$$
and the fact that $\int_{x\in \R^d} \int_{|t| \leq 2r} |\nabla u|^2 |u|^{q-2} dx \, \frac{dt}{|t|^{n-d-2}}$ goes to 0 as $r\to 0$ is an easy consequence of the fact that, thanks to Lemma \ref{lemS<Nd}, one has $\|S_{a,q}(u)\|_q \lesssim \left\| \wt N_{a,q}(u) \right\|_{q} < +\infty$. The claim \eqref{claim811} when $q\leq 2$ follows.
\ep

\begin{lem} \label{lem8.1'}
Let $\cL$ be an elliptic operator satisfying $(\cH^1_\kappa)$ for some constant $\kappa\geq 0$. Let $a>0$, $q\in (q_0,q'_0)$ where $q_0$ is given by Proposition \ref{prop1.1}. Choose $k\geq 1$. Let $g\in L^q$ and $u\in W^{1,2}_{loc}(\Omega)$ be a weak solution to $\cL u = 0$ that satisfies both $\|\wt N_{a,q}(u)\|_q < \infty$ and
\begin{equation} \label{lem8.1'a} 
   \lim_{(z,r)\in \Gamma_a(x) \atop {r\to 0}} \frac{1}{|\W_a(z,r)|} \iint_{W_a(z,r)} u(y,s) \, dy\, ds  = g(x)
\end{equation}
for almost every $x\in\R^d$.
We have
\begin{equation} \label{lem8.1'b}
\limsup_{\epsilon \to 0} \int_{y\in \R^d} \int_{s\in \R^{n-d}} |u|^q \dr_r [\Psi_\epsilon^{k}] \frac{ds}{|s|^{n-d-1}}\, dy \lesssim \|g\|_q^q,
\end{equation}
where the constant depends on $a$, $q$, $n$, $\lambda_q(\cA)$, $\|\cA\|_\infty$, $\|b^{-1}\|_\infty + \|b\|_\infty$, $\kappa$ and $k$.
\end{lem}

\begin{rmk}
	\begin{enumerate}
		\item The existence of such a $u$ is guaranteed by Lemma \ref{lem7.5}.
		\item Here the solution $u$ is only assumed to lie in $W^{1,2}_{loc}(\Omega)$, so its trace may not be defined. Instead we use the assumption \eqref{lem8.1'a} to describe ``$u=g$ on the boundary''. 
	In particular, if $u\in W$ is a weak solution to $\cL u=0$ such that $\Tr u = g$, then $u$ satisfies the assumption \eqref{lem8.1'a}.
		\item This lemma is an analogue and a generalization of Lemma \ref{lem7.1'}, which only holds for energy solutions in $W$.
	\end{enumerate}
\end{rmk}

\bp Observe first that we have the uniform bound
\begin{equation} \label{unq:RFasmp}
\begin{split} 
\int_{y\in \R^d} \int_{s\in \R^{n-d}} |u|^q \dr_r [\Psi_\epsilon^{k}] \frac{ds}{|s|^{n-d-1}}\, dy & \lesssim \int_{y\in \R^d} \frac{1}{|\W_a(y,\epsilon)|} \iint_{W_a(y,\epsilon)} |u(x,t)|^q \, dx\, dt\, dy \\
& \lesssim \int_{y\in \R^d } |\wt N_{a,q}(u)(y)|^q \, dy < +\infty.
\end{split}\end{equation}
By \eqref{lem8.1'a} and Lemma \ref{lem8.11}, we have
\[ \frac{1}{|\W_a(y,\epsilon)|} \iint_{W_a(y,\epsilon)} |u(x,t)|^q \, dx\, dt \to |g(y)|^q \quad \text{ as } \epsilon \to 0, \]
for almost every $y\in \R^d$.
Therefore, the reverse Fatou's lemma and the pointwise domination \eqref{unq:RFasmp} imply that
\[\begin{split}
\limsup_{\epsilon \to 0} \int_{y\in \R^d} \int_{s\in \R^{n-d}} |u|^q \dr_r [\Psi_\epsilon^{k}] \frac{ds}{|s|^{n-d-1}}\, dy & \lesssim \limsup_{\epsilon \to 0}  \int_{y\in \R^d} \frac{1}{\W_a(y,\epsilon)} \iint_{W_a(y,\epsilon)} |u(x,t)|^q \, dx\, dt\, dy \\
& \lesssim \int_{y\in \R^d} |g(y)|^q dy.
\end{split}\]
The lemma follows.
\ep

\begin{lem} \label{lem8.2}
Let $\cL$ be an elliptic operator satisfying $(\cH^1_\kappa)$ for some constant $\kappa\geq 0$. Let $a>0$, $q\in (q_0,q'_0)$ where $q_0$ is given by Proposition \ref{prop1.1}. For any $g\in L^q$, suppose $u\in W^{1,2}_{loc}(\Omega)$ is a weak solution to $\cL u = 0$ that satisfies both $\|\wt N_{a,q}(u)\|_q < \infty$ and
\begin{equation} \label{lem8.2a} 
   \lim_{(z,r)\in \Gamma_a(x) \atop {r\to 0}} \frac{1}{|\W_a(z,r)|} \iint_{W_a(z,r)} u(y,s) \, dy\, ds  = g(x) \qquad \text{ for almost every } x\in \R^d.
\end{equation}
Then
\begin{equation} \label{lem8.2b}
\|S_{a,q}(u)\|_q^q \leq C\kappa \|\wt N_{a,q}(u)\|_q^q + C \|g\|_q^q,
\end{equation}
where the constant $C>0$ depends only on $a$, $q$, $n$, $\lambda_q(\cA)$, $\|\cA\|_\infty$, $\|b^{-1}\|_\infty + \|b\|_\infty$.
\end{lem}

\bp
Let $\chi_l$ be the  function defined in the beginning of Section \ref{SExistence}. By Lemma \ref{lemS<Nb},
\begin{equation}\label{unq:apf}
	\|S_{a,q}(u|\chi_l^k)\|_q \lesssim \|\wt N_{a,q}(u|\chi_l^{k-2}) \|_q \leq \|\wt N_{a,q}(u) \|_q < + \infty.
\end{equation}  
We only proved Lemma \ref{lemS<Ne} for energy solutions in $W$ such that the trace operator is defined. But by the a priori estimate \eqref{unq:apf} in a similar manner we can get
\begin{equation} \label{lem8.2c} \begin{split}
c \|S_{a,q}(u|\chi_{l}^{k})\|_q^q
& \leq  C \kappa \|\wt N_{q,a} (u|\chi_{l}^{k-2})\|_{q}^q  + C \|g\|_q^q \\
& \hspace{3cm} + \, \frac1q \, \iint_{\Omega} \left(\frac{|u|^q}{|t|} - \dr_r [|u|^q] \right) \dr_r[\chi_l^k] \, \frac{dt}{|t|^{n-d-2}}\, dx.
\end{split} 
\end{equation}
Note that we use Lemma \ref{lem8.1'} for the second term on the right. 
The last term can be bounded as follows
\[\begin{split}
& \left|  \frac1q \, \iint_{\Omega} \left(\frac{|u|^q}{|t|} - \dr_r [|u|^q] \right) \dr_r[\chi_l^k] \, \frac{dt}{|t|^{n-d-2}}\, dx \right| \\
& \qquad \lesssim \int_{x\in \R^d} \fint_{l\leq a|t| \leq 2l} |u|^q \, dt\, dx + \int_{x\in \R^d} \int_{l\leq a|t| \leq 2l} |u|^{q-1} |\nabla u| \, \frac{dt}{|t|^{n-d-1}}\, dx \\
& \qquad \lesssim \int_{x\in \R^d} \fint_{l\leq a|t| \leq 2l} |u|^q \, dt\, dx + \left( \int_{x\in \R^d} \fint_{l\leq a|t| \leq 2l} |u|^q \, dt\, dx \right)^\frac12 \left( \iint_{(x,t)\in \Omega} |\nabla u|^2 |u|^{q-2} \, \frac{dt}{|t|^{n-d-2}} \, dx \right)^\frac12 \\
& \qquad  \lesssim \int_{x\in \R^d \setminus B_{l/4}(0)} |\wt N_{a,q}(u)|^q dx + \left( \int_{x\in \R^d \setminus B_{l/4}(0)} |\wt N_{a,q}(u)|^q dx \right)^\frac12 \|S_{a,q}(u)\|_q^{q/2} \\
& \qquad \rightarrow 0 \quad \text{ as } l\to +\infty.
\end{split}\]
Hence the lemma follows by taking $l\to +\infty$ in \eqref{lem8.2c}.
\ep

\begin{lem} \label{lem8.3}
Let $\cL$ be an elliptic operator satisfying $(\cH^1_\kappa)$ for some constant $\kappa\geq 0$. Let $a>0$, $q\in (q_0,q'_0)$ where $q_0$ is given by Proposition \ref{prop1.1}. There exist two values $\kappa_0 >0$ and $C>0$, both depending only on $a$, $q$, $n$, $\lambda_q(\cA)$, $\|\cA\|_\infty$ and $\|b^{-1}\|_\infty + \|b\|_\infty$, such that if $\kappa \leq \kappa_0$, then for any $g\in L^q$ and any weak solution $u\in W^{1,2}_{loc}(\Omega)$ to $\cL u = 0$ that satisfy both $\|\wt N_{a,q}(u)\|_q < \infty$ and
\begin{equation} \label{lem8.3a} 
   \lim_{(z,r)\in \Gamma_a(x) \atop {r\to 0}} \frac{1}{|\W_a(z,r)|} \iint_{W_a(z,r)} u(y,s) \, dy\, ds  = g(x) \qquad \text{ for almost every } x\in \R^d,
\end{equation}
we have
\begin{equation} \label{lem8.3b}
\|\wt N_{a,q}(u)\|_q \leq  C \|g\|_q.
\end{equation}
\end{lem}

\bp Let $\kappa_0 \leq 1$ to be fixed later and $\kappa \leq \kappa_0$.
Lemmas \ref{lem8.2} and \ref{lem8.1} give that
$$\|\wt N_{a,q}(u)\|_q^q  \leq C\kappa \|\wt N_{a,q}(u)\|_q^q + C \|g\|_q^q,$$
where the constant $C$ doesn't depend on $\kappa$ anymore (since $\kappa \leq 1$). We choose $\kappa_0$ such that $C\kappa_0 \leq \frac12$ and the lemma follows.
\ep

\begin{lem} \label{lem8.4}
Let $\cL$ be an elliptic operator satisfying $(\cH^1_\kappa)$ for some constant $\kappa\geq 0$. Let $a>0$, $q\in (q_0,q'_0)$ where $q_0$ is given by Proposition \ref{prop1.1}. There exists $\kappa_0 >0$ depending only on $a$, $q$, $n$, $\lambda_q(\cA)$, $\|\cA\|_\infty$ and $\|b^{-1}\|_\infty + \|b\|_\infty$ such that if $\kappa \leq \kappa_0$, then for any couple of weak solution $u_1,u_2\in W^{1,2}_{loc}(\Omega)$ to $\cL u_1 = \cL u_2 = 0$ that satisfy 
$$\|\wt N_{a,q}(u_1)\|_q + \|\wt N_{a,q}(u_2)\|_q < \infty,$$
and
\begin{equation} \label{lem8.4a} 
   \lim_{(z,r)\in \Gamma_a(x) \atop {r\to 0}} \frac{1}{|\W_a(z,r)|} \iint_{W_a(z,r)} u_1(y,s) \, dy\, ds  = \lim_{(z,r)\in \Gamma_a(x) \atop {r\to 0}} \frac{1}{|\W_a(z,r)|} \iint_{W_a(z,r)} u_1(y,s) \, dy\, ds
\end{equation}
 for almost every $x\in \R^d$, we have
\begin{equation} \label{lem8.4b}
u_1 = u_2 \text{ almost everywhere in } \R^n.
\end{equation}
\end{lem}

\bp
Set the weak solution $v=u_1 - u_2 \in W^{1,2}_{loc}(\Omega)$. The function $v$ satisfies
$$\lim_{(z,r) \in \Gamma_a(x) \atop{r\to 0}} \frac{1}{|\W_a(z,r)|} \iint_{W_a(z,r)} |v(y,s)|^q \, dy = 0 \quad \text{ for almost every } x\in \R^d$$
and
$$\|\wt N_{a,q}(v)\|_q < +\infty.$$
So if $\kappa_0$ is chosen as in Lemma \ref{lem8.3}, we have $\|\wt N_{a,q}(v)\|_q = 0$ and so $v \equiv 0$. The lemma follows.
\ep
\begin{rmk}
	It is easy to see that the same conclusion holds if we assume, in place of \eqref{lem8.4a}, that
	\[ \lim_{(z,r)\in \Gamma_a(x) \atop {r\to 0}} \frac{1}{|\W_a(z,r)|} \iint_{W_a(z,r)} u_i(y,s) \, dy\, ds  = g(x), \qquad \text{ for almost every } x\in \R^d, \]
	with $i=1,2$. Moreover by Proposition \ref{lemMoser22}, $\|\wt N_{a,q}(u_i)\|_q <+\infty$ if and only if $\|\wt N_{a,2}(u_i)\|_q <+\infty$. Therefore we finish the proof of the uniqueness of solutions to Dirichlet problem with given boundary value $g$.
\end{rmk}


\begin{thebibliography}{AAA}

\bibitem[Ax]{Ax} A. Axelsson. {\em Non-unique solutions to boundary value problems for non-symmetric divergence form equations,} Trans. Amer. Math. Soc. {\bf 362} (2010), no. 2, 661--672.

\bibitem[AA]{AA} P. Auscher, A. Axelsson. {\em Weighted maximal regularity estimates and solvability of non-smooth elliptic systems I.} Invent. Math. {\bf 184} (2011), no. 1, 47--115. 

\bibitem[AHLMT]{Kato} P. Auscher, S. Hofmann, M. Lacey, A. McIntosh, P. Tchamitchian. {\em The solution of the Kato square root problem for second order elliptic operators on $\R^n$}. Ann. of Math. (2) {\bf 156} (2002), no. 2, 633--654. 

\bibitem[AMM]{AMM} P. Auscher, A. McIntosh, M. Mourgoglou. {\em On $L^2$ solvability of BVPs for elliptic systems.} J. Fourier Anal. Appl. {\bf 19} (2013), no. 3, 478-494. 


\bibitem[BM]{BM} A. Barton, S. Mayboroda. {\em Layer potentials and boundary-value problems for second order elliptic operators with data in Besov spaces.} Mem. Amer. Math. Soc. {\bf 243} (2016), no. 1149, v+110 pp.


\bibitem[CD]{CD} A. Carbonaro, O. Dragi\v cevi\'c. \em Convexity of power functions and bilinear embedding for divergence-form operators with complex coefficients. \em Preprint, arXiv:1611.00653.

\bibitem[CM]{CM} A. Cialdea, V. Maz'ya. \em Criterion for the $L^p$-dissipativity of second order diffrential operators with complex coefficients. \em J. Math. Pures Appl. {\bf 84} (2005), no. 9, 1067--1100.

  \bibitem[Da]{Da}
B.~E.~J. Dahlberg, \emph{Estimates of harmonic measure}.  Arch. Rational Mech.
  Anal. \textbf{65} (1977), no.~3, 275--288. 
  
  \bibitem[Da2]{Da2}
B.~E.~J. Dahlberg, \emph{ On the Poisson integral for Lipschitz and C1-domains.} Studia Math. 66 (1979), no. 1, 13--24.

\bibitem[DFM1]{DFMprelim} G. David, J.  Feneuil, S. Mayboroda. {\em Elliptic theory for sets with higher co-dimensional boundaries}.  Preprint, arXiv:1702.05503.

\bibitem[DFM2]{DFMAinfty} G. David, J.  Feneuil, S. Mayboroda. {\em Dahlberg's theorem in higher co-dimension}.  Preprint, arXiv:1704.00667.

\bibitem[DFM3]{DFM3} G. David, J. Feneuil, S. Mayboroda. \textit{A new elliptic measure on lower dimensional sets}. Preprint, arXiv:1807.07035.

\bibitem[DPP]{DPP} M. Dindo\v s, S.\,Petermichl, J.\,Pipher. 
{\em The {$L^p$} {D}irichlet problem for second order elliptic operators and a {$p$}-adapted square function}. J. Funct. Anal. {\bf 249} (2007), no. 2, 372--392. 

\bibitem[DPP2]{DPP2} M. Dindo\v s, S.\,Petermichl, J.\,Pipher, {\it BMO solvability and the $A^\infty$ condition for second order parabolic operators}. Preprint, arXiv:1510.05813.

\bibitem[DP]{DP} M. Dindo{\v s}, J. Pipher. \em Regularity theory for solutions to second order elliptic operators with complex coefficients and the {$L^p$} {D}irichlet problem. \em Preprint, arXiv:1612.01568.

\bibitem[FJK1]{FJK} E. Fabes, D. Jerison, C. Kenig.
{\em The {W}iener test for degenerate elliptic equations}.
Ann. Inst. Fourier (Grenoble), 32 (1982), no. 3, 151--182.

\bibitem[FKP]{FKP} R. Fefferman, C. Kenig, J. Pipher. {\em The theory of weights and the Dirichlet problem for elliptic equations}. Ann. of
Math. (2) {\bf 134} (1991), no. 1, 65--124.

\bibitem[FKS]{FKS}  E. Fabes, C. Kenig, R. Serapioni.
{\it The local regularity of solutions of degenerate elliptic equations.} 
Comm. Partial Differential Equations, 7 (1982), no. 1, 77--116. 


\bibitem[HKMP]{HKMP}
S. Hofmann, C. Kenig, S. Mayboroda, J. Pipher. \emph{Square function\slash non-tangential maximal function estimates
  and the {D}irichlet problem for non-symmetric elliptic operators}, J. Amer.
  Math. Soc. \textbf{28} (2015), no.~2, 483--529. 
  




\bibitem[HMM]{HMaMo} S. Hofmann; S. Mayboroda, M. Mourgoglou. {\em Layer potentials and boundary value problems for elliptic equations with complex $L^\infty$ coefficients satisfying the small Carleson measure norm condition.} Adv. Math. {\bf 270} (2015), 480--564. 

\bibitem[JK]{JK} D. Jerison and C. Kenig. {\em The Dirichlet problem in nonsmooth domains.} Ann. of Math. (2) {\bf 113}  (1981), no. 2, 367--382.

\bibitem[Ken]{KenigB} C. E. Kenig.
{\em Harmonic analysis techniques for second order elliptic boundary value problems.} 
CBMS Regional Conference Series in Mathematics, {\bf 83}. Amer. Math. Soc., Providence, RI, 1994.

\bibitem[KKPT]{KKPT} C. Kenig, H. Koch, J. Pipher, T. Toro. 
{\it A new approach to absolute continuity of elliptic measure, with applications to non-symmetric equations.} 
Adv. Math., {\bf 153} (2000), no. 2, 231--298.

\bibitem[KKiPT]{KKiPT} C. Kenig, B. Kirchheim, J. Pipher, T. Toro. 
{\it Square Functions and the $A^\infty$ Property of Elliptic Measures.} 
J. Geom. Anal., {\bf 26} (2016), no. 3, 2383--2410.

\bibitem[KP]{KePiDrift} C. Kenig, J. Pipher.  
{\it The Dirichlet problem for elliptic equations with drift terms.} 
Publ. Mat., {\bf 45} (2001), no. 1, 199--217. 

\bibitem[Ma]{Mattila} P. Mattila. {\em Geometry of sets and measures in Euclidean spaces. Fractals and rectifiability.} Cambridge Studies in Advanced Mathematics, {\bf 44}. Cambridge University Press, Cambridge, 1995. xii+343 pp.

\bibitem[M]{M} S. Mayboroda. {\em The connections between Dirichlet, regularity and Neumann problems for second order elliptic operators with complex bounded measurable coefficients.} Adv. Math. {\bf 225} (2010), no. 4, 1786--1819. 

\bibitem[MZ]{MZ} S. Mayboroda, Z. Zhao. 
{\em Square function estimates, {BMO} {D}irichlet problem, and absolute continuity of harmonic measure on lower-dimensional sets}. Preprint, arXiv:1802.09648.


\bibitem[St1]{Stein93} E. M. Stein.
{\em Harmonic analysis: real-variable methods, orthogonality, and oscillatory integrals}.
Princeton Mathematical Series, {\bf 43}. Princeton University Press, Princeton, N.J., 1993.

\bibitem[St2]{SteinSI} E. M. Stein.
{\em Singular integrals and differentiability properties of functions}.
Princeton Mathematical Series, {\bf 30}. Princeton University Press, Princeton, N.J., 1970.

\end{thebibliography}
\end{document}